\documentclass[
11pt, 
a4paper, 
oneside, 
headinclude,footinclude, 
]{article}

\usepackage[
nochapters, 
beramono, 
eulermath,
pdfspacing, 
dottedtoc 
]{classicthesis} 

\usepackage{arsclassica} 

\usepackage[T1]{fontenc} 

\usepackage[utf8]{inputenc} 

\usepackage{graphicx} 

\usepackage{enumitem} 

\usepackage{amsmath,amssymb,amsthm} 

\usepackage{varioref} 

\hypersetup{
colorlinks=true, breaklinks=true, bookmarks=true,bookmarksnumbered,
urlcolor=webbrown, linkcolor=RoyalBlue, citecolor=webgreen, 
pdftitle={}, 
pdfauthor={\textcopyright}, 
pdfsubject={}, 
pdfkeywords={}, 
pdfcreator={pdfLaTeX}, 
pdfproducer={LaTeX with hyperref and ClassicThesis} 
} 

\usepackage{tikz}

\usepackage{bm, stmaryrd}											
\usepackage{a4wide}
\usepackage{mathdots}
\usepackage[english]{babel}
\usepackage{xcolor}
\usepackage{float}

\newtheorem{thm}{Theorem}[section]
\newtheorem{lem}[thm]{Lemma}
\newtheorem{prop}[thm]{Proposition}
\newtheorem{conj}[thm]{Conjecture}
\newtheorem{coro}[thm]{Corollary}

\theoremstyle{definition}				
\newtheorem{rem}[thm]{Remark}

\numberwithin{equation}{section}

\newcommand{\R}{\mathbb{R}}

\newcommand{\N}{\mathbb{N}}
\newcommand{\Z}{\mathbb{Z}}

\newcommand{\K}{\mathbb{K}}

\newcommand{\F}{\mathbb{F}}

\newcommand{\rank}{\textbf{rank}}

\title{\normalfont\spacedallcaps{On the $P(t)$--adic Littlewood Conjecture in Characteristics $\ell\equiv 3\pmod{4}$}} 

\author{\spacedlowsmallcaps{Faustin ADICEAM\textsuperscript{1} \& \spacedlowsmallcaps{Dzmitry BADZIAHIN\textsuperscript{2}}} }

\date{} 

\begin{document}

\maketitle


\renewcommand{\sectionmark}[1]{\markright{\spacedlowsmallcaps{#1}}} 
\lehead{\mbox{\llap{\small\thepage\kern1em\color{halfgray} \vline}\color{halfgray}\hspace{0.5em}\rightmark\hfil}} 

\pagestyle{scrheadings} 

\begin{flushright}
\emph{To William Frederick Lunnon, \\ who introduced so many of us to\\ the universe of Number Walls.}\\
\end{flushright}

\begin{abstract}

\noindent Given a prime $p$, the \emph{$p$-adic Littlewood Conjecture} stands as a well--known arithmetic variant of the celebrated \emph{Littlewood Conjecture} in Diophantine Approximation. 
In the same way as the latter, it admits a natural function field analogue depen\-ding on the choice of an irreducible polynomial $P(t)$ with coefficients in a field $\K$. This analogue is referred to as the \emph{$P(t)$-adic Littlewood Conjecture} (\emph{$P(t)$--LC} for short). \\



\noindent 
$P(t)$--LC is proved to fail for \emph{any} choice of irreducible polynomial $P(t)$ over \emph{any} ground field $\K$  with characteristic $ \ell \equiv 3\pmod{4}$. The counterexample refuting it 
is shown to present a local arithmetic obstruction  emerging from the fact that -1 is not a quadratic residue  modulo a prime $ \ell\equiv 3\pmod{4}$.\\

\noindent  The theory developed elucidates and generalises all previous approaches towards refuting the conjecture. They were all based on the  computer--assisted method initiated by Adiceam, Nesharim and Lunnon (2021) which has been  able to establish that  $P(t)$--LC fails in some small cha\-ra\-cteristics (essentially up to 11). This computer--assisted  method is, however,  unable to provide a general statement  as it relies on  \emph{ad hoc} computer verifications which, provided they terminate, refute $P(t)$--LC in a given characteristic. This limitation is overcome  by exhibiting an arithmetic obstruction to the validity of $P(t)$--LC in infinitely many characteristics. \\

\noindent  The existence of arithmetic obstructions within the context of $P(t)$--LC  leaves the remaining case of odd charac\-teristics $ \ell\equiv 1\pmod{4}$ dependent on their determination. This is shown to hold in an effective and explicit way.
\end{abstract}

\tableofcontents

\let\thefootnote\relax\footnotetext{\textsuperscript{1} {Laboratoire d’analyse et de mathématiques appliquées (LAMA), Université Paris-Est Créteil, France,} \texttt{faustin.adiceam@u-pec.fr}}
\let\thefootnote\relax\footnotetext{\textsuperscript{2}The University of Sydney, Camperdown 2006, NSW, Australia, \texttt{dzmitry.badziahin@sydney.edu.au}.\\}

\section{Introduction}

\subsection{Littlewood--type Problems in Multiplicative Diophantine Approximation.} \noindent In 2004, de Mathan and Teulié~\cite{dMT} stated  the \emph{$p$--adic Littlewood Conjecture}. It claims that for any given prime number $p$, the equation
\begin{equation*}\label{pLC}\tag{$p$--LC}
 \inf_{q\in\Z\backslash\{0\}}\left| q\right| \cdot\left|q\right|_p \cdot \left| \left\langle q\alpha\right\rangle\right|\;=\;0
\end{equation*}
holds for all  reals $\alpha$. Here and throughout, $\left|\;\cdot\;\right|$ stands for the usual absolute value and $ \left|\left\langle\;\cdot\;\right\rangle\right|$ for the distance from a real parameter to the set of integers\textsuperscript{3}
\let\thefootnote\relax\footnotetext{\textsuperscript{3}This slightly unusual notation for the distance to the nearest integer function  is meant to achieve a notational analogy with the function field setup introduced below. It is needed in this introductory subsection only.}. Furthermore, $\left| q\right|_p$ denotes the $p$-adic absolute value defined for any integer $q\in\Z\backslash\{0\}$ as $\left|q\right|_p = p^{-\nu_p(q)}$, where $\nu_p(q)$ stands for the $p$-adic valuation of $q$, and where consequently $p^{\nu_p(q)}$ is the largest power of $p$ dividing $q$.\\

\noindent The statement of $p$--LC should be seen as an arithmetic variant of \emph{Littlewood's Conjecture} (1930's) which asserts that for any real numbers $\alpha, \beta\in\R$,
\begin{equation*}\label{LC}\tag{LC}
 \inf_{q\in\Z\backslash\{0\}}\left| q\right| \cdot\left|\left\langle q\alpha\right\rangle\right|\cdot \left|\left\langle q\beta\right\rangle\right|\;=\;0.
\end{equation*}
In fact,  $p$--LC likely stemmed from the assumption that the arithmetic function $q\mapsto \left|q\right|_p$ should behave, within the context of multiplicative Diophantine approximation, in the same way as the analytic map $q\mapsto  \left|\left\langle q\beta\right\rangle\right|$ intervening in the statement of LC.\\

\noindent Despite a wealth of partial contributions to $p$--LC and to LC (especially from a metric standpoint), the two conjectures remain open. A quite comprehensive summary of the state of the art regarding these two conjectures can be found in the survey~\cite{B} and in (the introduction of) the paper~\cite{BBEK}. For the purpose of the present work, it is worth mentioning  the 2006 result by Einsiedler, Katok and Lindenstrauss~\cite{EKL} who, through measure classification, established that the set of putative counterexamples to LC has Hausdorff dimension zero. In 2007, Einsiedler and Kleinbock~\cite{EK}  proved the same  conclusion for   the set of putative counterexamples to $p$--LC. \\

\noindent Both LC and $p$--LC admit analogues over functions fields motivated by a classical dictionary between the real and function field setups. According to it, the sets of integers,  rationals and reals  correspond to the sets of polynomials, rational functions and  Laurent series, respectively, over a given field $\K$. The analogue of the prime numbers (the irreducible elements over the ring of integers) is then the set of irreducible polynomials over $\K$. \\

\noindent  More formally, let $\K[t]$ be the ring of polynomials over $\K$ and let $\K(t)$ be the field of rational functions over it. The valuation over  $\K[t]$  defined by the degree extends to a valuation over $\K(t)$  so as to determine an absolute value,   referred to as a \emph{norm}, given for any $\Gamma(t)\in\K(t)$ by
\begin{equation}\label{laurser}
\left| \Gamma(t) \right|\;=\; 2^{\deg \Gamma (t)}.
\end{equation}
The completion of the field of rational functions with respect to this norm gives rise to the field of formal Laurent series in the variable $t^{-1}$, denoted by $\K\left(\left( t^{-1}\right)\right)$. This is  the set of formal expansions with only finitely many nonzero terms corresponding to positive exponents, namely expansions of the form $$\Gamma(t)\;=\; \sum_{k=-n_0}^{\infty}\gamma_k t^{-k}.$$ Here, $n_0\in\Z$ and $\left(\gamma_k\right)_{k\ge -n_0}$ is a sequence in $\K$ such that $\gamma_{-n_0}\neq 0$. This decomposition is unique, and the integer $-n_0$ determines the \emph{degree $\deg \Gamma(t)$ of the Laurent series $\Gamma (t)$}. Its norm is then given by
\begin{equation}\label{defnorm}
\left|\Gamma(t)\right|=2^{n_0}=2^{\deg \Gamma (t)},
\end{equation}
and its \emph{fractional part} is  the subsum in the expansion containing only negative exponents. It is denoted by
\begin{equation*}
\left\langle\Gamma(t)\right\rangle\;=\; \sum_{k=1}^{\infty}\gamma_k t^{-k}
\end{equation*}
and is obtained by subtracting from $\Gamma(t)$ its \emph{polynomial part} $\sum_{0\le k\le -n_0} \gamma_k t^{-k}$ (which is nontrivial precisely when $n_0\le 0$). From now on, a Laurent series such as in~\eqref{laurser} is identified with the sequence of its coefficients $\bm{\gamma}=\left(\gamma_k\right)_{k\ge -n_0}$.\\

\noindent With the above analogy, notations and terminology, Davenport and Lewis~\cite{DL} stated in 1963 the \emph{Littlewood Conjecture over Function Fields} as the claim that all $\Gamma(t), \Delta(t)\in \K\left(\left( t^{-1}\right)\right)$ meet the equation
\begin{equation*}\label{pLC}\tag{FFLC}
 \inf_{Q(t)\in\K[t]\backslash\{\bm{0}\}}\; \left| Q(t)\right| \cdot\left|\left\langle Q(t)\cdot \Gamma(t)\right\rangle\right|\cdot \left|\left\langle Q(t)\cdot\Delta(t)\right\rangle\right|\;=\;0.
\end{equation*}
Similarly, in the same paper as where they stated $p$--LC, De Mathan and Teulié~\cite{dMT} dealt with its function field analogue by asserting that, given an irreducible polynomial $P(t)\in\K[t]$, the \emph{$P(t)$--adic Littlewood Conjecture} should hold. This is saying that
\begin{equation*}\label{P(t)--LC}\tag{$P(t)$--LC}
 \inf_{Q(t)\in\K[t]\backslash\{\bm{0}\}}\; \left| Q(t)\right| \cdot\left|Q(t)\right|_{P(t)} \cdot \left| \left\langle Q(t)\cdot\Gamma(t)\right\rangle\right|\;=\;0
\end{equation*}
for all $\Gamma(t)\in \K\left(\left( t^{-1}\right)\right)$. In this relation, the $P(t)$--adic norm $\left|Q(t)\right|_{P(t)}$ of a polynomial  $Q(t)\in\K[t]$ is defined analogously to the real case~: it is the inverse of the largest power of $P(t)$ dividing $Q(t)$.\\

\noindent Davenport and Lewis~\cite{DL} on the one hand and de Mathan and Teulié~\cite{dMT}  on the other noticed when stating  FFLC and $P(t)$--LC, respectively,  that both conjectures fail when the ground field $\K$ is infinite. The situation when $\K$ is finite (which should be seen as the actual function field analogue of the real case) is not as well understood. For instance, compared with the real case, the measure classification results undertaken by Einsiedler, Lindenstrauss, and Mohammadi~\cite{ELM}  in positive characteristics
are not conclusive to determine whether the  respective sets of counterexamples to FFLC and $P(t)$--LC should have Hausdorff dimension zero. As far as FFLC is concerned, this essentially sums up the state of the knowledge  when it comes to establishing or refuting it in full generality. Particular  examples of pairs of algebraic power series satisfying the conjecture nontrivially are, however, also constructed in~\cite{AB}.

\subsection{The State of the Art in $P(t)$--LC.} The situation in the case of $P(t)$--LC is rather different 
thanks to recent advances. Indeed, Adiceam, Nesharim and Lunnon~\cite{ANL} established in 2021 that when $P(t)=t$, the $t$-adic Littlewood Conjecture ($t$--LC) fails  over fields with characteristic $ \ell=3$. Their proof amounts to showing that an infinite array recording all Hankel determinants formed from the paperfolding sequence (seen as a sequence over $\F_3$), the so-called \emph{Number Wall} of this sequence, has bounded zero regions. This is achieved showing that it  enjoys a two-dimensional automatic structure, which is established by computer inspection. The formal definition of a Number Wall is given in Section~\ref{sec2} below. Also, for a precise definition of the property of automaticity, the reader is referred to~\cite{AS}. The paperfolding sequence is introduced therein in~\cite[p.155, Example~5.1.6]{AS} and some of  its many properties are further surveyed in~\cite{MFvdP}. For the purpose of the present work, the charac\-terising property of discrete self-similarity enjoyed by this sequence is pivotal and is thus taken as its definition. \\

\noindent More precisely, given a sequence $\bm{v}=\left(v_n\right)_{n\ge 0}$ over  $\Z$ and a pair $\bm{\sigma}=\left(a,b\right)\in\Z^2$, let $\bm{v}\!\left[\bm{\sigma}\right]$ be the sequence obtained by inserting alternatively the values $a$ and $b$ before every element of the sequence $\bm{v}$; that is,
\begin{equation*}
\bm{v}\!\left[\bm{\sigma}\right]=\left( a,\; v_0, \; b, \; v_1,\; a, \; v_2,\; b,\; v_3,\; a, \; v_4, \dots\right)\;\quad \textrm{whenever}\;\quad  \bm{v}=\left(v_0, v_1,  v_2,  v_3,  v_4,  \dots\right).
\end{equation*}
\noindent The  \emph{paperfolding sequence} over the two-letter alphabet $\bm{\sigma}$ is then the sequence  $\bm{\lambda}=\left(\lambda_n\right)_{n\ge 0}$ taking values in $\{a,b\}$ characterised by the identity
\begin{equation}\label{pfdalphbab}
\bm{\lambda}\;=\;   \bm{\lambda}\!\left[\bm{\sigma}\right].
\end{equation}
Thus, $\bm{\lambda}=\left(a, \; a,\;  b,\; a,\; a,\; b,\; b,\; a,\; a,\; a,\; b,\; b,\; a,\; b,\; b,\; ...\right)$. \\

\noindent The counterexample to $t$--LC over fields of characteristic 3 found by Adiceam, Nesharim and Lunnon~\cite{ANL} corresponds to the Laurent series whose coefficients are the paperfolding sequence over the alphabet $\{0,1\}$. The authors conjecture that this Laurent series should remain a counterexample in all characteristics $ \ell\equiv 3\pmod{4}$. This has been verified following the same computer--assisted method by Garrett and Robertson~\cite{GR} in characteristics $ \ell=7 \textrm{ and }  \ell=11$. \\

\noindent Robertson also proves in~\cite{R} through an elementary and elegant  argument that a counterexample to $t$--LC induces a counterexample to $P(t)$--LC for any choice of irreducible polynomial $P(t)$ over a given field $\K$. This feature, specific to the function field setup, reduces the refutation of $P(t)$--LC to that of $t$--LC.

\subsection{The Main Theorem.}\label{mainres} The purpose of the present work is to confirm the aforementioned conjecture made by Adiceam, Nesharim and Lunnon~\cite{ANL} in the following strong sense~:

\begin{thm}[The $P(t)$--adic Littlewood Conjecture fails over fields of characteristics $ \ell\equiv 3\pmod{4}$]\label{mainthm}
Let $\K$ be a field of characteristic $ \ell\equiv 3\pmod{4}$ and let  $P(t)$ be an irreducible polynomial with coefficients in $\K$. Then, $P(t)$--LC fails over $\K$. More precisely, consider  the Laurent series
\begin{equation}\label{laurlambda}
\Lambda(t)\;=\; \sum_{n=1}^{\infty}\lambda_nt^{-n}
\end{equation}
whose coefficients are given by the paperfolding sequence  $\bm{\lambda}=\left(\lambda_n\right)_{n\ge 0}$   over the alphabet $(1,-1)$.
Then, $\Lambda(P(t))$ is a counterexample to $P(t)$--LC over $\K$ and, furthermore,
\begin{equation}\label{inflplitpfd}
\inf_{Q(t)\in\K[t]\backslash\{\bm{0}\}}\; \left| Q(t)\right| \cdot\left|Q(t)\right|_{P(t)} \cdot \left| \left\langle Q(t)\cdot\Lambda(t)\right\rangle\right|\;=\; 2^{-4\cdot\deg P(t)}.
\end{equation}
\end{thm}

\noindent The choice of the alphabet $(1,-1)$
(as opposed to the choice (0,1)
made by Adiceam, Nesharim and Lunnon~\cite{ANL}) is essentially for convenience~: it considerably simplifies the calculations in the proof. \\

\noindent To be more specific, following classical reductions, the proof boils down to working over a given finite field $\F_\ell$ with a prime order $\ell\equiv 3\pmod{4}$. It then comes down  to showing that a doubly indexed sequence of Hankel determinants formed from the paperfolding sequence, say $\left(H_{\bm{\lambda}}(m,n)\right)_{m, n\ge 0}$, cannot vanish along arbitrary large consecutive values of the index set  (in a sense made precise in Section~\ref{sec2}). The approach undertaken here to establish this property is in nature very different from the previous computer--assisted method developed by Adiceam, Nesharim and Lunnon~\cite{ANL}.\\

\noindent Indeed, the  method developed in~\cite{ANL} starts with recording all values of the Hankel determinants in the Number Wall of the paperfolding sequence over $\F_\ell$ (recall that this is an infinite two-dimensional array). It breaks down into  three main steps which can be sketched as follows~: (1) build a large enough portion of the Number Wall to detect by computer inspection (essentially by brute force) the existence of an underlying aperiodic structure; namely, a tiling given by some fixed rule prescribing how finitely many tiles should be combined ---- the algorithm fails if the computer is unable to find such a tiling; (2) consider the infinite matrix generated by the tiling rule and show that it is indeed a Number Wall  by checking finitely many known necessary and sufficient conditions on the finitely many tiles; (3) finally, check that within the (finitely many) generating tiles, the vanishing of the values of the Hankel determinants can only occur in the interior. The conclusion is that the Number Wall cannot contain arbitrarily large regions of vanishing entries. \\

\noindent This sketch makes it clear that the method does not provide a theoretical reason explaining the failure of $t$--LC. It should also be noted that the assistance of a computer in this task is fully necessary~: already over $\F_3$, the number of tiles needed in~\cite{ANL} is 2353. This number is significantly reduced by Garett and Robertson to 390 in~\cite{GR}. In fact,  a theoretical formalisation of this strategy is likely to rely on the proof of  the following conjecture, which is also stated in~\cite[\S 5]{GR}~:

\begin{conj}\label{NWconj}
The Number Wall of an automatic sequence over a given finite field is a (two-dimensional) automatic sequence.
\end{conj}

\noindent The approach undertaken here  to establish Theorem~\ref{mainthm} is of a different nature~: it  relies on the determination of a \emph{local arithmetic obstruction}  for the Laurent series~\eqref{laurlambda} to satisfy $t$--LC over prime fields of order $ \ell\equiv 3\pmod{4}$. \\

\noindent The nonvanishing property is indeed proven by establishing a set of recurrence relations between the values of the sequence $\left(H_{\bm{\lambda}}(m,n)\right)_{m, n\ge 0}$ (see Corollary~\ref{labelhenkelbis} in Section~\ref{sec3}). The key observation is that, given integers $m,n\ge 1$, one of these relations expresses, up to a sign factor, the determinant $H_{\bm{\lambda}}(2m,2n)$ as a sum of two squares, one of them being $H_{\bm{\lambda}}(m,n)$. Under the assumption that $ \ell\equiv 3\pmod{4}$, the fact that $-1$ is not a quadratic residue modulo $ \ell$ then forces the relation \mbox{$H_{\bm{\lambda}}(2m,2n)\equiv 0\pmod{ \ell}$} to imply that $H_{\bm{\lambda}}(m,n)\equiv 0 \pmod{ \ell}$. This allows for the setting up of an elaborate inductive argument showing that the determinant sequence can only  vanish  over $\F_ \ell$ along a string of at most three consecutive indices in specific configurations (see Section~\ref{sec4} for details). 
The main theoretical tool to achieve the proof of this claim is the theory of continued fractions in positive characteristics. This requires the specification  of a dictionary between  properties of Number Walls and properties of best approximants so as to rely on the better adapted theory at each step of the proof.\\



\noindent In the language of Number Walls, the above--described vanishing property translates into saying that, when $ \ell\equiv 3 \pmod{4}$,  the Number Wall of the paperfolding sequence over $\F_\ell$ does not contain square regions with zero entries which have  side-lengths   strictly bigger than 3. This is illustrated with pictures when $\ell\in\{3,7,11,19\}$ in the Appendix. Furthermore, the way the square regions are generated is explicitly determined in Section~\ref{sec4} in terms of divisibility properties of denominators of convergents of sui\-table Laurent series.
In this respect, the proof of  Theorem~\ref{mainthm}  can be seen as a first step towards elucidating the automatic structure of the Number Wall of the paperfolding sequence over finite fields, and thus as a first step towards Conjecture~\ref{NWconj}.

\subsection{On $P(t)$--LC in the Remaining Characteristics $ \ell\equiv1\pmod{4}$, and $ \ell=2$.}

The existence of a local arithmetic obstruction to the validity of $P(t)$--LC as unveiled by  the proof of Theorem~\ref{mainthm} is not expected to be specific to characteristics $ \ell\equiv 3\pmod{4}$. The situation when $ \ell\equiv 1\pmod{4}$ is, however, more intricate as there is no evidence whatsoever that the paperfolding Laurent series~\eqref{laurlambda} should fail $P(t)$--LC in this case. In fact, the case $\ell\equiv 3$ (mod~4) should
rather be seen as the first level in a tower of counterexamples exhausting all odd  characteristics. \\

\noindent To be precise, following the terminology introduced by Garrett and Robertson~\cite{GR} and given an integer $s\ge 1$,   the construction of the paperfolding sequence  
 is generalised so as to define the \emph{$s^{\textrm{th}}$--level paperfolding sequence} in the following way~: it is the integer sequence $\bm{u}^{(s)}=\left(u^{(s)}_n\right)_{n\ge 0}$ obtained by in\-ser\-ting periodically the successive values of  the finite set of integers $1, -1, 2, -2, 3, -3, \dots, s, -s$ before each element of the sequence itself. Thus, when $s=1$, this reduces to  the sequence $\bm{\lambda}[\bm{\sigma}]$ introduced in~\eqref{pfdalphbab}. When  $s=2$  for instance, this gives
\begin{equation*}
\left(u_0^{(2)}, u_1^{(2)}, u_2^{(2)}, u_3^{(2)}, u_4^{(2)},  u_5^{(2)}, u_6^{(2)}, u_7^{(2)}, \dots\right)\;=\; \left(1, u_0^{(2)}, -1, u_1^{(2)} , 2, u_2^{(2)} , -2, u_3^{(2)} , 1, \dots\right).
\end{equation*}
\noindent Garret and Robertson\textsuperscript{4}\let\thefootnote\relax\footnotetext{\textsuperscript{4}To be more precise, Garret and Robertson's definition of the $s^{th}$--level paperfolding sequence amounts to inserting periodically the successive values of the finite set of integers $0, 1, 2, \dots, 2s-1$ rather than the set $1, -1, 2, -2, 3, -3, \dots, s, -s$ as above. The choice adopted here, in compliance with Theorem~\ref{mainthm}, seems more natural in view of the proof of this theorem (that being said, it is not expected that the conclusion should differ in either case). Also, the  statement of Conjecture~\ref{mainconj} corrects what appears to be a typo in~\cite[Conjecture~1.8]{GR}~: the exponent of the power of 2 on the right--hand side of~\eqref{inflplitpfdbis} is stated therein as $s$, which would not include the case $s=1$ already established over $\F_3$ in~\cite{ANL} and over $\F_7$ and $\F_{11}$ in~\cite{GR}.}~\cite{GR} posed the following conjecture.

\begin{conj}[Garret \& Robertson]\label{mainconj}
Let $ \ell$ be a prime number and let \mbox{$s=\nu_2( \ell-1)$}, where $\nu_2(k)$ denotes the $2$-adic valuation of an integer $k\ge 1$. Then, given a field $\K$ with chara\-cteristic $ \ell$ and given an irreducible polynomial $P(t)\in\K[t]$, the above--defined $s^\textrm{th}$--level paper\-fol\-ding sequence fails $P(t)$--LC. Furthermore, denoting by $\Lambda^{(s)}(t)\in\K\left(\left(t^{-1}\right)\right)$ the Laurent series it determines, it holds that
\begin{equation}\label{inflplitpfdbis}
\inf_{Q(t)\in\K[t]\backslash\{\bm{0}\}}\; \left| Q(t)\right| \cdot\left|Q(t)\right|_{P(t)} \cdot \left| \left\langle Q(t)\cdot\Lambda^{(s)}(t)\right\rangle\right|\;=\; 2^{-2^{s+1}\cdot\deg P(t)}.
\end{equation}
\end{conj}

\noindent Applying once again the computer--assisted method developed by Adiceam, Nesharim and Lunnon~\cite{ANL}, Garrett and Robertson~\cite{GR} confirmed the above conjecture over  the finite field $\F_5$. This case is not covered by Theorem~\ref{mainthm}. However, as mentioned above, the existence of a local arithmetic obstruction underlying its proof is expected to still hold  for any value of the integer $s=\nu_2( \ell-1)$. The natural generalisation of  the observation holding in the $s=1$ case, namely that $-1$ is not a quadratic residue over $\F_\ell$, is the fact that the integer $s$ is the largest value of $k\ge 1$ for which there exists a $\left(2^k\right)^{\textrm{th}}$--root of unity over $\F_\ell$. Unfortunately, this obstruction alone is insufficient to extend the proof of Theorem~\ref{mainthm} to valuations $s\ge 2$. \\

\noindent The identification of the general form of the arithmetic obstruction when $s\ge 2$ might rely on class field theory following the theory developed in~\cite{C}. Establishing it seems a theoretical jump at the same scale as the one unveiled in the present work where, from the case of \emph{finitely} many odd characteristics dealt with the help of computer verifications, $P(t)$--LC is shown to fail in \emph{infinitely} many odd characteristics, which leads one to a clear picture of the way the conjecture should then fail in \emph{all} odd characteristics. \\

\noindent As for the case of the even characteristic, it remains mysterious.

\paragraph{Organisation of the paper.} Section~\ref{sec2} is devoted to the reduction of Theorem~\ref{mainthm} to a standard statement about the vanishing of Hankel determinants of the paperfol\-ding sequence over a finite field $\F_\ell$ such that  $ \ell\equiv 3\pmod{4}$. This property is rephrased in terms of the boundedness of the sizes of square regions (referred to as windows) in the Number Wall of this sequence. The following Section~\ref{sec3} is concerned with establi\-shing a set  of  recurrence relations enjoyed by the Hankel determinants. 
The nature of the local arithmetic obstruction for the Laurent series of the paperfolding sequence $\Lambda(t)$ to meet $P(t)$--LC  is showcased therein (this is Corollary~\ref{keyobstr}). Section~\ref{sec4} relies on this obstruction to determine explicitly the process of generation of the windows in the Number Wall under consideration (see Proposition~\ref{genwindow}). The boundedness of their sizes is then shown to derive from a property of divisibility of the denominators of convergents associated to the paperfolding Laurent series in Corollary~\ref{arithCNS}. The paper concludes in Section~\ref{secverif} with checking that this divisibility property is indeed verified.


\paragraph{Notations.}  Most of the notations intervening in the proof are introduced when needed as they are mostly standard. It is enough to mention here that a vector $\bm{x}\in\R^k$ is seen as a $1\times k$ matrix (i.e.~as a row vector). Its transpose is denoted by $\bm{x}^T
$. When it depends on a parameter, say an integer $n$, the vector is more conveniently denoted by $\bm{x}(n)$ and its transpose by $\bm{x}^T(n)$. \\

\noindent  To make it easier to refer to the numerous quantities involved in the calculations, the following convention is adopted~: a notation such as $x:=y$, where $x$ is introduced for the first time, means that it is defined by its identity with $y$.\\

\noindent  Finally, bullet points are employed in lengthy proofs to break them into subparts that can be followed successively in blocks.

\paragraph{Acknowledgement.}  This paper results from the visit of the first-named author to the second one in May 2025 as a laureate of the International Visitor Program funded by the Sydney Mathematical Research Institute (SMRI). The authors would like to thank the SMRI for excellent working conditions which made it possible to carry out the present work in a very efficient way. The authors  also   thank Steven Robertson for technical assistance.

\section{Reduction of the Problem}\label{sec2}

Throughout this section, fix a Laurent series $\Gamma(t)$ over a field $\K$. Given an irreducible polynomial $P(t)$ over $\K$, the statement of $P(t)$--LC makes it clear that it may be assumed, without loss of generality, that the polynomial part of $\Gamma(t)$ vanishes, i.e.~that $\Gamma(t)$  is of the form $$\Gamma(t)\;=\;\sum_{k=1}^{\infty}\gamma_k t^{-k}.$$ 


\noindent To avoid delicate boundary considerations to initialise recurrences, it is convenient to first extend the sequence of coefficients $\left(\gamma_k\right)_{k\ge 1}$ to a  doubly infinite sequence
\begin{equation}\label{defdoublyextendedeq}
\bm{\gamma}=\left(\gamma\right)_{k\in\Z}.
\end{equation}
Then, associate to it the  family of Laurent series
\begin{equation}\label{doublylaurent}
\Gamma_m(t)\;=\;\sum_{k=1}^\infty \gamma_{k+m-1}\cdot t^{-k}
\end{equation}
parametrised by the integer $m\in\Z$. The particular choice of the doubly infinite extension will be made precise later in the case of the paperfolding sequence. For the time being, it is enough to notice that, when $m=1$, the Laurent series $\Gamma_1(t)$ reduces to $\Gamma(t)$.\\

\noindent Given integers $m\in\Z$, and $p\ge n\ge 1$, define the Hankel matrix of size $n\times p$

\begin{equation}\label{defhankelmnp}
\mathcal{H}_{\bm{\gamma}}(m,n,p)=
\begin{pmatrix}
\gamma_m  \!& \!\cdots  \!& \! \gamma_k  \!& \! \cdots  \!& \! \gamma_{m+n-1}  \!& \! \cdots  \!& \! \gamma_r  \!&  \!\cdots  \!& \! \gamma_{m+p-1}\\
\vdots  \!& \! \iddots  \!& \!  \!& \! \iddots  \!& \!  \!& \! \iddots  \!& \!  \!& \! \iddots  \!& \! \vdots\\
\gamma_k \!& \! & \! \gamma_{m+n-1}  \!& \!  \!& \! \gamma_{r}  \!& \!  \!& \! \gamma_{m+p-1}  \!& \!  \!& \! \gamma_s\\
\vdots  \!& \! \iddots  \!& \!  \!& \! \iddots \!& \!  \!& \!\iddots  \!& \!  \!& \! \iddots  \!& \! \vdots\\
\gamma_{m+n-1}  \!& \! \cdots  \!& \! \gamma_r \!& \!\cdots  \!& \! \gamma_{m+p-1} \!& \!\cdots \!& \!\gamma_s \!& \!\cdots  \!& \!\gamma_{m+n+p-2}
\end{pmatrix}.
\end{equation}

\noindent When $n=p$,  the corresponding determinant is denoted by $H_{\bm{\gamma}}(m,n)=\det \left(\mathcal{H}_{\bm{\gamma}}(m,n, n)\right)$. In other words,

\begin{equation}\label{defhankelmn}
H_{\bm{\gamma}}(m,n)\;=\;
\begin{vmatrix}
\gamma_m &\cdots & \gamma_k& \cdots & \gamma_{m+n-1}\\
\vdots & \iddots && \iddots &\vdots\\
\gamma_k& & \gamma_{m+n-1} &&\gamma_{r}\\
\vdots & \iddots &&\iddots&\vdots\\
\gamma_{m+n-1} & \cdots &\gamma_r&\cdots & \gamma_{m+2n-2}
\end{vmatrix}.
\end{equation}

\noindent The  $n$--dimensional vector
\begin{equation}\label{defespivec}
 \bm{w}_{n}(t)=\left(1, t, t^2, \dots, t^{n-1} \right)\in\R^n.
\end{equation}
will be ubiquitous in the various statements and proofs.

\subsection{Vanishing of Hankel Determinants and $P(t)$--LC}

As mentioned in the introduction, the refutation of $P(t)$--LC can be reduced without loss of generality to the case where $P(t)=t$. The precise statement conveying this observation is due to Robertson~\cite[Theorem~1.0.3]{R}~:

\begin{thm}[Robertson]\label{red1}
With the above notations, assume that $\Gamma(t)\in\K\left(\left(t^{-1}\right)\right)$ fails $t$--LC in the sense that there exists an integer $k\ge 0$ such that
\begin{equation*}\label{tLCred}
\inf_{Q(t)\in\K[t]\backslash\{\bm{0}\}}\; \left| Q(t)\right| \cdot\left|Q(t)\right|_{t} \cdot \left| \left\langle Q(t)\cdot\Gamma(t)\right\rangle\right|\;=\; 2^{-(k+1)}.
\end{equation*}
Then, the Laurent series $\Gamma \left(P(t)\right)\in \K\left(\left(P(t)^{-1}\right)\right)\subset \K\left(\left(t^{-1}\right)\right)$ fails $P(t)$--LC; indeed, the following equation then holds~:
\begin{equation*}\label{tLCred}
\inf_{Q(t)\in\K[t]\backslash\{\bm{0}\}}\; \left| Q(t)\right| \cdot\left|Q(t)\right|_{P(t)} \cdot \left| \left\langle Q(t)\cdot\Gamma\left(P(t)\right)\right\rangle\right|\;=\; 2^{-(k+1)\cdot \deg P(t)}.
\end{equation*}
\end{thm}

\noindent The study of $t$--LC then reduces to the determination of the vanishing of Hankel determinants formed from the sequence of coefficients of a given Laurent series.  The following was essentially already observed by de Mathan and Teulié~\cite[\S 4]{dMT} when stating $P(t)$--LC, and is  formalised in~\cite[Theorem~2.2]{ANL}~:

\begin{thm}\label{red2}
Let $k\ge 0$ be a integer and let  $\Gamma(t)\in\K\left(\left(t^{-1}\right)\right)$ be a Laurent series as above. The following are equivalent~:
\begin{itemize}
\item[\textbf{(i)}] it holds that
\begin{equation}\label{tLCred1}
\inf_{Q(t)\in\K[t]\backslash\{\bm{0}\}}\; \left| Q(t)\right| \cdot\left|Q(t)\right|_{t} \cdot \left| \left\langle Q(t)\cdot\Gamma(t)\right\rangle\right|\;=\; 2^{-(k+1)}
\end{equation}
(in particular, $\Gamma(t)$ fails $t$--LC);

\item[\textbf{(ii)}] the integer $k$ is the maximal one for which there exist integers $m,n\ge 1$  such that all terms in the sequence of $k$ determinants  $$H_{\bm{\gamma}}(m, n), \qquad H_{\bm{\gamma}}(m, n+1)\;,  \qquad \dots \qquad,  \qquad H_{\bm{\gamma}}(m, n+k-1)$$  vanish.  (When  $k=0$, this condition is to be interpreted as the nonvanishing of all determinants $H_{\bm{\gamma}}(m, n)$ when $m,n\ge 1$).
\end{itemize}
\noindent If any of these conditions holds, then the infimum in~\eqref{tLCred1} is attained at a polynomial whose coefficients can be chosen in the minimal subfield of $\K$ containing the coefficients of the Laurent series $\Gamma(t)$.
\end{thm}

\begin{proof} Only the last claim is not explicitly established\textsuperscript{5}\let\thefootnote\relax\footnotetext{\textsuperscript{5} In all previous works --- see~\cite{ANL, GR, R}, the invariance of the infimum by field extension was said to follow just from the invariance of the norm of a Laurent series by field extension. However, something more is needed to reach the sought conclusion as the set of polynomials in the infimum is then taken over a larger set.} in~\cite[Theorem~2.2]{ANL}. It follows from two observations. \\

\noindent Firstly, and plainly, the norm of a Laurent series only depends on the smallest field which its coefficients belong to  (see relation~\eqref{defnorm}). \\

\noindent Secondly,  the observations made in~\cite[Remarks~2.4 \& 2.5]{ANL} show that the infimum in~\eqref{tLCred1} can be taken over the set of denominators of the convergents of the family of Laurent series  $\Gamma_m(t)$ introduced in~\eqref{doublylaurent}. From the Euclidean algorithm, such polynomials are defined over the same field as the coefficients of $\Gamma_m(t)$. Alternatively, this second observation can also be taken as this one~:  it is shown in~\cite[Equation~(2.5)]{ANL} that the coefficients of a (nonzero) polynomial achieving the infimum~\eqref{tLCred1} coincide with any nonzero element lying in the kernel of a matrix formed from the coefficients of $\Gamma(t)$. When nontrivial, it is a standard fact in Linear Algebra --- see~\cite[Theorem~51, p.159]{McC} --- that such a kernel element can be chosen in the same field as the one containing the coefficients of the matrix.  This completes the proof of Theorem~\ref{red2}.
\end{proof}

\noindent The doubly indexed sequence $\left(H_{\bm{\gamma}}(m,n)\right)_{n\ge 1, m\in\Z}$ is from now on recorded in the form of an infinite array referred to as the \emph{Number Wall of the sequence $\bm{\gamma}$} (\textsuperscript{6})\let\thefootnote\relax\footnotetext{\textsuperscript{6}This definition slightly differs from the one used in~\cite[\S 3]{ANL} and~\cite{L} inasmuch as the $(j,k)^{\textrm{th}}$ entry of a Number Wall in this paper corresponds to $(-1)^{k(k-1)/2}$ times the $(j-1, k-1)^{\textrm{st}}$ entry  of a Number Wall as defined in these references. The sign factor is useful to simplify the expressions of the algebraic  relations structuring the entries of a Number Wall (\emph{ibid.}). Since the goal here is to detect zero entries and that it is achieved without relying on these algebraic relations, this additional sign factor is dropped.\\}. The entries of this array, represented with matrix conventions\textsuperscript{7}\let\thefootnote\relax\footnotetext{\textsuperscript{7}This means that rows are successively ordered going down and columns from left to right.},  are labelled with indices $(j,k)$ where $k\in\Z$ and $j\ge 1$. The  $j^{\textrm{th}}$ row records the values of all determinants $H_{\bm{\gamma}}(m,n)$  of fixed size $n=j$ and the $k^{\textrm{th}}$ column  records the values of all determinants $H_{\bm{\gamma}}(m,n)$  along the set of diagonal  indices $n+m=k$. \\

\noindent To be more specific, consider
the indexation map
\begin{equation}\label{indexbij}
\iota\; :\; (m,n)\in\Z\times \N_{\ge 1} \;\;\mapsto\;\; (j,k)\;=\;\left(n, m+n\right)\in \N_{\ge 0}\times \Z
\end{equation}
with inverse
\begin{equation*}\label{indexbijinv}
\iota^{-1}\; :\; (j,k)\in  \N_{\ge 0}\times \Z \; \;\mapsto\;\; (m,n)\;=\; \left(k-j, j\right)\in \Z\times \N_{\ge 1}.
\end{equation*}
\noindent The Number Wall $NW_{\bm{\gamma}}$ is then a two dimensional array
$\left(NW_{\bm{\gamma}}(j,k)\right)_{k\in\Z, j\ge 1}$. Its  $(j,k)^{\textrm{th}}$ entry is determined by the integers $m\in\Z$ and $n\ge 1$ such that $(j,k)=\iota(m,n)$ in the sense that
\begin{equation*}\label{NWhankel}
NW_{\bm{\gamma}}(j,k)\; = \; H_{\bm{\gamma}}\left(\iota^{-1}(j,k)\right) \; = \; H_{\bm{\gamma}}\left(k-j, j\right)
\end{equation*}
and that
\begin{equation*}\label{NWhankelbis}
H_{\bm{\gamma}}\left(m, n\right) \; = \; NW_{\bm{\gamma}}\left(\iota(m,n)\right)\; = \; NW_{\bm{\gamma}}\left(n, m+n\right).
\end{equation*}
\noindent  In particular,  in view of the definition of the Hankel determinants in~\eqref{defhankelmn}, the first row of the Number Wall just records the doubly infinite sequence $\bm{\gamma}$ itself.  
Figure~\ref{fig1} below displays a portion of the Number Wall of a generic sequence $\bm{\gamma}$ and the corresponding indexation. The Appendix contains some illustrations of portions of the Number Wall of the paperfolding sequence over some finite fields.\\

\begin{figure}[h!]

\begin{center}
\scalebox{.7}{
\begin{tikzpicture}[ultra thick]

\draw [->]  (9,0) to (17,0);
\node  (k) at (18,0) {\Large{$k$}};

\draw [->]  (0,-9) to (0,-17);
\node  (j) at (0, -18) {\Large{$j$}};

\draw [-, dashed]  (4,-2) to (4,-4);
\draw [-, dashed]  (8,-2) to (8,-4);
\draw [-, dashed]  (12,-2) to (12,-4);
\draw [-, dashed]  (16,-2) to (16,-4);

\draw [-, dashed]  (4,-16) to (4,-18);
\draw [-, dashed]  (8,-16) to (8,-18);
\draw [-, dashed]  (12,-16) to (12,-18);
\draw [-, dashed]  (16,-16) to (16,-18);

\draw [-, dashed]  (2,-4) to (4,-4);
\draw [-, dashed]  (2,-8) to (4,-8);
\draw [-, dashed]  (2,-12) to (4,-12);
\draw [-, dashed]  (2,-16) to (4,-16);

\draw [-, dashed]  (16,-4) to (18,-4);
\draw [-, dashed]  (16,-8) to (18,-8);
\draw [-, dashed]  (16,-12) to (18,-12);
\draw [-, dashed]  (16,-16) to (18,-16);

\draw [-]  (4,-4) to (4,-16);
\draw [-]  (16,-4) to (16,-16);
\draw [-]  (8,-4) to (8,-16);
\draw [-]  (12,-4) to (12,-16);

\draw [-]  (4,-4) to (16,-4);
\draw [-]  (4,-8) to (16,-8);
\draw [-]  (4,-12) to (16,-12);
\draw [-]  (4,-16) to (16,-16);

\node  (c1) at (10,-9.25) {\Large{$\color{blue}NW_{\bm{\gamma}}(j,k)\color{black}$}};
\node  (c2) at (10,-10) {\Large{$=$}};
\node  (c3) at (10,-10.75) {\Large{$\color{red}H_{\bm{\gamma}}(m,n)\color{black}$}};

\node  (a1) at (14,-9.25) {\Large{$\color{blue}NW_{\bm{\gamma}}(j,k+1)\color{black}$}};
\node  (a2) at (14,-10) {\Large{$=$}};
\node  (a3) at (14,-10.75) {\Large{$\color{red}H_{\bm{\gamma}}(m+1,n)\color{black}$}};

\node  (b1) at (6,-9.25) {\Large{$\color{blue}NW_{\bm{\gamma}}(j,k-1)\color{black}$}};
\node  (b2) at (6,-10) {\Large{$=$}};
\node  (b3) at (6,-10.75) {\Large{$\color{red}H_{\bm{\gamma}}(m-1,n)\color{black}$}};

\node  (e1) at (6,-5.25) {\large{$\color{blue}NW_{\bm{\gamma}}(j-1,k-1)\color{black}$}};
\node  (e2) at (6,-6) {\Large{$=$}};
\node  (e3) at (6,-6.75) {\Large{$\color{red}H_{\bm{\gamma}}(m,n-1)\color{black}$}};

\node  (f1) at (10,-5.25) {\Large{$\color{blue}NW_{\bm{\gamma}}(j-1,k)\color{black}$}};
\node  (f2) at (10,-6) {\Large{$=$}};
\node  (f3) at (10,-6.75) {\large{$\color{red}H_{\bm{\gamma}}(m+1,n-1)\color{black}$}};

\node  (g1) at (14,-5.25) {\large{$\color{blue}NW_{\bm{\gamma}}(j-1,k+1)\color{black}$}};
\node  (g2) at (14,-6) {\Large{$=$}};
\node  (g3) at (14,-6.75) {\large{$\color{red}H_{\bm{\gamma}}(m+2,n-1)\color{black}$}};

\node  (h1) at (14,-13.25) {\large{$\color{blue}NW_{\bm{\gamma}}(j+1,k+1)\color{black}$}};
\node  (h2) at (14,-14) {\Large{$=$}};
\node  (h3) at (14,-14.75) {\Large{$\color{red}H_{\bm{\gamma}}(m,n+1)\color{black}$}};

\node  (i1) at (10,-13.25) {\Large{$\color{blue}NW_{\bm{\gamma}}(j+1,k)\color{black}$}};
\node  (i2) at (10,-14) {\Large{$=$}};
\node  (i3) at (10,-14.75) {\large{$\color{red}H_{\bm{\gamma}}(m-1,n+1)\color{black}$}};

\node  (j1) at (6,-13.25) {\large{$\color{blue}NW_{\bm{\gamma}}(j+1,k-1)\color{black}$}};
\node  (j2) at (6,-14) {\Large{$=$}};
\node  (j3) at (6,-14.75) {\large{$\color{red}H_{\bm{\gamma}}(m-2,n+1)\color{black}$}};

\end{tikzpicture}
}
\end{center}

\caption{A portion of the Number Wall of a sequence  $\bm{\gamma}$ centered at the indices $(j,k)=\iota(m,n)$, and the corresponding index set.}
\label{fig1}
\end{figure}

\noindent The goal here  is not  to develop a structural theory of Number Walls as is done in~\cite{ANL}. In fact, this concept and the corresponding indexation~\eqref{indexbij} are introduced solely for the sake of obtaining an elegant formulation of the following statement first proved by Lunnon in~\cite[p.9]{L}~:

\begin{lem}[Square Window Lemma]\label{sqwinthe}
Given a Number Wall with entries taking values in an integral domain, the zero entries can only occur in the form of \emph{windows}; that is, in square shapes with horizontal and vertical edges.
\end{lem}

\noindent Lemma~\ref{sqwinthe} implies for instance that if the diagonal $H_{\bm{\gamma}}(m,n-1)$, $H_{\bm{\gamma}}(m,n)$, \mbox{$H_{\bm{\gamma}}(m,n+1)$} vanishes in Figure~\ref{fig1}, then all entries in the $3\times 3$ square depicted therein also vanish.

\subsection{On the Polynomial Approximation Sequence  of a Laurent Series}\label{secpolyapprox}


\noindent Denote by
{\allowdisplaybreaks
\begin{align*}\label{contfrac}
\Gamma(t)\; &=\; \left[A_{\bm{\gamma}}(0; t)\;\; ;\; \;A_{\bm{\gamma}}(1; t)\;\;,\;\; \dots\;\; ,\;\; A_{\bm{\gamma}}(h; t)\;\; ,\; \; \dots\right] \\
&=\; A_{\bm{\gamma}}(0; t)+\cfrac{1}{A_{\bm{\gamma}}(1; t)+\cfrac{1}{\ddots\qquad +\cfrac{1}{A_{\bm{\gamma}}(h; t)+\cfrac{1}{\qquad \ddots}}}}\cdotp
\end{align*}
}
\noindent the continued fraction expansion of the Laurent series $\Gamma(t)$, and by
\begin{equation*}\label{conv}
 \frac{P_{\bm{\gamma}}(h; t)}{Q_{\bm{\gamma}}(h; t)}\;=\; \left[A_{\bm{\gamma}}(0; t)\; \; ;\; \;A_{\bm{\gamma}}(1; t)\;\;,\;\; \dots\;\;,\;\; \;A_{\bm{\gamma}}(h; t)\right]
\end{equation*}
its $h^{\textrm{th}}$ convergent. Here, the partial quotients \sloppy $\left(A_{\bm{\gamma}}(h;t)\right)_{h\ge 0}$, the numerators $\left(P_{\bm{\gamma}}(h;t)\right)_{h\ge 0}$ and the denominators $\left(Q_{\bm{\gamma}}(h;t)\right)_{h\ge 0}$ are all polynomial sequences with values in $\K\left[t\right]$, and one has  that $\deg\left(A_{\bm{\gamma}}(h;t)\right)\ge 1$ for all $h\ge 1$. 
For standard references on the theory of continued fractions in positive characteristics, see, e.g., \cite{BN} and the references within. \\

\noindent For each $h\ge 0$, The polynomials $P_{\bm{\gamma}}(h;t)$ and $Q_{\bm{\gamma}}(h;t)$ are coprime  and meet the pro\-per\-ty that
\begin{equation}\label{inegapproxcinv}
\left| Q_{\bm{\gamma}}\left(h;t\right) \cdot \Gamma(t)-P_{\bm{\gamma}}\left(h;t\right)\right|\;=\;\frac{1}{\left|Q_{\bm{\gamma}}\left(h+1;t\right) \right|}\cdotp
\end{equation}
Furthermore, given polynomials $\left(P(t), Q(t)\right)\in \K[t]\times \K[t]$, it holds that
\begin{equation}\label{convcharac}
\left|Q(t)\cdot \Gamma(t)-P(t)\right|\;<\;  \frac{1}{\left|Q( t)\right|} \qquad \Longleftrightarrow \qquad \exists h\ge 0\; :\; \frac{P(t)}{Q(t)}=\frac{P_{\bm{\gamma}}(h; t)}{Q_{\bm{\gamma}}(h; t)}\cdotp
\end{equation}

\noindent Throughout, given an integer $m\in\Z$,
\begin{quote}
\emph{$\left(P_{\bm{\gamma}}(m,h;t)/Q_{\bm{\gamma}}(m,h;t)\right)_{k\ge 0}$ denotes the sequence of  the convergents of the Laurent series $\Gamma_m(t)$ defined in~\eqref{doublylaurent}. }
\end{quote}
In particular when $m=1$, this sequence coincides with $\left(P_{\bm{\gamma}}(h; t)/Q_{\bm{\gamma}}(h; t)\right)_{h\ge 1}$. The \emph{$\left(m,h\right)^{\textrm{th}}$ normal order} of the Laurent series $\Gamma(t)$ is then defined as the integer \mbox{$q_{\bm{\gamma}}\left(m,h\right)\ge 0$} such that
\begin{equation}\label{defdegnorord}
q_{\bm{\gamma}}\left(m,h\right)\;=\; \deg Q_{\bm{\gamma}}(m,h;t).
\end{equation}
\noindent In the above notations, a distinction must be made between the role played by the indices $h$ and $(m,h)$ on the one hand and by the variable $t$ on the other. \\

\noindent The proof of Main Theorem~\ref{mainthm}  relies on the properties of the polynomial
{\allowdisplaybreaks
\begin{align}\label{detqmnt}
S_{\bm{\gamma}}(m,n;t)=
\begin{vmatrix}
\mathcal{H}_{\bm{\gamma}}(m,n,n+1)\vspace{2mm} \\
\bm{w}_{n+1}(t)
\end{vmatrix}
&=
\begin{vmatrix}
\gamma_m &\cdots & \gamma_k& \cdots & \gamma_{m+n-1}&  \gamma_{m+n}\\
\vdots & \iddots && \iddots &\vdots& \vdots\\
\gamma_k& & \gamma_{m+n-1} &&\gamma_{r}& \gamma_{r+1}\\
\vdots & \iddots &&\iddots&\vdots& \vdots\\
\gamma_{m+n-1} & \cdots &\gamma_r&\cdots & \gamma_{m+2n-2}&  \gamma_{m+2n-1}\vspace{3mm}\\
1&\cdots &t^k& \cdots &t^{n-1}&t^n
\end{vmatrix},
\end{align}
}
\!\!which is well-defined for all integers $m\in\Z$ and $n\ge 1$. 
It is from now on referred to as the \emph{polynomial approximation of order $\left(m,n\right)$}. The following statement collects more or less standard properties of this polynomial in relation with the sequence $\left(q_{\bm{\gamma}}\left(m,h\right)\right)_{m\in\Z,h\ge 0}$. 

\begin{lem}[\emph{Properties of the polynomial approximation sequence in relation with the normal orders}]\label{closedromlem} Keep the above notations and definitions. Then~:\\
\begin{itemize}
\item[\textbf{(1)}] \emph{[Values of Hankel determinants by specialisation]}. Given integers $m\in\Z$ and $n\ge 1$, the coefficient of $t^n$  in the polynomial  approximation $S_{\bm{\gamma}}(m,n;t)$ is the determinant $H_{\bm{\gamma}}(m,n)$, and its constant term is the signed determinant \mbox{$(-1)^{n}\cdot H_{\bm{\gamma}}(m+1,n)$}.\\

\item[\textbf{(2)}] \emph{[Closed determinant form  for convergents in the case that the polynomial appro\-xi\-mant has maximal rank]}. Given integers $m\in \Z$ and $n\ge 1$ such that $S_{\bm{\gamma}}(m,n;t)\neq 0$ (which amounts to claiming that $\rank(\mathcal{H}_{\bm{\gamma}}(m,n,n+1))=n$), there exists an index $h\ge 0$ for which the relation
\begin{equation}\label{closedform}
Q_{\bm{\gamma}}(m,h;t)\;=\; c\cdot S_{\bm{\gamma}}(m, n; t)
\end{equation}
is met for some nonzero scalar $c=c(m,n,h)\in\K$.    Moreover, the inequalities
\begin{equation}\label{lem24_2}
q_{\bm{\gamma}}(m,h)\;\le\; n\qquad \textrm{and}\qquad |Q_{\bm{\gamma}}(m,h;t)\cdot \Gamma_m(t) - P_{\bm{\gamma}}(m,h;t)|\;\le\; 2^{-n-1}
\end{equation}
then also hold.  

\item[\textbf{(3)}] \emph{[Vanishing of Hankel determinants, quality of approximation \& normal orders]}. Given integers $m\in\Z$ and $n\ge 1$, there exists an index $h\ge 0$ such that $q_{\bm{\gamma}}(m,h)=n$ (that is, such that $n$ is the $(m,h)^{\textrm{th}}$ normal order) if, and only if, the Hankel determinant $H_{\bm{\gamma}}\left(m, n\right)$ does not vanish. In this case,
\begin{equation}\label{inegapproxcinvbis}
\left| Q_{\bm{\gamma}}\left(m,h;t\right) \cdot   \Gamma_m(t)-P_{\bm{\gamma}}\left(m,h;t\right)\right|\;=\; 2^{-q_{\bm{\gamma}}(m, h+1)}.\\
\end{equation}

\item[\textbf{(4)}] \emph{[Normal order and quality of approximation when  the polynomial approximant has a deficient rank]} Given integers $m\in \Z$ and $n\ge 1$ such that $S_{\bm{\gamma}}(m,n;t)= 0$ (which amounts to claiming that $\rank(\mathcal{H}_{\bm{\gamma}}(m,n,n+1))<n$), one has that
\begin{equation*}
H_{\bm{\gamma}}(m,n)\; =\; H_{\bm{\gamma}}(m+1,n) \; =\;  H_{\bm{\gamma}}(m-1,n+1) \; =\;  H_{\bm{\gamma}}(m,n+1) \; =\;  0.
\end{equation*}
\end{itemize}
\color{black}
\end{lem}

\begin{rem}\label{rem0}
Throughout, the coefficients of the denominators of convergents are normalised in such a way that  $c=1$ in~\eqref{closedform} when $q_{\bm{\gamma}}(m,h)=n$.\qed
\end{rem}

\begin{proof}[Proof of Lemma~\ref{closedromlem}]
Property~(1) is immediate from identity~\eqref{detqmnt} defining the polynomial approximation of order $(m,n)$.

\paragraph{$\bullet$} As for Property~(2), define the index $h$ as the largest one for which $q_{\bm{\gamma}}(m,h)\le n$ (recall here the definition of the normal order $q_{\bm{\gamma}}(m,h)$ in~\eqref{defdegnorord}). By maximality of this choice, equation~\eqref{inegapproxcinv} implies that the second inequality in~\eqref{lem24_2} is also verified. \\

\noindent  In order to establish relation~\eqref{closedform}, 
 expand 
 the polynomial  $Q_{\bm{\gamma}}(m,h;t)$  
in the form
\begin{equation}\label{defdenomexp}
Q_{\bm{\gamma}}(m,h;t)\;=\; \sum_{s=0}^{n}\rho_{\bm{\gamma}}\left(m,s\right)\cdot t^s,
\end{equation}
and set $\bm{\overline{\rho}}_{\bm{\gamma}}\left(m,n\right)=\left(\rho_{\bm{\gamma}}\left(m,0\right), \dots, \rho_{\bm{\gamma}}\left(m,n\right)\right)$. From the second inequality in~\eqref{lem24_2},  the coefficients of $t^{-1}, t^{-2}, \dots, t^{-n}$ vanish in the pro\-duct $Q_{\bm{\gamma}}(m,h;t)\cdot \Gamma_m(t)$. This 
yields the system of equations
\begin{equation}\label{systcvgt}
\sum_{s=0}^{n}\rho_{\bm{\gamma}}\left(m,s\right)\cdot \gamma_{r+s}=0\qquad \textrm{for all }\qquad  r\in\llbracket m, m+n-1\rrbracket.
\end{equation}
Adding relation~\eqref{defdenomexp} as an $(n+1)^{\textrm{st}}$ equation to this system  leads one to the matrix equation
\begin{equation*}
\begin{pmatrix}
\mathcal{H}_{\bm{\gamma}}(m, n , n+1)\\
\bm{w}_{n+1}(t)
\end{pmatrix}
\cdot  \bm{\overline{\rho}}_{\bm{\gamma}}\left(m,n\right)^T\;=\; \left(\underbrace{0, \dots, 0}_{n \textrm{ times}}\;, \; Q_{\bm{\gamma}}(m,h;t)\right)^T.
\end{equation*}
Given an integer $0\le s\le n$, Cramer's rule provides a closed form expression for the coefficients $\rho_{\bm{\gamma}}\left(m,s\right)$, namely
\begin{equation}\label{rhoqht}
\rho_{\bm{\gamma}}\left(m,s\right)\;=\; (-1)^{n+s}\cdot Q_{\bm{\gamma}}(m,h;t)\cdotp \frac{ H^{(s)}_{\bm{\gamma}}(m,n)}{\begin{vmatrix}
\mathcal{H}_{\bm{\gamma}}(m,n,n+1)\\
\bm{w}_{n+1}(t)
\end{vmatrix}}\cdotp
\end{equation}
Here, $H^{(s)}_{\bm{\gamma}}(m,n)$ is the minor of the matrix $\mathcal{H}_{\bm{\gamma}}(m,n,n+1)$ obtained after removing the $(s+1)^{\textrm{st}}$ column. Upon setting $s$ to be the normal order $q_{\bm{\gamma}}(m,k)$,  the polynomial $Q_{\bm{\gamma}}(m,h;t)$ happens to be a nonzero scalar multiple of $S_{\bm{\gamma}}(m,n;t)$. 
This completes the proof of Property~(2).

\paragraph{$\bullet$}  As far as Property~(3) is concerned, relation~\eqref{inegapproxcinvbis} is nothing but a rephra\-sing of identity~\eqref{inegapproxcinv} in terms of normal orders. \\

\noindent To establish the equivalence stated therein, assume first that $H_{\bm{\gamma}}(m,n)\neq 0$. This determinant equals the minor $H^{(n)}_{\bm{\gamma}}(m,n)$ appearing in~\eqref{rhoqht} (obtained when setting $s=n$). Since, from Property~(1), the leading coefficient of the polynomial approximant $S_{\bm{\gamma}}(m,n;t)$ is $H_{\bm{\gamma}}(m,n)$, the (nonzero) denominator of convergent $Q_{\bm{\gamma}}(m,h;t)$ involved in this relation is readily seen to be of degree $n$. \\

\noindent Conversely, assume that $H_{\bm{\gamma}}(m,n)=\det\left(\mathcal{H}_{\bm{\gamma}}(m,n)\right)= 0$ and rewrite  system~\eqref{systcvgt} as $$\mathcal{H}_{\bm{\gamma}}(m,n)\cdot \bm{\overline{\rho}}_{\bm{\gamma}}\left(m,n-1\right)^T+\rho_{\bm{\gamma}}(m,n)\cdot \bm{\gamma}\left[m,n\right]^T\;=\;\bm{0},$$ where $ \bm{\gamma}\left[m,n\right]=\left(\gamma_{m+n-1}, \dots, \gamma_{m+2n-1}\right)$. Choose then a nontrivial element $ \bm{\overline{\rho}}_{\bm{\gamma}}\left(m,n-1\right)$ in the kernel of the matrix $\mathcal{H}_{\bm{\gamma}}(m,n)$ and set $\rho_{\bm{\gamma}}(m,n)=0$. The polynomial $Q_{\bm{\gamma}}(m,h;t)$ defined in~\eqref{defdenomexp} from the coefficients $\bm{\overline{\rho}}_{\bm{\gamma}}\left(m,n\right)=\left(\bm{\overline{\rho}}_{\bm{\gamma}}\left(m,n-1\right), 0\right)$ has thus degree strictly less than $n$ and meets by construction the second inequality in~\eqref{lem24_2}. From identity~\eqref{inegapproxcinv}, one deduces that $q_{\bm{\gamma}}(m,h)<n<q_{\bm{\gamma}}(m,h+1)$, which shows that the integer $n$ cannot be a normal order and thus concludes the proof of the sought claim.



\paragraph{$\bullet$} To establish Property~(4), note that under the assumption that the rank of the matrix $\mathcal{H}_{\bm{\gamma}}(m,n,n+1)$ is not maximal,   the two maximal minors $H_{\bm{\gamma}}(m,n)$ and $H_{\bm{\gamma}}(m+1,n)$ both vanish. Also, a square matrix obtained by appending a row to the matrix $\mathcal{H}_{\bm{\gamma}}(m,n,n+1)$ is singular (as can be seen either by expanding the determinant along this row or by considering a nontrivial vanishing linear combination of the rows in $\mathcal{H}_{\bm{\gamma}}(m,n,n+1)$, which exists by assumption). Since the two square matrices $\mathcal{H}_{\bm{\gamma}}(m,n+1,n+1)$ and $\mathcal{H}_{\bm{\gamma}}(m-1,n+1,n+1)$ can be achieved this way, the corresponding determinants   $H_{\bm{\gamma}}(m,n+1)$ and  $H_{\bm{\gamma}}(m-1,n+1)=0$ vanish, as was to be verified.

\end{proof}

\section[Recurrence Relations in the Polynomial Approxi\-mation Sequence]{Recurrence Relations in the Polynomial Appro\-ximation Sequence of the Paperfolding Laurent Series}\label{sec3}

\noindent The main goal in this section is to establish (doubly--indexed) recurrence relations for the polynomial approximation sequence $\left(S_{\bm{\lambda}}\left(m,n;t\right)\right)_{m\in\Z,n\ge 1}$ (defined in~\eqref{detqmnt}), where $\bm{\lambda}$ is the paperfolding sequence over the alphabet $(1,-1)$. Relations for the sequence of Hankel determinants $\left(H_{\bm{\lambda}}\left(m,n\right)\right)_{m\in\Z,n\ge 1}$  (defined in~\eqref{defhankelmn}) are then derived from them.

\subsection{Preliminaries}\label{prelim}

\noindent The recurrence relations for the sequences $\left(S_{\bm{\lambda}}\left(m,n;t\right)\right)_{m\in\Z,n\ge 1}$  and  $\left(H_{\bm{\lambda}}\left(m,n\right)\right)_{m\in\Z,n\ge 1}$  comprise a set of four   equations in each case. They involve not only the coefficients of the Laurent series  $\Lambda(t)$, but also those of the two shifted series  $(t+1)\cdot\Lambda(t)$ and  $(t+1)^2\cdot\Lambda(t)$. Part of the argument to establish them is valid in the framework of an arbitrary sequence $\bm{\gamma}$. The corresponding theory is then developed at this level of generality. To this end, some additional notations and definitions are first introduced. \\


\noindent Given an integer $n\ge 1$, set
\begin{equation*}\label{defespivecbis}
\bm{\varepsilon}(n)\; :=\; \left((-1)^{n-1}, (-1)^{n-2}, \dots, -1,1\right)\in\K^n.
\end{equation*}

\noindent Let $\left(\gamma^*_k\right)_{k\ge 1}$ and $\left(\gamma^{**}_k\right)_{k\ge 1}$ be the coefficients of the  fractional parts of the Laurent series \mbox{$(t+1)\cdot\Gamma(t)$} and  \mbox{$(t+1)^2\cdot\Gamma(t)$}, respectively.
This is saying that for all $k\ge 1$,
\begin{equation*}
\gamma^*_k\;=\; \gamma_k+\gamma_{k+1} \qquad \qquad \textrm{ and }\qquad \qquad \gamma^{**}_k\;=\; \gamma_k+2\gamma_{k+1}+\gamma_{k+2}.
\end{equation*}
\noindent These sequences are extended to doubly infinite sequences $\bm{\gamma^*}=\left(\gamma^*_k\right)_{k\in\Z}$ and $\bm{\gamma^{**}}= \left(\gamma^{**}_k\right)_{k\in\Z}$ obtained from the values for $k\le 0$ of the extended sequence  $\bm{\gamma}= \left(\gamma_k\right)_{k\in\Z}$ introduced in~\eqref{defdoublyextendedeq}.\\

\noindent To simplify the notation, given integers $m\in\Z$ and $n\ge 1$, the Hankel determinant $H_{\bm{\gamma^*}}(m,n)$ of index $(m,n)$ of the sequence $\bm{\gamma^*}$ is denoted by $F_{\bm{\gamma}}(m,n)$, and the Hankel determinant $H_{\bm{\gamma^{**}}}(m,n)$   of the sequence $\bm{\gamma^{**}}$ is denoted by $G_{\bm{\gamma}}(m,n)$. These two quantities admit  alternative representations as determinants of respective sizes $n+1$ and $n+2$, namely
\begin{equation}\label{deterformFG}
F_{\bm{\gamma}}(m,n)=
\begin{vmatrix}
\mathcal{H}_{\bm{\gamma}}(m,n,n+1)\\
\bm{\varepsilon}(n+1)
\end{vmatrix}
\quad \textrm{and}\quad
G_{\bm{\gamma}}(m,n)= -
\begin{vmatrix}
\mathcal{H}_{\bm{\gamma}}(m,n+1, n+1)& \bm{\varepsilon}^T(n+1)\\
\bm{\varepsilon}(n+1)& 0
\end{vmatrix}.
\end{equation}
Here, the matrices $\mathcal{H}_{\bm{\gamma}}(m,n,n+1)$ and $\mathcal{H}_{\bm{\gamma}}(m,n+1, n+1)$ are defined in~\eqref{defhankelmnp}. 
These identities are easily obtained by elementary row and columns operations. (For $F_{\bm{\gamma}}(m,n)$, it suffices to add to each column the following one starting from the left and then to expand the determinant along the last row; for $G_{\bm{\gamma}}(m,n)$, it suffices to add to each column the following one starting from the left, to add to each row the following one starting from the top, and then to successively expand along the last row and then the last column.) \\

\noindent Similarly, the polynomial approximation sequences $\left(S_{\bm{\gamma^*}}(m,n;t)\right)_{m\in\Z,n\ge 1}$ and $\left(S_{\bm{\gamma^{**}}}(m,n;t)\right)_{m\in\Z,n\ge 1}$ defined in~\eqref{detqmnt}  and associated to the sequences $\bm{\gamma^*}$ and $\bm{\gamma^{**}}$ are denoted by $\left(R_{\bm{\gamma}}(m,n;t)\right)_{m\in\Z,n\ge 1}$ and $\left(V_{\bm{\gamma}}(m,n;t)\right)_{m\in\Z,n\ge 1}$, respectively. In the same way as above, given integers $m\in\Z,n\ge 1$, they admit alternative representations as determinants of respective sizes $n+2$ and $n+3$,  namely
\begin{equation}\label{detf}
 R_{\gamma}\left(m,n;t\right)\;=\; -
\begin{vmatrix}
\mathcal{H}_{\bm{\gamma}}(m,n+1,n+1)& \bm{\varepsilon}^T(n+1)\vspace{2mm}\\
\bm{w}_{n+1}\left(t\right) & 0
\end{vmatrix}
\end{equation}
and
\begin{equation}\label{detg}
 \left(t+1\right)\cdot V_{\gamma}\left(m,n;t\right)\;=\; -
\begin{vmatrix}
\mathcal{H}_{\bm{\gamma}}(m,n+1,n+2)& \bm{\varepsilon}^T(n+1)\vspace{2mm}\\
\bm{\varepsilon}(n+2)& 0\vspace{2mm}\\
\bm{w}_{n+2}\left(t\right) & 0
\end{vmatrix}.
\end{equation}

\noindent The quantities  $G_{\bm{\gamma}}(m,n)$, $R_{\gamma}\left(m,n;t\right)$ and $\left(t+1\right)\cdot V_{\gamma}\left(m,n;t\right)$ are introduced here for values of $n\ge 1$. However, their definitions also make sense for $n=0$. 
They shall thus be assumed to be defined for all $n\ge 0$.\\

\noindent  Given $n,p\ge 1$, denote by
\begin{equation}\label{defEmatr}
\bm{E}(n,p)\; :=\;
\begin{pmatrix}
\bm{\varepsilon}(p);\\
-\bm{\varepsilon}(p);\\
\vdots\\
(-1)^{n-1}\cdot \bm{\varepsilon}(p)
\end{pmatrix}
\in\K^{n\times p}
\end{equation}
the $n\times p$ matrice with $\pm 1$ entries and by  $\bm{0}(n, p)$ the $n\times p$ zero matrix. Let also
\begin{equation}\label{def1n}
\bm{0}(n)\; :=\;\left(0, \dots, 0\right)\in\K^n \qquad \textrm{and} \qquad \bm{1}(n)\; :=\;\left(1, \dots, 1\right)\in\K^n.
\end{equation}




\subsection{The Recurrence Relations and some of their Consequences}


It is easily checked from its definition in~\eqref{pfdalphbab} that the sequence $\bm{\lambda}[1,-1]$ obeys the recurrence relations
\begin{equation}\label{reclambda}
\lambda_{4k}\;=\; 1,\qquad \lambda_{4k+1}\;=\;\lambda_{2k}, \qquad \lambda_{4k+2}\;=\;-1, \qquad \textrm{and}\qquad  \lambda_{4k+3}\;=\;\lambda_{2k+1}
\end{equation}
for all $k\ge 0$. These relations are extended to negative values of the index upon setting $\lambda_{-1}=1$. The resulting doubly infinite sequence $\left(\lambda_k\right)_{k\in\Z}$  is still denoted $\bm{\lambda}$ (as in~\eqref{pfdalphbab}). Thus,
\begin{align}\label{symmpaperf}
&\left(\dots, \; \lambda_{-4}, \;\lambda_{-3},  \;\lambda_{-2},  \;\lambda_{-1},  \;\lambda_0,   \;\lambda_1,  \;\lambda_2,  \;\lambda_3,  \;\lambda_4,  \;\dots\right)\nonumber \\
& \qquad\qquad \qquad \qquad \qquad\qquad =\; \left(\dots, \; 1,  \; \lambda_{-2},  \; -1, \; \lambda_{-1} ,  \;1, \; \lambda_0 ,  \;-1, \; \lambda_1 , \; 1, \; \dots\right).
\end{align}
\noindent In particular, $\lambda_k=1$ whenever $k\equiv 0 \pmod{4}$ and $\lambda_k=-1$ whenever $k\equiv 2 \pmod{4}$.

\begin{prop}[Recurrence relations for the polynomial approximation sequence of the paperfolding Laurent series]\label{proprec} Keep the above notations and conventions. Then, the polynomial approximation sequence of the paperfolding sequence meets the following recurrence relations~:
\begin{itemize}
\item[\textbf{(A)}] \; \;  when $m\in\Z$ and $n\ge 1$,
\begin{align*}
S_{\bm{\lambda}}(2m+1,2n;t)\;=\;  &H_{\bm{\lambda}}\left(m+1,n\right)\cdot S_{\bm{\lambda}}\left(m,n;t^2\right)\\
& \qquad\;\; \qquad - G_{\bm{\lambda}}\left(m+1, n-1\right) \cdot \left(t^2+1\right)\cdot V_{\bm{\lambda}}\left(m,n-1; t^2\right)\\
&\qquad\;\; \qquad + (-1)^{m+1}\cdot F_{\bm{\lambda}}(m,n)\cdot t\cdot R_{\bm{\lambda}}\left(m+1,n-1;t^2\right);
\end{align*}
\item[\textbf{(B)}] \; \;  when  $m\in\Z$ and $n\ge 2$,
\begin{align*}
S_{\bm{\lambda}}(2m+1,2n-1;t)\;=\;  &\left(-1\right)^{m}\cdot F_{\bm{\lambda}}\left(m+1, n-1\right)\cdot R_{\lambda}\left(m,n-1;t^2\right)\\
&\qquad +  H_{\bm{\lambda}}\left(m,n\right)\cdot t\cdot S_{\bm{\lambda}}\left(m+1, n-1; t^2\right)\\
&\qquad - G_{\bm{\lambda}}\left(m, n-1\right)\cdot t\cdot \left(t^2+1\right)\cdot  V_{\bm{\lambda}}\left(m+1, n-2; t^2\right);
\end{align*}
\item[\textbf{(C)}]\;\; when  $m\in\Z$ and $n\ge 1$,
\begin{align*}
S_{\bm{\lambda}}(2m,2n;t)\;=\;  &(-1)^n\cdot\left(\left(-1\right)^{m} \cdot F_{\bm{\lambda}}\left(m, n\right)\cdot t\cdot R_{\bm{\lambda}}\left(m, n-1; t^2\right)\right.\\
&\quad+H_{\bm{\lambda}}\left(m,n\right)\cdot S_{\bm{\lambda}}\left(m,n;t^2\right)\\
& \quad \left. + G_{\bm{\lambda}}\left(m,n-1\right)\cdot \left(t^2+1\right)\cdot V_{\bm{\lambda}}\left(m, n-1; t^2\right)\right);
\end{align*}
\item[\textbf{(D)}]\;\;  when $m\in\Z$ and $n\ge 2$,
\begin{align*}
S_{\bm{\lambda}}(2m,2n-1;t)\;=\;  &(-1)^n\cdot\left(\left(-1\right)^{m+1}\cdot F_{\bm{\lambda}}\left(m, n-1\right)\cdot t\cdot R_{\bm{\lambda}}\left(m, n-1; t^2\right)\right.\\
&\qquad\quad\quad+H_{\bm{\lambda}}\left(m,n\right)\cdot S_{\bm{\lambda}}\left(m, n-1; t^2\right)\\
&\left. \qquad\quad\quad + G_{\bm{\lambda}}\left(m,n-1\right)\cdot \left(t^2+1\right)\cdot V_{\bm{\lambda}}\left(m, n-2; t^2\right)\right).
\end{align*}
\end{itemize}
\end{prop}

\noindent The above four recurrences are established in the next Section~\ref{detsubsec}. They should  be seen as a compact formulation of a  number of other recurrences involving the sequences $\left(H_{\bm{\lambda}}(m,n)\right)_{m\in\Z,n\ge 1}$, $\left(F_{\bm{\lambda}}(m,n)\right)_{m\in\Z,n\ge 1}$ and $\left(G_{\bm{\lambda}}(m,n)\right)_{m\in\Z,n\ge 1}$ upon identifying coefficients or specialising the values of the variable $t$. As such, Proposition~\ref{proprec} constitutes an alternative to the usually large number of relations needed to establish recurrences for Hankel determinants of automatic sequences. For instance, the foundational work by Allouche \emph{et.~al.}~\cite{APWW} on Hankel determinants modulo~2 formed from the Thue-Morse sequence relies on a system of 16 recurrence relations. Similarly, Guo \emph{et.~al.}~\cite{GWW} introduce a system of 18 recurrence relations to analyse properties of the Hankel determinants of the paperfolding sequence over the alphabet $(0,1)$ reduced modulo 2. \\


\noindent The following is a crucial and  immediate consequence of the above proposition. It illustrates the kind of the identities that can be deduced from it upon identification of the coefficients. Recall that given an integer $m\in\Z$, the definitions of the quantities $G_{\bm{\gamma}}(m,n)$, $R_{\gamma}\left(m,n;t\right)$ and $\left(t+1\right)\cdot V_{\gamma}\left(m,n;t\right)$ have been extended in the previous section to all values of $n\ge 0$.

\begin{coro}[Recurrence relations for the Hankel determinants of the paperfolding sequence over the integers]\label{labelhenkelbis} The Hankel determinants of the paperfolding sequence meet the following recurrence relations~:
\begin{itemize}
\item[\textbf{(W)}]\;\;  when $m\in\Z$ and $n\ge 1$, \begin{align*}H_{\bm{\lambda}}(2m+&1, 2n)\\ &= H_{\bm{\lambda}}(m, n)\cdot H_{\bm{\lambda}}(m+1, n) -G_{\bm{\lambda}}(m, n-1)\cdot G_{\bm{\lambda}}(m+1, n-1);\end{align*}  
\item[\textbf{(X)}]\;\;  when $m\in\Z$ and $n\ge 2$, \begin{align*}H_{\bm{\lambda}}(2m+&1, 2n-1)\\&= H_{\bm{\lambda}}(m, n)\cdot H_{\bm{\lambda}}(m+1, n-1) - G_{\bm{\lambda}}(m, n-1)\cdot G_{\bm{\lambda}}(m+1, n-2);\end{align*} 
\item[\textbf{(Y)}] \; \; when $m\in\Z$ and $n\ge 1$, $$H_{\bm{\lambda}}(2m, 2n)\; =\;  (-1)^n\cdot\left( \left(H_{\bm{\lambda}}(m, n)\right)^2+\left(G_{\bm{\lambda}}(m, n-1)\right)^2\right);$$ 
\item[\textbf{(Z)}] \; \; when $m\in\Z$ and $n\ge 2$, $$H_{\bm{\lambda}}(2m, 2n-1)\; =\; \left(-1\right)^{n+m+1}\cdot \left(F_{\bm{\lambda}}(m, n-1)\right)^2.$$
\end{itemize}
\end{coro}

\begin{proof}
Identify the coefficients of the leading powers in Proposition~\ref{proprec}. The claims then follow from Property~(1) in Lemma~\ref{closedromlem}.
\end{proof}

\noindent Identity (Y) is the key to the obstruction preventing the paperfolding Laurent series $\Lambda(t)$ from satisfying $P(t)$--LC over finite fields with characteristics $\ell \equiv 3\pmod{4}$.
Given its importance, its main implication is recorded as a statement on its own~:

\begin{coro}[Arithmetic Obstruction]\label{keyobstr}
Let $\ell$ be a prime congruent to 3 modulo 4, and let $\K$ be a field with characteristic $\ell$. Assume that $m\in\Z$ and $n\ge 1$ are integers such that
$$H_{\bm{\lambda}}(2m, 2n)\;=\; 0.$$ Then, $$H_{\bm{\lambda}}(m, n)\;=\;0 \;\;\qquad \textrm{and}\;\;\qquad G_{\bm{\lambda}}(m, n-1)\;=\; 0.$$
\end{coro}

\begin{proof}
Since the paperfolding sequence $\bm{\lambda}$ is defined over the alphabet $(-1,1)$,  the determinants $H_{\bm{\lambda}}(2m, 2n), H_{\bm{\lambda}}(m, n)$ and $G_{\bm{\lambda}}(m, n-1)$ all take their values in the finite field $\F_\ell$. The statement then
follows immediately from identity (Y) and from  the fact that -1 is not a quadratic residue modulo $\ell$.
\end{proof}

\noindent Corollary~\ref{labelhenkelbis} also provides a first piece of information on the sizes (i.e.~the side lengths) of the windows in the number wall of the paperfolding sequence over the  fields under consideration. 

\begin{coro}[Unevenness of the sizes of the windows]\label{odwind}
Let $\ell$ be a prime congruent to 3 modulo 4, and let $\K$ be a field with characteristic $\ell$.  Recall that  $NW_{\bm{\lambda}}$ denotes the Number Wall of the paperfolding sequence $\bm{\lambda}$.
Then, all windows in $NW_{\bm{\lambda}}$ have odd size. More precisely, if for some integers $m\in\Z$ and $n\ge 1$, it holds that $H_{\bm{\lambda}}(2m,2n)= 0$,
then  all entries in the $3\times 3$ square centered at this entry also vanish.
\end{coro}

\begin{proof}
\noindent The last claim is  established first. To this end, consider a window of size at least two in $NW_{\bm{\lambda}}$.
Then, one of its entries is of the form $H_{\bm{\lambda}}(2m,2n)$
for some integers $m\in\Z,n\ge 1$ (this is  compatible with the choice of indexation --- see Figure~\ref{fig1}). Let $k\in\Z$ and $j\ge 1$ be the indices induced by the indexation map $\iota$ defined in~\eqref{indexbij} such that $(j,k)=\iota(2m,2n)$. Thus, $NW_{\bm{\lambda}}(j,k)=H_{\bm{\lambda}}(2m,2n)$.\\

\noindent  It is a consequence of  Corollary~\ref{keyobstr} that $H_{\bm{\lambda}}(m,n)= G_{\bm{\lambda}}(m, n-1) = 0$.
Therefore, from Corollary~\ref{labelhenkelbis},
\begin{itemize}
\item identity (W) implies that $NW_{\bm{\lambda}}(j,k+1) = H_{\bm{\lambda}}(2m+1,2n) = 0$
and also that $NW_{\bm{\lambda}}(j,k-1)= H_{\bm{\lambda}}(2m-1,2n)= 0$;
\item identity (X) implies that $NW_{\bm{\lambda}}(j-1,k)= H_{\bm{\lambda}}(2m+1,2n-1)= 0$
and also that $NW_{\bm{\lambda}}(j+1,k)= H_{\bm{\lambda}}(2m-1,2n+1)= 0$.
\end{itemize}
The Square Window Lemma~\ref{sqwinthe} then yields that all entries in the $3\times 3$ square centered at $(j,k)$  vanish, which establishes the claim.\\

\noindent The unevenness of the sizes of the windows easily follows. To see this, it is enough to surround each entry corresponding to a Hankel determinant of the form   $H_{\bm{\lambda}}(2m,2n)$ in a given window by a $3\times 3$ zero square.
\end{proof}

\subsection{Determinant Calculations}\label{detsubsec}

\begin{proof}[Proof of Proposition~\ref{proprec}] Each of the four identities is established following the same three steps applied to the definition~\eqref{detqmnt} of the polynomial approximation $S_{\bm{\lambda}}(m,n;t)$~:
\begin{enumerate}
\item group columns and rows with the same parity together. From the self-similarity property~\eqref{reclambda} of the paperfolding sequence, this induces the presence of the $\pm 1$ matrix $\bm{E}(n,p)$ defined in~\eqref{defEmatr} as a submatrix in the (matrix defining the) determinant under con\-si\-de\-ration;
\item introduce auxiliary rows and columns which do not change the value of the determinant. By performing elementary row and column operations with them and with the submatrix $\bm{E}(n,p)$, this creates many zero entries in the (matrix defining the) determinant;
\item a Laplace expansion along chosen columns then shows that the determinant  is the sum of three terms, which yields the sought identities.
\end{enumerate}
\noindent The calculations in the first case (A) are given in full. The other cases are treated along the same lines without detailing each step as meticulously. Recall in what follows the definition of the polynomial vector $ \bm{\widehat{w}}_{n}(t)$ given in~\eqref{defespivec}.

\paragraph{$\bullet$\textbf{Case (A)}~:} let $m\in\Z$ and $n\ge 1$. From the self-similarity property~\eqref{reclambda} 
and by permutation of rows and columns,
\begin{align*}
S_{\bm{\lambda}}&(2m+1,2n;t)\;=\;
(-1)^n\cdot\begin{vmatrix}
 \mathcal{H}_{\bm{\lambda}}\left(m, n, n+1\right) & (-1)^{m+n}\cdot\bm{E}\left(n, n\right) \vspace{2mm} \\
 (-1)^{m+n+1}\cdot \bm{E}\left(n, n+1\right) & \mathcal{H}_{\bm{\lambda}}\left(m+1, n, n\right) \vspace{2mm} \\
 \bm{w}_{n+1}\left(t^2\right) &  t\cdot \bm{w}_{n}\left(t^2\right)
\end{vmatrix}. 
\end{align*}
\noindent This determinant is left unchanged by the addition of the following two rows and columns~:
\begin{align*}
S_{\bm{\lambda}}&(2m+1,2n;t) \;=\;
(-1)^n\cdot\begin{vmatrix}
 \mathcal{H}_{\bm{\lambda}}\left(m, n, n+1\right) &  (-1)^{m+n}\cdot\bm{E}\left(n, n\right)  & \bm{0}^T\left(n\right)& \bm{0}^T\left(n\right) \vspace{2mm} \\
 (-1)^{m+n+1}\cdot \bm{E}\left(n, n+1\right)  & \mathcal{H}_{\bm{\lambda}}\left(m+1, n, n\right)  & \bm{0}^T\left(n\right)& \bm{0}^T\left(n\right)  \vspace{2mm} \\
 \bm{w}_{n+1}\left(t^2\right) &  t\cdot \bm{w}_{n}\left(t^2\right) & 0 & 0 \vspace{2mm}\\
 \bm{0}\left(n+1\right)& \bm{\varepsilon}\left(n\right)  &1&  0\vspace{2mm}\\
   \bm{\varepsilon}\left(n+1\right)  & \bm{0}\left(n\right) &0&  1
\end{vmatrix}.
\end{align*}
\noindent Upon successively adding the last row and its opposite to the ones determined by the matrix \mbox{$ (-1)^{m+n+1}\cdot \bm{E}\left(n, n+1\right) $}, and upon successively adding the second last row and its opposite to the ones determined by the matrix $ (-1)^{m+n}\cdot \bm{E}\left(n, n\right) $, this yields that $S_{\bm{\lambda}}(2m+1,2n;t)$ equals
\begin{align*}
(-1)^n\cdot\begin{vmatrix}
 \mathcal{H}_{\bm{\lambda}}\left(m, n, n+1\right) &  \bm{0}\left(n, n\right)  &  (-1)^{m}\cdot \bm{\varepsilon}^T\left(n\right)&  \bm{0}^T\left(n\right) \vspace{2mm} \\
  \bm{0}\left(n, n+1\right)  & \mathcal{H}_{\bm{\lambda}}\left(m+1, n, n\right)  &    \bm{0}^T\left(n\right)     &  (-1)^{m+1}\cdot \bm{\varepsilon}^T\left(n\right)  \vspace{2mm} \\
 \bm{w}_{n+1}\left(t^2\right) &  t\cdot \bm{w}_{n}\left(t^2\right) & 0 & 0 \vspace{2mm}\\
  \bm{0}\left(n+1\right)  &  \bm{\varepsilon}\left(n\right)   &1&  0\vspace{2mm}\\
\bm{\varepsilon}\left(n+1\right) &  \bm{0}\left(n\right)&0&  1
\end{vmatrix}.
\end{align*}
\noindent Apply now a Laplace expansion along $n+2$ columns, namely the first $n+1$ determined by $ \mathcal{H}_{\bm{\gamma}}\left(m, n, n+1\right)$ and also the second last containing $(-1)^{m}\cdot \bm{\varepsilon}^T\left(n\right)$ in its first rows. This expresses  $S_{\bm{\lambda}}(2m+1,2n;t)$ as the sum of three nonvanishing terms. The first one is
\begin{align*}
-
&\begin{vmatrix}
\mathcal{H}_{\bm{\lambda}}\left(m, n, n+1\right)  &  (-1)^{m}\cdot \bm{\varepsilon}^T\left(n\right)  \vspace{2mm} \\
\bm{w}_{n+1}\left(t^2\right) & 0  \vspace{2mm} \\
\bm{\varepsilon}\left(n+1\right) & 0
\end{vmatrix}
\cdotp
\begin{vmatrix}
\mathcal{H}_{\bm{\lambda}}\left(m+1, n, n\right)  &  (-1)^{m+1}\cdot \bm{\varepsilon}^T\left(n\right) \vspace{2mm} \\
\bm{\varepsilon}\left(n\right) & 0
\end{vmatrix} \vspace{2mm} \\
&\underset{\eqref{deterformFG} \& \eqref{detg}}{=}\; 
-\left(t^2+1\right)\cdot V_{\bm{\lambda}}\left(m, n-1; t^2\right)\cdot G_{\bm{\lambda}}\left(m+1,n-1\right);
\end{align*}
the second one is
\begin{align*}
&\begin{vmatrix}
\mathcal{H}_{\bm{\lambda}}\left(m, n, n+1\right)  &  (-1)^{m}\cdot \bm{\varepsilon}^T\left(n\right)  \vspace{2mm} \\
\bm{w}_{n+1}\left(t^2\right) & 0  \vspace{2mm} \\
\bm{0}\left(n+1\right) & 1
\end{vmatrix}
\cdot
\begin{vmatrix}
\mathcal{H}_{\bm{\lambda}}\left(m+1, n, n\right)  &  (-1)^{m+1}\cdot \bm{\varepsilon}^T\left(n\right)  \vspace{2mm} \\
\bm{0}\left(n\right) & 1
\end{vmatrix}\\
& \underset{\eqref{detqmnt}}{=}\; 
S_{\bm{\lambda}}\left(m,n;t^2\right)\cdot H_{\bm{\lambda}}\left(m+1,n\right);
\end{align*}
and the third one  is
\begin{align*}
&\begin{vmatrix}
\mathcal{H}_{\bm{\lambda}}\left(m, n, n+1\right)  &  (-1)^{m}\cdot \bm{\varepsilon}^T\left(n\right)  \vspace{2mm} \\
\bm{0}\left(n+1\right) & 1 \vspace{2mm} \\
\bm{\varepsilon}\left(n+1\right) & 0
\end{vmatrix}
\cdot
\begin{vmatrix}
\mathcal{H}_{\bm{\lambda}}\left(m+1, n, n\right)  &  (-1)^{m+1}\cdot \bm{\varepsilon}^T\left(n\right)  \vspace{2mm} \\
t\cdot\bm{w}_n\left(t^2\right) & 0
\end{vmatrix}\\
& \underset{\eqref{deterformFG} \& \eqref{detf}}{=}\; 
(-1)^{m+1}\cdot t\cdot R_{\bm{\lambda}}\left(m+1,n-1;t^2\right)\cdot F_{\bm{\lambda}}\left(m,n\right).
\end{align*}

\noindent This concludes the proof in this case.


\paragraph{$\bullet$\textbf{Case (B)}~:} let $m\in\Z$ and $n\ge 2$. Proceeding as above, one obtains
{\allowdisplaybreaks
\begin{align*}
S_{\bm{\lambda}}(2m+1, & 2n-1;t)\;\underset{\eqref{reclambda}}{=}\;
\begin{vmatrix}
 \mathcal{H}_{\bm{\gamma}}\left(m, n, n\right) & (-1)^{m+n}\cdot\bm{E}\left(n, n\right) \vspace{2mm} \\
 (-1)^{m+n}\cdot \bm{E}\left(n-1, n\right) & \mathcal{H}_{\bm{\gamma}}\left(m+1, n-1, n\right) \vspace{2mm} \\
 \bm{w}_{n}\left(t^2\right) &  t\cdot \bm{w}_{n}\left(t^2\right)
\end{vmatrix} \vspace{3mm} \\
&\;=\;
\begin{vmatrix}
 \mathcal{H}_{\bm{\gamma}}\left(m, n, n\right) &  (-1)^{m+n}\cdot\bm{E}\left(n, n\right)  & \bm{0}^T\left(n\right)& \bm{0}^T\left(n\right) \vspace{2mm} \\
 (-1)^{m+n}\cdot \bm{E}\left(n-1, n\right)  & \mathcal{H}_{\bm{\gamma}}\left(m+1, n-1, n\right)  & \bm{0}^T\left(n-1\right)& \bm{0}^T\left(n-1\right)  \vspace{2mm} \\
 \bm{w}_{n}\left(t^2\right) &  t\cdot \bm{w}_{n}\left(t^2\right) & 0 & 0 \vspace{2mm}\\
 \bm{0}\left(n\right)  &  \bm{\varepsilon}\left(n\right) &1&  0\vspace{2mm}\\
\bm{\varepsilon}\left(n\right)  & \bm{0}\left(n\right)  &0&  1
\end{vmatrix}.
\end{align*}
}
By elementary row and column operations, this yields that $S_{\bm{\lambda}}(2m+1,2n-1;t)$ equals
\begin{align*}
\begin{vmatrix}
 \mathcal{H}_{\bm{\gamma}}\left(m, n, n\right) &  \bm{0}\left(n, n\right)  &    (-1)^{m}\cdot \bm{\varepsilon}^T\left(n\right)     &   \bm{0}^T\left(n\right)\vspace{2mm} \\
  \bm{0}\left(n-1, n\right)  & \mathcal{H}_{\bm{\gamma}}\left(m+1, n-1, n\right)  &  \bm{0}^T\left(n-1\right) & (-1)^{m+1}\cdot \bm{\varepsilon}^T\left(n-1\right)       \vspace{2mm} \\
 \bm{w}_{n}\left(t^2\right) &  t\cdot \bm{w}_{n}\left(t^2\right) & 0 & 0 \vspace{2mm}\\
 \bm{0}\left(n\right)   & \bm{\varepsilon}\left(n\right)  &1&  0\vspace{2mm}\\
 \bm{\varepsilon}\left(n\right)& \bm{0}\left(n\right)  &0&  1
\end{vmatrix}.
\end{align*}
Apply the Laplace expansion along $n+1$ columns, namely  the first $n$ determined by $ \mathcal{H}_{\bm{\gamma}}\left(m, n, n\right)$ and the second last containing the vector $(-1)^{m}\cdot \bm{\varepsilon}^T\left(n\right)$ in its top rows. This expresses $S_{\bm{\lambda}}(2m+1,2n-1;t)$ as the sum of three nonvanishing terms. The first one is
\begin{align*}
&-\begin{vmatrix}
\mathcal{H}_{\bm{\gamma}}\left(m, n, n\right)  &  (-1)^{m}\cdot \bm{\varepsilon}^T\left(n\right)  \vspace{2mm} \\
\bm{w}_{n}\left(t^2\right) & 0  \vspace{2mm}
\end{vmatrix}
\cdotp
\begin{vmatrix}
\mathcal{H}_{\bm{\gamma}}\left(m+1, n-1, n\right)  &  (-1)^{m+1}\cdot \bm{\varepsilon}^T\left(n-1\right) \vspace{2mm} \\
\bm{\varepsilon}\left(n\right) & 0\vspace{2mm}\\
\bm{0}\left(n\right) & 1
\end{vmatrix} \vspace{2mm} \\
&\underset{\eqref{deterformFG} \& \eqref{detf}}{=}\; (-1)^{m}
\cdot F_{\bm{\lambda}}\left(m+1,n-1\right) \cdot  R_{\bm{\lambda}}\left(m, n-1; t^2\right);
\end{align*}
the second one is
\begin{align*}
&
\begin{vmatrix}
\mathcal{H}_{\bm{\gamma}}\left(m, n, n\right)  &  (-1)^{m}\cdot \bm{\varepsilon}^T\left(n\right)  \vspace{2mm} \\
\bm{0}\left(n\right) & 1  \vspace{2mm}
\end{vmatrix}
\cdot
\begin{vmatrix}
\mathcal{H}_{\bm{\gamma}}\left(m+1, n-1, n\right)  &  (-1)^{m+1}\cdot \bm{\varepsilon}^T\left(n-1\right) \vspace{2mm} \\
t\cdot\bm{w}_n\left(t^2\right) & 0\vspace{2mm}\\
\bm{0}\left(n\right) & 1
\end{vmatrix} \vspace{2mm} \\
& \underset{\eqref{detqmnt}}{=}\; t \cdot H_{\bm{\lambda}}\left(m,n\right)\cdot S_{\bm{\lambda}}\left(m+1, n-1; t^2\right);
\end{align*}
and the third one is
\begin{align*}
& -\begin{vmatrix}
\mathcal{H}_{\bm{\lambda}}\left(m, n, n\right)  &  (-1)^{m}\cdot \bm{\varepsilon}^T\left(n\right)  \vspace{2mm} \\
\bm{\varepsilon}(n) & 0
\end{vmatrix}
\cdotp
\begin{vmatrix}
\mathcal{H}_{\bm{\lambda}}\left(m+1, n-1, n\right)  &  (-1)^{m+1}\cdot \bm{\varepsilon}^T\left(n-1\right) \vspace{2mm} \\
t\cdot \bm{w}_n\left(t^2\right) & 0\vspace{2mm}\\
\bm{\varepsilon}\left(n\right) & 0
\end{vmatrix} \vspace{2mm}\\
&\underset{\eqref{deterformFG} \& \eqref{detg}}{=}\;  -t\cdot \left(t^2+1\right)\cdot G_{\bm{\lambda}}\left(m, n-1\right)\cdot V_{\bm{\lambda}}\left(m+1, n-2; t^2\right).
\end{align*}
This concludes the proof in this case.


\paragraph{$\bullet$\textbf{Case (C)}~:} let $m\in\Z$ and $n\ge 1$. Proceeding along the same way, one obtains that
{\allowdisplaybreaks
\begin{align*}
S_{\bm{\lambda}}(2m,  2n;t)&\;=\;
(-1)^n\cdot\begin{vmatrix}
(-1)^{m+n}\cdot \bm{E}\left(n, n+1\right) & \mathcal{H}_{\bm{\lambda}}\left(m, n, n\right)\vspace{2mm} \\
\mathcal{H}_{\bm{\lambda}}\left(m, n, n+1\right) & (-1)^{m+n}\cdot \bm{E}\left(n, n\right) \vspace{2mm} \\
 \bm{w}_{n+1}\left(t^2\right) &  t\cdot \bm{w}_{n}\left(t^2\right)
\end{vmatrix} \vspace{3mm} \\
&\;=\;
(-1)^n\cdot\begin{vmatrix}
(-1)^{m+n}\cdot \bm{E}\left(n, n+1\right) & \mathcal{H}_{\bm{\gamma}}\left(m, n, n\right) & \bm{0}^T\left(n\right)& \bm{0}^T\left(n\right) \vspace{2mm} \\
\mathcal{H}_{\bm{\gamma}}\left(m, n, n+1\right) & (-1)^{m+n}\cdot \bm{E}\left(n, n\right) & \bm{0}^T\left(n\right)& \bm{0}^T\left(n\right)  \vspace{2mm} \\
 \bm{w}_{n+1}\left(t^2\right) &  t\cdot \bm{w}_{n}\left(t^2\right) & 0 & 0 \vspace{2mm}\\
\bm{\varepsilon}\left(n+1\right) & \bm{0}\left(n\right)&1&  0\vspace{2mm}\\
  \bm{0}\left(n+1\right) & \bm{\varepsilon}\left(n\right) &0&  1
\end{vmatrix}.
\end{align*}
}
By elementary row and column operations, this yields that $S_{\bm{\lambda}}(2m,2n;t)$ equals
\begin{align*}
(-1)^n\cdot\begin{vmatrix}
\bm{0}\left(n, n+1\right) & \mathcal{H}_{\bm{\gamma}}\left(m, n, n\right) & (-1)^{m}\cdot \bm{\varepsilon}^T\left(n\right)& \bm{0}^T\left(n\right) \vspace{2mm} \\
\mathcal{H}_{\bm{\gamma}}\left(m, n, n+1\right) &  \bm{0}\left(n, n\right) & \bm{0}^T\left(n\right)& (-1)^{m}\cdot \bm{\varepsilon}^T\left(n\right) \vspace{2mm} \\
 \bm{w}_{n+1}\left(t^2\right) &  t\cdot \bm{w}_{n}\left(t^2\right) & 0 & 0 \vspace{2mm}\\
\bm{\varepsilon}\left(n+1\right) & \bm{0}\left(n\right)&1&  0\vspace{2mm}\\
  \bm{0}\left(n+1\right) & \bm{\varepsilon}\left(n\right) &0&  1
\end{vmatrix}.
\end{align*}
\noindent Apply a Laplace expansion along $n+2$ columns, namely the first $n+1$ determined by $ \mathcal{H}_{\bm{\gamma}}\left(m, n, n+1\right)$ and the last one determined by $(-1)^{m}\cdot \bm{\varepsilon}^T\left(n\right)$. This expresses $S_{\bm{\lambda}}(2m, 2n;t)$ as the sum of three nonvanishing terms. The first one is
\begin{align*}
(-1)^{n+1}\cdot&\begin{vmatrix}
\mathcal{H}_{\bm{\gamma}}\left(m, n, n+1\right)  &  (-1)^{m}\cdot \bm{\varepsilon}^T\left(n\right)  \vspace{2mm} \\
\bm{w}_{n+1}\left(t^2\right) & 0  \vspace{2mm} \\
\bm{\varepsilon}\left(n+1\right) & 0
\end{vmatrix}
\cdotp
\begin{vmatrix}
\mathcal{H}_{\bm{\gamma}}\left(m, n, n\right)  &  (-1)^{m}\cdot \bm{\varepsilon}^T\left(n\right) \vspace{2mm} \\
\bm{\varepsilon}\left(n\right) & 0
\end{vmatrix} \vspace{2mm} \\
&\underset{\eqref{deterformFG} \& \eqref{detg}}{=}\; (-1)^n\cdot\left(t^2+1\right)\cdot V_{\bm{\lambda}}\left(m, n-1; t^2\right)\cdot G_{\bm{\lambda}}\left(m,n-1\right);
\end{align*}
the second one is
\begin{align*}
(-1)^{n+1}\cdot&\begin{vmatrix}
\mathcal{H}_{\bm{\gamma}}\left(m, n, n+1\right)  &  (-1)^{m}\cdot \bm{\varepsilon}^T\left(n\right)  \vspace{2mm} \\
\bm{\varepsilon}\left(n+1\right) & 0  \vspace{2mm} \\
\bm{0}\left(n+1\right) & 1
\end{vmatrix}
\cdot
\begin{vmatrix}
\mathcal{H}_{\bm{\gamma}}\left(m, n, n\right)  &  (-1)^{m}\cdot \bm{\varepsilon}^T\left(n\right)  \vspace{2mm} \\
t\cdot\bm{w}_n\left(t^2\right) & 0  
\end{vmatrix}\\
& \underset{\eqref{deterformFG} \& \eqref{detf}}{=}\; \left(-1\right)^{m+n} \cdot t \cdot F_{\bm{\lambda}}\left(m,n\right)\cdot R_{\bm{\lambda}}\left(m, n-1; t^2\right);
\end{align*}
and  the third one is
\begin{align*}
(-1)^n\cdot&\begin{vmatrix}
\mathcal{H}_{\bm{\lambda}}\left(m, n, n+1\right)  &  (-1)^{m}\cdot \bm{\varepsilon}^T\left(n\right)  \vspace{2mm} \\
\bm{w}_{n+1}\left(t^2\right) & 0  \vspace{2mm} \\
\bm{0}\left(n+1\right) & 1
\end{vmatrix}
\cdotp
\begin{vmatrix}
\mathcal{H}_{\bm{\lambda}}\left(m, n, n\right)  &  (-1)^{m}\cdot \bm{\varepsilon}^T\left(n\right) \vspace{2mm} \\
\bm{0}\left(n\right) & 1
\end{vmatrix} \vspace{2mm} \\
&\underset{\eqref{detqmnt}}{=}\;    (-1)^n\cdot S_{\bm{\lambda}}\left(m,n;t^2\right) \cdot H_{\bm{\lambda}}\left(m, n\right).
\end{align*}
This concludes the proof in this case.


\paragraph{$\bullet$\textbf{Case (D)}~:} let $m\in\Z$ and $n\ge 2$. Proceeding along the same way as before in this final case also, one obtains that
{\allowdisplaybreaks
\begin{align*}
S_{\bm{\gamma}}(2m,&2n-1;t)\; =\;
\begin{vmatrix}
(-1)^{m+n+1}\cdot \bm{E}\left(n, n\right) & \mathcal{H}_{\bm{\gamma}}\left(m, n, n\right)\vspace{2mm} \\
\mathcal{H}_{\bm{\gamma}}\left(m, n-1, n\right) & (-1)^{m+n}\cdot \bm{E}\left(n-1, n\right) \vspace{2mm} \\
 \bm{w}_{n}\left(t^2\right) &  t\cdot \bm{w}_{n}\left(t^2\right)
\end{vmatrix} \vspace{3mm} \\
&=
\begin{vmatrix}
(-1)^{m+n+1}\cdot \bm{E}\left(n, n\right) & \mathcal{H}_{\bm{\gamma}}\left(m, n, n\right) & \bm{0}^T\left(n\right)& \bm{0}^T\left(n\right) \vspace{2mm} \\
\mathcal{H}_{\bm{\gamma}}\left(m, n-1, n\right) & (-1)^{m+n}\cdot \bm{E}\left(n-1, n\right) & \bm{0}^T\left(n-1\right)& \bm{0}^T\left(n-1\right)  \vspace{2mm} \\
 \bm{w}_{n}\left(t^2\right) &  t\cdot \bm{w}_{n}\left(t^2\right) & 0 & 0 \vspace{2mm}\\
\bm{\varepsilon}\left(n\right) & \bm{0}\left(n\right)&1&  0\vspace{2mm}\\
  \bm{0}\left(n\right) & \bm{\varepsilon}\left(n\right) &0&  1
\end{vmatrix}.
\end{align*}
}
By elementary row and column operations, this yields that $S_{\bm{\lambda}}(2m,2n-1;t)$ equals
\begin{align*}
\begin{vmatrix}
\bm{0}\left(n, n\right)& \mathcal{H}_{\bm{\gamma}}\left(m, n, n\right) & (-1)^{m+1}\cdot \bm{\varepsilon}^T\left(n\right)& \bm{0}^T\left(n\right) \vspace{2mm} \\
\mathcal{H}_{\bm{\gamma}}\left(m, n-1, n\right) &   \bm{0}\left(n-1, n\right) & \bm{0}^T\left(n-1\right)& (-1)^{m+1}\cdot \bm{\varepsilon}^T\left(n-1\right) \vspace{2mm} \\
 \bm{w}_{n}\left(t^2\right) &  t\cdot \bm{w}_{n}\left(t^2\right) & 0 & 0 \vspace{2mm}\\
\bm{\varepsilon}\left(n\right) & \bm{0}\left(n\right)&1&  0\vspace{2mm}\\
  \bm{0}\left(n\right) & \bm{\varepsilon}\left(n\right) &0&  1
\end{vmatrix}.
\end{align*}
Apply a Laplace expansion along $n+1$ columns, namely the first $n$ determined by $ \mathcal{H}_{\bm{\gamma}}\left(m, n-1, n\right)$ and the last one determined by $(-1)^{m+1}\cdot \bm{\varepsilon}^T\left(n-1\right)$. This expresses $S_{\bm{\lambda}}(2m, 2n-1;t)$ as the sum of three nonvanishing terms. The first one is
\begin{align*}
\left(-1\right)^{n+1}\cdot
&\begin{vmatrix}
\mathcal{H}_{\bm{\gamma}}\left(m, n-1, n\right)  &  (-1)^{m+1}\cdot \bm{\varepsilon}^T\left(n-1\right)  \vspace{2mm} \\
\bm{w}_{n}\left(t^2\right) & 0  \vspace{2mm} \\
\bm{\varepsilon}\left(n\right) & 0
\end{vmatrix}
\cdotp
\begin{vmatrix}
\mathcal{H}_{\bm{\gamma}}\left(m, n, n\right)  &  (-1)^{m+1}\cdot \bm{\varepsilon}^T\left(n\right) \vspace{2mm} \\
\bm{\varepsilon}\left(n\right) & 0
\end{vmatrix} \vspace{2mm} \\
&\underset{\eqref{deterformFG} \& \eqref{detg}}{=}\; \left(-1\right)^{n} \cdot \left(t^2+1\right)\cdot V_{\bm{\gamma}}\left(m, n-2; t^2\right)\cdot G_{\bm{\gamma}}\left(m,n-1\right);
\end{align*}
the second one is
\begin{align*}
\left(-1\right)^{n}\cdot
&\begin{vmatrix}
\mathcal{H}_{\bm{\gamma}}\left(m, n-1, n\right)  &  (-1)^{m+1}\cdot \bm{\varepsilon}^T\left(n-1\right)  \vspace{2mm} \\
\bm{w}_{n}\left(t^2\right) & 0  \vspace{2mm} \\
\bm{0}\left(n\right) & 1
\end{vmatrix}
\cdotp
\begin{vmatrix}
\mathcal{H}_{\bm{\gamma}}\left(m, n, n\right)  &  (-1)^{m+1}\cdot \bm{\varepsilon}^T\left(n\right) \vspace{2mm} \\
\bm{0}\left(n\right) & 1
\end{vmatrix} \vspace{2mm} \\
&\underset{\eqref{detqmnt}}{=}\; \left(-1\right)^{n}\cdot S_{\bm{\lambda}}\left(m, n-1; t^2\right) \cdot H_{\bm{\lambda}}\left(m,n\right);
\end{align*}
and the third one is
\begin{align*}
\left(-1\right)^{n+1}\cdot
&\begin{vmatrix}
\mathcal{H}_{\bm{\lambda}}\left(m, n-1, n\right)  &  (-1)^{m+1}\cdot \bm{\varepsilon}^T\left(n-1\right)  \vspace{2mm} \\
\bm{\varepsilon}\left(n\right) & 0  \vspace{2mm} \\
\bm{0}\left(n\right) & 1
\end{vmatrix}
\cdotp
\begin{vmatrix}
\mathcal{H}_{\bm{\gamma}}\left(m, n, n\right)  &  (-1)^{m+1}\cdot \bm{\varepsilon}^T\left(n\right) \vspace{2mm} \\
t\cdot\bm{w}_n\left(t^2\right) & 0
\end{vmatrix} \vspace{2mm} \\
&\underset{\eqref{deterformFG}\&\eqref{detf}}{=}\;  (-1)^{n+m+1}\cdot F_{\bm{\lambda}}\left(m, n-1\right)\cdot t\cdot R_{\bm{\lambda}}\left(m, n-1; t^2\right).
\end{align*}
This concludes the proof of Case (D) and thus  of Proposition~\ref{proprec} also.

\end{proof}

\section{Arithmetic Properties  of the Convergents of the Paperfolding Sequence}
\label{sec4}

\subsection{Convergent of a Window}\label{cvgtwinsbsc}

\noindent Consider the doubly infinite sequence $\bm{\gamma}$  defined in~\eqref{defdoublyextendedeq} over a given field $\K$ and let $NW_{\bm{\gamma}}$  be its Number Wall over this field. A window $W$ contained in it has a size denoted by $\rho(W)\ge 1$ and  is identified with its top--left corner. The coordinates of this top--left corner shall  be assumed to be of the form  $(j(W)+1,k(W)+1)$, where $j(W)\ge 1$ and $k(W)\in\Z$ are integers\textsuperscript{8}\let\thefootnote\relax\footnotetext{\textsuperscript{8}This notation is adopted with a view towards specialising $\bm{\gamma}$ to the paperfolding sequence over the alphabet $\pm 1$~: then, the first row of the Number Wall (obtained when $j=1$) contains the values taken by the sequence and can thus not be the location of the corner of a window.}. Set
\begin{equation}\label{defnmW}
\left(m(W), \; n(W)\right)\;=\; \iota^{-1}\left(j(W), k(W)\right)  \;\underset{\eqref{indexbij}}{=}\;  \left(k(W)-j(W), \; j(W)\right),
\end{equation}
where $\iota$ is the indexation map defined in~\eqref{indexbij} (thus, $m(W)\in\Z$ and $n(W)\ge 1$). The denominator and numerator of the \emph{convergent associated to the window $W$} are then defined as
\begin{equation*}\label{convpoly}
Q_{\bm{\gamma}}\left(W; t\right)\;=\; Q_{\bm{\gamma}}\left(m(W), \; h(W); t\right)\qquad \textrm{and}\qquad P_{\bm{\gamma}}\left(W; t\right)\;=\; P_{\bm{\gamma}}\left(m(W), \; h(W); t\right).
\end{equation*}
Here,  $h(W)\ge 0$ is the unique integer for which
\begin{equation}\label{defqlamnw}
q_{\bm{\gamma}}(W)\;:=\;\deg Q_{\bm{\gamma}}\left(W;t\right) \;=\; j(W) \;\underset{\eqref{defnmW}}{=}\; n\left(W\right).
\end{equation}
These quantities are well--defined from Point~(3) in Lemma~\ref{closedromlem} since  $ H_{\bm{\gamma}}\left(m(W), \; n(W)\right) = NW_{\bm{\gamma}}(j(W),k(W))\neq 0$.
(In fact, they can also be defined  when $j(W)=0$ upon taking $Q_{\bm{\gamma}}\left(W;t\right)$ to be  the constant polynomial $Q_{\bm{\gamma}}\left(m(W), 0; t\right)$ obtained when $h(W)=0$, and thus upon also setting $P_{\bm{\gamma}}\left(W;t\right)=P_{\bm{\gamma}}\left(m(W), 0; t\right)$.)

\begin{prop}[Properties of the convergent of a window]\label{propconvwind}
Keep the above notations and terminology. Then, the convergent of a window $W$ with size $\rho(W)\ge 1$ and top--left corner at the entry $(j(W)+1, k(W)+1)$, where $j(W)\ge 1$, satisfies the relation
\begin{equation}\label{approxcovpolywin}
\left| Q_{\bm{\gamma}}\left(W;t\right) \cdot   \Gamma_{m(W)}(t)-P_{\bm{\gamma}}\left(W;t\right)\right|\;=\; 2^{-q_{\bm{\gamma}}(W)-\rho(W)-1}.
\end{equation}
Furthermore, neither the polynomial $P_{\bm{\gamma}}\left(W;t\right)+\gamma_{\left(m(W)-1\right)}\cdot Q_{\bm{\gamma}}\left(W;t\right)$ nor the polynomial $Q_{\bm{\gamma}}\left(W;t\right)$ is divisible by  $t$.
\end{prop}

\noindent The proof of the statement is based on the following equation valid for all integers $m\in\Z$~:
\begin{equation}\label{shiftgammam}
\Gamma_m(t)\;=\; t\cdot\Gamma_{m-1}(t)-\gamma_{m-1}.
\end{equation}
It is readily derived from the definition of the Laurent series $\Gamma_m(t)$ in~\eqref{doublylaurent}.

\begin{proof}
\noindent Equation~\eqref{approxcovpolywin} is an immediate consequence of Point~(3) in Lemma~\ref{closedromlem}.

\paragraph{$\bullet$}  To show that $Q_{\bm{\gamma}}\left(W;t\right)$ cannot be divisible by $t$, argue by contradiction and decompose it as $$Q_{\bm{\gamma}}\left(W; t\right) = t\cdot \widetilde{Q}_{\bm{\gamma}}\left(W; t\right),$$ where  $\widetilde{Q}_{\bm{\gamma}}\left(W; t\right)$ is a nonzero polynomial. Then,
\begin{align}
&\left| Q_{\bm{\gamma}}\left(W;t\right) \cdot   \Gamma_{m(W)}(t)-P_{\bm{\gamma}}\left(W;t\right)\right| \;=\; \left|  \widetilde{Q}_{\bm{\gamma}}\left(W; t\right)\cdot  t\cdot \Gamma_{m(W)}(t)-P_{\bm{\gamma}}\left(W;t\right)\right|  \nonumber\\
&\qquad\qquad \qquad\underset{\eqref{shiftgammam}}{=}\; \left| \widetilde{Q}_{\bm{\gamma}}\left(W; t\right)\cdot   \Gamma_{(m(W)+1)}(t)-\left(P_{\bm{\gamma}}\left(W;t\right)-\gamma_{m(W)}\cdot \widetilde{Q}_{\bm{\gamma}}\left(W; t\right)\right)\right| \nonumber\\
&\qquad\qquad \qquad\underset{\eqref{approxcovpolywin}}{=}\;  2^{-q_{\bm{\gamma}}(W)-\rho(W)-1}. \label{expqgammaW}
\end{align}
If the polynomials $\widetilde{Q}_{\bm{\gamma}}\left(W; t\right)$ is constant (meaning that $q_{\bm{\gamma}}(W)=1= j(W)=n(W)$), then the definition of the norm of a Laurent series forces the coefficients of $t^{-1}$ 
in $ \Gamma_{(m(W)+1)}(t)$ to vanish. This coefficient is  $$\gamma_{(m(W)+1)}\;=\;H_{\bm{\gamma}}(m(W)+1, 1) \; =\; NW_{\bm{\gamma}}(1, m(W)+2)$$ (see Figure~\ref{fig1} for the indexation). 
From the Square Window Lemma~\ref{sqwinthe}, 
this nevertheless contradicts the assumption that the entry with coordinates $\left(j(W)+1, k(W)+1\right) = (2, m(W)+2)$ is the top--left corner of the window $W$.\\

\noindent Assume therefore from now on that,  upon setting $\widetilde{q}_{\bm{\gamma}}(W)=\deg \widetilde{Q}_{\bm{\gamma}}\left(W;t\right)$, it holds that
\begin{equation}\label{qWqw'}
n(W)\;\underset{\eqref{defqlamnw}}{=}\;q_{\bm{\gamma}}(W)\; =\; \widetilde{q}_{\bm{\gamma}}(W)+1\;\ge\; 2.
\end{equation}
In particular,
\begin{align*}
\left| \widetilde{Q}_{\bm{\gamma}}\left(W; t\right)\cdot   \Gamma_{(m(W)+1)}(t)-\left(P_{\bm{\gamma}}\left(W;t\right)-\gamma_{m(W)}\cdot \widetilde{Q}_{\bm{\gamma}}\left(W; t\right)\right)\right|\;&\underset{\eqref{expqgammaW}}{=}\; 2^{-\widetilde{q}_{\bm{\gamma}}(W)-\rho(W)-2}\\
&= \; \frac{1}{2^{\rho(W)+2}\cdot\left|\widetilde{Q}_{\bm{\gamma}}\left(W; t\right)\right|}\cdotp
\end{align*}

\noindent From the equivalence~\eqref{convcharac}, $\widetilde{Q}_{\bm{\gamma}}\left(W; t\right)$ is the denominator of a convergent to the Laurent series $ \Gamma_{(m(W)+1)}(t)$. One then infers from
Point~(3) in Lemma~\ref{closedromlem} 
 that $$H_{\bm{\gamma}}\left(m(W)+1, \widetilde{q}_{\bm{\gamma}}(W)\right)\;\underset{\eqref{qWqw'}}{=}\;H_{\bm{\gamma}}\left(m(W)+1, \;n(W)-1\right)\;\neq\; 0$$
and that $$H_{\bm{\gamma}}\left(m(W)+1, \widetilde{q}_{\bm{\gamma}}(W)+\rho(W)+2\right)\underset{\eqref{qWqw'}}{=} H_{\bm{\gamma}}\left(m(W)+1,\; n(W)+\rho(W)+1\right)\neq 0$$ but that for all  $1\le s\le \rho(W)+1$,
\begin{equation}\label{hlmqlws}
H_{\bm{\gamma}}\left(m(W)+1, \widetilde{q}_{\bm{\gamma}}(W)+s\right)\;\underset{\eqref{qWqw'}}{=}\;H_{\bm{\gamma}}\left(m(W)+1, n(W)+s-1\right)\;=\; 0.
\end{equation}
Recall that  $\left(j(W), k(W)\right)$ are the coordinates associated to the pair $\left(m(W), n(W)\right)$ according to the indexation map $\iota$ defined in~\eqref{indexbij}. Then, the pair $\left(m(W)+1, n(W))\right)$ is associated to the coordinates $\left(j(W), k(W)+1\right)$ (see Figure~\ref{fig1}). From the above relation~\eqref{hlmqlws} evaluated at $s=1$, this is a zero entry, which contradicts the assumption that $\left(j(W)+1, k(W)+1\right)$  is a corner of the window $W$. This shows that the polynomial $Q_{\bm{\gamma}}\left(W;t\right)$ cannot be divisible by $t$.

\paragraph{$\bullet$}  The argument for the polynomial  $D_{\bm{\gamma}}(W; t):=P_{\bm{\gamma}}\left(W;t\right)+\gamma_{\left(m(W)-1\right)}\cdot Q(W;t)$ runs si\-mi\-larly~: argue by contradiction and define the polynomial $\widetilde{D}_{\bm{\gamma}}(W; t)$ by the di\-vi\-si\-bi\-lity relation $D_{\bm{\gamma}}(W; t)=t\cdot \widetilde{D}_{\bm{\gamma}}(W; t)$. Then,
\begin{align*}
\left| Q_{\bm{\gamma}}\left(W;t\right) \cdot   \Gamma_{m(W)}(t)-P_{\bm{\gamma}}\left(W;t\right)\right|\; &\underset{\eqref{shiftgammam}}{=}\;  \left| Q_{\bm{\gamma}}\left(W; t\right)\cdot t\cdot  \Gamma_{(m(W)-1)}(t)-t\cdot \widetilde{D}_{\bm{\gamma}}\left(W;t\right)\right| \nonumber\\
& \underset{\eqref{approxcovpolywin}}{=}\;  2^{-q_{\bm{\gamma}}(W)-\rho(W)-1},
\end{align*}
whence $$\left| Q_{\bm{\gamma}}\left(W; t\right)\cdot  \Gamma_{(m(W)-1)}(t)- \widetilde{D}_{\bm{\gamma}}\left(W;t\right)\right|\; = \; \frac{1}{2^{\rho(W)+2}\cdot\left|Q_{\bm{\gamma}}\left(W; t\right)\right|}\cdotp$$


\noindent The polynomial $Q_{\bm{\gamma}}\left(W; t\right)$ is therefore the denominator of a convergent of the Laurent series  $\Gamma_{(m(W)-1)}(t)$. In the same way as above, one infers from Points~(3) in Lemma~\ref{closedromlem}  that
\begin{align*}
H_{\bm{\gamma}}(m(W)-1, q_{\bm{\gamma}}(W)) \;&\underset{\eqref{defqlamnw}}{=}\; H_{\bm{\gamma}}(m(W)-1, n(W))
&\underset{ \eqref{indexbij} \& \eqref{defnmW}}{=} NW_{\bm{\gamma}}(j(W), k(W)-1)\;\neq\; 0
\end{align*}
but that
\begin{align*}
H_{\bm{\gamma}}(m(W)-1, q_{\bm{\gamma}}(W)) +1)&\underset{\eqref{defqlamnw}}{=}H_{\bm{\gamma}}(m(W)-1, n(W)+1)=NW_{\bm{\gamma}}(j(W)+1, k(W))\; =\; 0
\end{align*}
(see Figure~\ref{fig1} for the indexation). This contradicts the fact that the entry with coordinates $(j(W)+1, k(W)+1)$ is a corner of the window $W$ and thus   completes the  proof of Proposition~\ref{propconvwind}.
\end{proof}

\noindent The convergents associated to the entries in the top and in the left parts of  the  frame of a window can be derived from the convergent of the window~:

\begin{prop}[Expression for the convergents associated to the frame of a window]\label{convframe}
Let $W$ be a window of size $\rho(W)\ge 1$ with top--left corner at the entry $(j(W)+1, k(W)+1)$, where $j(W)\ge 1$ is an integer and
where $(j(W), k(W))= \iota(m(W), n(W))$ according to the indexation map $\iota$  introduced in~\eqref{indexbij}. Let $P_{\bm{\gamma}}(W; t)/Q_{\bm{\gamma}}(W; t)$ be its convergent. Then,
\begin{itemize}
\item[\textbf{(U)}] \emph{[Browsing the left of the frame]}. Given an integer $s\in\left\llbracket 0, \rho(W)\right\rrbracket$, the convergent of the Laurent series $\Gamma_{\left(m(W)-s\right)}(t)$ achieving  the normal order $n(W)+s$ and thus corresponding to the entry with coordinates
\begin{equation}\label{indexationbrowsin}
\left(j(W)+s,k(W)\right)\;=\; \iota\left(m(W)-s, n(W)+s\right)
\end{equation}
is the irreducible fraction
\begin{equation}\label{convbrowframepoi}
\frac{P_{\bm{\gamma}}(W; t)+Q_{\bm{\gamma}}(W; t)\cdot \left(\sum_{r=1}^{s}\gamma_{\left(m(W)-r\right)}\cdot t^{r-1}\right)}{t^s\cdot Q_{\bm{\gamma}}(W; t)}\cdotp
\end{equation}
\item[\textbf{(V)}] \emph{[Browsing the top of the frame]}. Given an integer $s\in\left\llbracket 0, \rho(W)\right\rrbracket$, the convergent of the Laurent series $\Gamma_{\left(m(W)+s\right)}(t)$ achieving  the normal order $n(W)$ and thus corresponding to the entry with coordinates $$\left(j(W),k(W)+s\right)\;=\;  \iota\left(m(W)+s, n(W)\right)$$  is the irreducible fraction
\begin{equation}\label{convbrowframedzdz}
\frac{t^s\cdot P_{\bm{\gamma}}(W; t)-Q_{\bm{\gamma}}(W; t)\cdot \left(\sum_{r=1}^{s}\gamma_{\left(m(W)+s-r\right)}\cdot t^{r-1}\right)}{Q_{\bm{\gamma}}(W; t)}\cdotp
\end{equation}
\end{itemize}
\noindent Conventionally, the above sums vanish when $s=0$.
\end{prop}

\begin{proof}
\noindent Each  claim is established by induction on the integer $s$.

\paragraph{$\bullet$\textbf{ Claim (U).}} The statement is clearly true when $s=0$. Assume it holds for a given index  $0\le s\le \rho(W)-1$  
and define   for the sake of simplicity $P_{\bm{\gamma}}^{(s)}(W; t)/Q_{\bm{\gamma}}^{(s)}(W; t)$ as the  fraction~\eqref{convbrowframepoi}, assumed to be irreducible by induction hypothesis. Since the window $W$ has size $\rho(W)\ge 1$ and has a top--left corner located at $(j(W)+1, k(W)+1)$, where $(j(W), k(W))= \iota(m(W), n(W))$,  it follows from Point~(3) in Lemma~\ref{closedromlem} and from the indexation~\eqref{indexationbrowsin}  that $$\left| Q_{\bm{\gamma}}^{(s)}(W; t) \cdot   \Gamma_{\left(m(W)-s\right)}(t)-P_{\bm{\gamma}}^{(s)}(W; t)\right|\;=\; 2^{-q_{\bm{\gamma}}^{(s)}(W) -(\rho(W)-s)-1},$$ where $q_{\bm{\gamma}}^{(s)}(W):= \deg Q_{\bm{\gamma}}^{(s)}(W; t) = n(W)+s.$ Upon setting
\begin{equation}\label{defrsgamw}
D_{\bm{\gamma}}^{(s)}(W; t)\;=\;P_{\bm{\gamma}}^{(s)}(W; t)+\gamma_{\left(m(W)-s-1\right)}\cdot Q_{\bm{\gamma}}^{(s)}(W; t),
\end{equation}
identity~\eqref{shiftgammam} then yields that
\begin{equation}\label{relrecconvwin}
 \left| Q_{\bm{\gamma}}^{(s)}(W; t) \cdot t\cdot   \Gamma_{\left(m(W)-s-1\right)}(t)-D_{\bm{\gamma}}^{(s)}(W; t)\right| \;=\; 2^{-q_{\bm{\gamma}}^{(s)}(W) -(\rho(W)-s)-1}.
\end{equation}
In this relation, the polynomial $D_{\bm{\gamma}}^{(s)}(W; t)$ cannot be divisible by $t$. Indeed, otherwise, upon dividing~\eqref{relrecconvwin} by $t$, the equivalence characterising the convergents of a Laurent series stated in~\eqref{convcharac} and  Point~(3) in Lemma~\ref{closedromlem} yield on the one hand that the entry $NW_{\bm{\gamma}}\left(j(W)+s, k(W)-1\right)$ with coordinates $$\left(j(W)+s, k(W)-1\right)\;=\; \iota\left(m(W)-s-1, q_{\bm{\gamma}}^{(s)}(W)\right)\;=\;\iota\left(m(W)-s-1, n(W)+s\right)$$  does not vanish. On the other, one also gets that $q_{\bm{\gamma}}^{(s)}(W)+1$ cannot be a normal order for the Laurent series $ \Gamma_{\left(m(W)-s-1\right)}(t)$; in other words that $ H_{\bm{\gamma}}\left(m(W)-s-1, n(W)+s+1\right)=NW_{\bm{\gamma}}\left(j(W)+s+1, k(W)\right)=  0$.
The Square Window Lemma~\ref{sqwinthe} then prevents the entry with coordinates $\left(j(W)+1, k(W)+1\right)$ from being the top--left corner of the window $W$. This confirms that  $t$ can indeed not divide  $D_{\bm{\gamma}}^{(s)}(W; t)$.\\

\noindent Recall now that the convergents of a Laurent series have been normalised by prescribing the leading coefficients of their denominators (see Remark~\ref{rem0}). Since $t$ does not divide $D_{\bm{\gamma}}^{(s)}(W; t)$, it is coprime with $t\cdot Q_{\bm{\gamma}}^{(s)}(W; t)$. Consequently, one gets from Equation~\eqref{relrecconvwin}
that the normal order corresponding to the indices $\left(m(W)-s-1, q_{\bm{\gamma}}^{(s)}(W)+1\right)= \left(m(W)-s-1, n(W)+s+1\right)$ is achieved by the convergent $D_{\bm{\gamma}}^{(s)}(W; t)/ (t\cdot Q_{\bm{\gamma}}^{(s)}(W; t))$. By uniqueness, one thus obtains that $$Q_{\bm{\gamma}}^{(s+1)}(W; t)\;=\; t\cdot Q_{\bm{\gamma}}^{(s)}(W; t))\qquad \quad \textrm{and}\quad\qquad P_{\bm{\gamma}}^{(s+1)}(W; t)\;\underset{\eqref{defrsgamw}}{=}\; D_{\bm{\gamma}}^{(s)}(W; t).$$ This completes the proof of the claim.

\paragraph{$\bullet$\textbf{ Claim (V).}} The argument runs similarly~: only essential differences are put forward and the details  left to the reader. \\

\noindent The statement is clearly true when $s= 0$. Assume then that it holds for a given index $0\le s\le \rho(W)-1$ and define   for the sake of simplicity $\widehat{P}_{\bm{\gamma}}^{(s)}(W; t)$  as the numerator of the fraction~\eqref{convbrowframedzdz} (assumed to be irreducible by induction hypothesis). Then,
\begin{align*}
&\left| Q_{\bm{\gamma}}(W; t) \cdot   \Gamma_{\left(m(W)+s+1\right)}(t)-\widehat{P}_{\bm{\gamma}}^{(s+1)}(W; t)\right|\\
&\qquad \qquad\qquad \underset{\eqref{shiftgammam} }{=} \left|t\cdot Q_{\bm{\gamma}}(W; t) \cdot   \Gamma_{\left(m(W)+s\right)}(t)-\left(\widehat{P}_{\bm{\gamma}}^{(s+1)}(W; t)+\gamma_{\left(m(W)+s\right)}\cdot Q_{\bm{\gamma}}(W; t) \right)\right|\\
&\qquad \qquad\qquad \underset{\eqref{convbrowframedzdz} }{=}\left|t\cdot Q_{\bm{\gamma}}(W; t) \cdot   \Gamma_{\left(m(W)+s\right)}(t)-t\cdot \widehat{P}_{\bm{\gamma}}^{(s)}(W; t)\right|
\end{align*}
whence, calling upon the induction hypothesis, 
\begin{align*}
&\left| Q_{\bm{\gamma}}(W; t) \cdot   \Gamma_{\left(m(W)+s+1\right)}(t)-\widehat{P}_{\bm{\gamma}}^{(s+1)}(W; t)\right|\;=\; \frac{1}{2^{\rho(W)-s}\cdot \left|Q_{\bm{\gamma}}(W; t) \right|}\;<\; \frac{1}{ \left|Q_{\bm{\gamma}}(W; t) \right|}\cdotp
\end{align*}
By Proposition~\ref{propconvwind},
the polynomial $Q_{\bm{\gamma}}(W; t)$ is not divisible by $t$, and is therefore coprime with $\widehat{P}_{\bm{\gamma}}^{(s+1)}(W; t)$. One thus infers  that $\widehat{P}_{\bm{\gamma}}^{(s+1)}(W; t)/Q_{\bm{\gamma}}(W; t)$ is the convergent achieving the normal order corresponding to the indices $\left(m(W)+s+1, q_{\bm{\gamma}}(W)\right)= \left(m(W)+s+1, n(W)\right)=\iota^{-1}(j(W), k(W)+s+1)$. This completes the proof of the claim and thus of the proposition.
\end{proof}

\noindent The next corollary derives a converse to Proposition~\ref{propconvwind}. Taken together, these two statements provide a necessary and sufficient condition for an irreducible fraction to be the convergent of a window in the Number Wall of a given sequence.

\begin{coro}[Sufficient condition for an irreducible fraction to be  the convergent of a window]\label{CSconvwin}
Given two coprime polynomials $P(t)$ and $Q(t)$, there exists a window $W$ in the Number Wall of the sequence $\bm{\gamma}$ such that $\left(P(t), Q(t)\right)=\left(P_{\bm{\gamma}}\left(W;t\right), Q_{\bm{\gamma}}\left(W;t\right) \right)$  if there exist integers $m\in\Z$ and $\rho\ge 1$ meeting the following two conditions~: neither the polynomial $P(t)+\gamma_{m-1}\cdot Q(t)$   nor the polynomial $Q(t)$ is divisible by $t$, and
\begin{equation}\label{condiwindow}
\left| Q(t) \cdot   \Gamma_{m}(t)-P(t)\right|\;=\; 2^{-\deg Q(t)-\rho-1}.
\end{equation}
In this case, the window $W$ is of size $\rho(W)=\rho$ and its top--left corner is located at the entry with coordinates $$(j(W)+1, k(W)+1)\; =\; (\deg Q(t)+1, m+\deg Q(t)+1).$$
\end{coro}

\begin{proof}
\noindent Let $(j,k):=(\deg Q(t)+1, m+\deg Q(t)+1)=\iota(m, \deg Q(t)+1)$ be the coordinates of the point appearing in the conclusion of the statement. Equation~\eqref{condiwindow} implies from Point~(3) in Lemma~\ref{closedromlem} that the diagonal $\left\{\left(j+s, k+s\right)\right\}_{0\le s\le \rho-1}$  in the Number Wall $NW_{\bm{\gamma}}$ vanishes, and that the integer $\rho\ge 1$ is the maximal one for this property to hold. The proof thus boils down to establishing that, under the divisibility assumptions,  $(j,k)$ is the top--left corner of a window. This is obvious when $j=1$; assume therefore that  $j\ge 2$. \\

\noindent The claim  is a direct consequence of the explicit expressions in Proposition~\ref{convframe} for the convergents achieving the normal orders in the top and left parts of a window. Indeed, these expressions imply that  if the entry with coordinates $(j-1, k)$  were to also be part of the  frame of the window, then the denominator $Q(t)$ would be divisible by a power of $t$;  if the entry with coordinates $(j, k-1)$  were to also be part of the  frame of the window, then the polynomial $P(t)+\gamma_{m-1}\cdot Q(t)$ would be divisible by a power of $t$. This suffices to conclude the proof.
\end{proof}

%

\subsection{The Process of Generation of Windows in the Number Wall of the Paperfolding Sequence}\label{subprocess}

\noindent With the statement of Corollary~\ref{keyobstr} in mind, in the case where $H_{\bm{\gamma}}(2m,2n)=H_{\bm{\gamma}}(m,n)=0$ over $\K$ for some integers $m\in\Z$ and  $n\ge 1$, the window containing the entry with coordinates $\iota(m,n)$ is said to \emph{induce} the one containing the entry with coordinates $\iota(2m,2n)$.

\begin{rem} \label{rem1}
It is plain from the definition of the indexation map $\iota$ in~\eqref{indexbij} that any window $W'$ with size $\rho(W')\ge 2$ contains an entry with coordinates of the form  $\iota(2m,2n)$ ($m\in\Z, n\ge 1$). From Corollary~\ref{keyobstr}, any such window is thus induced by another one containing the entry with coordinates $\iota(m,n)$. \qed
\end{rem}

\noindent It is convenient to conventionally extend the definitions of a window and of the  inducement relation  to avoid unnecessary distinction of cases in the follwing way~: four nonzero entries located at $(j(W), k(W)), (j(W), k(W)+1), (j(W)+1,k(W))$ and $\left(j(W)+1, k(W)+1\right), $ are seen to be the frame an \emph{empty window} $W$ with size $\rho(W)=0$.  \\

\noindent What is more, given an empty window $W$,   the denominator $Q_{\bm{\gamma}}\left(W; t\right)$ is conventionally set to equal $Q_{\bm{\gamma}}\left(m(W), h(W); t\right)$, where the integer $h(W)$ is here also chosen so that  $q_{\bm{\gamma}}(W):=\deg Q_{\bm{\gamma}}(W;t) = j(W) = n(W)$ (as in relations~\eqref{defnmW} and~\eqref{defqlamnw}). 

\begin{rem}\label{empywinddiv}
The convergent of an empty window still meets the conclusion of Proposition~\ref{propconvwind}. Indeed, equation~\eqref{approxcovpolywin} (specialised to the case $\rho(W)=0$) is just a consequence of  the equivalence~\eqref{convcharac} characterising the convergents of a Laurent series and of Point~(3) in Lemma~\ref{closedromlem}. As for the property that none of the polynomials $P_{\bm{\gamma}}\left(W;t\right)+\gamma_{\left(m(W)-1\right)}\cdot Q_{\bm{\gamma}}\left(W;t\right)$ or $Q_{\bm{\gamma}}\left(W;t\right)$ is divisible by  $t$, the proof of Proposition~\ref{propconvwind} is easily seen to cover the case of empty windows also. \qed
\end{rem}

\noindent In the case where a window $W'$ is not induced by another one (which imposes that  $\rho(W')=1$ from Remark~\ref{rem1}), it is said to be induced by the empty window $W$ as soon as
\begin{equation}\label{frameinduc}
\left((j(W'), k(W'\right) \;=\; \left(2\cdot j(W), \; 2\cdot k(W)\right)\qquad\mbox{and}\qquad t+1\mid Q(W;t).\\
\end{equation}
This purely formal definition is made with a view towards the statement of Proposition~\ref{genwindow} below. It enables one to conveniently describe the process of generating all the  windows in the Number Wall of the paperfolding sequence. 

\begin{rem}\label{rem2}
One of the  purposes behind   the definition of the denominator $Q_{\bm{\gamma}}\left(W; t\right)$ for empty wondows $W$ is made apparent  when combining it with the expression, given in Proposition~\ref{convframe}, for the convergents associated to the frame of a window. Indeed, it is then plain that  any of the polynomials $Q_{\bm{\gamma}}\left(m, h; t\right)$, where $m\in\Z$ and $h\ge 1$ are integers, can be expressed as the product of the denominator $Q_{\bm{\gamma}}\left(W; t\right)$ of the convergent of a (possibly empty) window $W$ by a (possibly trivial) power of the variable $t$.\qed
\end{rem}

\noindent All definitions and notations are now specialised to the case where $\bm{\gamma}$ is the doubly infinite paperfolding sequence $\bm{\lambda}$ over the alphabet $(-1, 1)$. Let then $\left\{\Lambda_m(t)\right\}_{m\in\Z}$ be the family of Laurent series   defined as in~\eqref{doublylaurent}. The recurrence characterisation~\eqref{reclambda} of the sequence   $\bm{\lambda}$ yields  the functional equations
\begin{alignat*}{3}\label{funceqt}
\Lambda_{4s}(t) &\quad= \quad  &&\Lambda_{2s}\left(t^2\right)+\frac{t}{1+t^2};\\
\Lambda_{4s+2}(t) &\quad= \quad &&\Lambda_{2s+1}\left(t^2\right)-\frac{t}{1+t^2};\\
\Lambda_{4s+1}(t) &\quad= \quad  &&t\cdot \Lambda_{2s}\left(t^2\right)-\frac{1}{1+t^2}\qquad\qquad\textrm{and}\\
\Lambda_{4s+3}(t) &\quad= \quad  &&t\cdot \Lambda_{2s+1}\left(t^2\right)+\frac{1}{1+t^2}
\end{alignat*}
valid for all $s\in\Z$. In particular, given $k\in\Z$, the first two identities amount to
\begin{equation}\label{funceqpair}
\Lambda_{2k}\left(t\right)\;=\; \Lambda_k\left(t^2\right)+(-1)^k\cdot \frac{t}{1+t^2}\cdotp\\
\end{equation}

\noindent From now on, \begin{quotation} \begin{center} \emph{fix a prime number $\ell\equiv 3\pmod{4}$}.\end{center}\end{quotation} From Theorem~\ref{red2}, it suffices to prove the Main Theorem~\ref{mainthm} over $\F_\ell$ . Consequently, the Number Wall $NW_{\bm{\lambda}}$ of the paperfolding sequence is  from now on seen as a doubly infinite array over this field. 
The process of generation of windows in $NW_{\bm{\lambda}}$  is determined in the  following statement.

\begin{prop}[Process of generation of windows in the Number Wall]\label{genwindow}
Keep the above notations and assumptions and let $W$ be a (possibly empty) window of size $\rho(W)\ge 0$ in the Number Wall $NW_{\bm{\lambda}}$ over $\F_\ell$. Then,
\begin{itemize}
\item[\textbf{(a)}] under the assumption that $t+1$ divides $Q_{\bm{\lambda}}\left(W; t\right)$, the window $W$ induces another window $W'$ of size $\rho(W')=2\cdot \rho(W)+1$ such that $$ \left(j\left(W'\right), \; k\left(W'\right)\right)=  2\cdot\left( j\left(W\right), \;  k\left(W\right)\right).$$ What is more,
\begin{equation}\label{denomacseb0}
Q_{\bm{\lambda}}\left(W'; t\right)\;=\; Q_{\bm{\lambda}}\left(W; t^2\right).
\end{equation}
\item[\textbf{(b)}] under the assumptions that $t+1$ does not divide $Q_{\bm{\lambda}}\left(W; t\right)$ and that $\rho(W)\ge 3$, the window $W$ induces another window $W'$ of size $\rho(W')=2\cdot \rho(W)-3$ such that $$ \left(j\left(W'\right), \; k\left(W'\right)\right)=  2\cdot \left(j\left(W\right)+1, \;  k\left(W\right)+1\right).$$ What is more,
\begin{equation}\label{denomacseb}
Q_{\bm{\lambda}}\left(W'; t\right)\;=\; \left(t^2+1\right)\cdot Q_{\bm{\lambda}}\left(W; t^2\right).
\end{equation}
\item[\textbf{(c)}] in the remaining case where $t+1$ does not divide $Q_{\bm{\lambda}}\left(W; t\right)$ and where $\rho(W) \in\{0,1\}$, 
the window $W$ does not induce any other window.
\end{itemize}
\end{prop}

\noindent The situation where $\rho(W)=2$ is excluded from Case (c) as it is known from Corollary~\ref{odwind} that the size of a  nonempty window is necessarily odd. 

\begin{proof}[Proof of Proposition~\ref{genwindow}]
In order to establish Parts (a) and (b) of the statement, fix a window $W$ in the Number Wall with size $\rho(W)\ge 0$ and with top--left corner located at $\left(j(W)+1, k(W)+1\right)$, where $j(W)\ge 1$ and where $\left(j(W), k(W)\right)=\iota\left(m(W), n(W)\right)=\left(n(W), m(W)+n(W)\right) $ (note that a corner of a nonempty window cannot be located in the first row since it contains the values of the paperfolding sequence over the alphabet $\pm 1$). When nonempty (i.e.~when $\rho(W)\ge 1$), the window $W$ is determined  by the set of entries with coordinates $(j,k)=\iota\left(m, n\right)$ such that
\begin{align}\label{coord1}
j(W)+1 \;\le\; j\;& \le\; j(W)+\rho(W)
\end{align}
and
\begin{align}\label{coord2}
k(W)+1\;\le\; k \;&\le\; k(W)+\rho(W).
\end{align}
The proof relies on two properties of the paperfolding sequence. The first one is a particular case of the recurrence relations~\eqref{reclambda}, namely that
\begin{equation}\label{2mfctmpapfld}
\lambda_{2m+1}\;=\; \lambda_m \qquad \textrm{for all } m\in\Z.
\end{equation}
The second one is the functional equation~\eqref{funceqpair}, which implies the following set of relations~:
\begin{align}
&\left| Q_{\bm{\lambda}}\left(W;t^2\right) \cdot   \Lambda_{m(W)}(t^2)-P_{\bm{\lambda}}\left(W;t^2\right)\right|\nonumber \\
&\qquad \qquad \qquad  \underset{\eqref{approxcovpolywin}}{=}\; 2^{-2\cdot q_{\bm{\lambda}}(W)-2\cdot \rho(W)-2}\nonumber\\
&\qquad  \qquad \qquad \underset{\eqref{funceqpair}}{=}\;  \left| Q_{\bm{\lambda}}\left(W;t^2\right) \cdot   \left(\Lambda_{2\cdot m(W)}(t)-(-1)^k\cdot \frac{t}{t^2+1}\right) -P_{\bm{\lambda}}\left(W;t^2\right)\right|. \label{approxcovpolywinbister}
\end{align}
They lead one  to a natural distinction of cases, each corresponding to one of the situations described in Parts (a) and (b) of Proposition~\ref{genwindow}.

\paragraph{$\bullet$\textbf{ Case 1.}}  Assume first that   $t+1$ divides $Q_{\bm{\lambda}}\left(W;t\right) $ (as  in Part (a) of the statement) and set $$\widehat{Q}_{\bm{\lambda}}\left(W;t\right)\;:=\;\frac{Q_{\bm{\lambda}}\left(W;t\right)}{t+1}\in\F_\ell\left[t\right].$$
The goal is to show that the fraction
\begin{equation}\label{fraccase1irr}\frac{P_{\bm{\lambda}}\left(W;t^2\right)+(-1)^k\cdot t\cdot \widehat{Q}_{\bm{\lambda}}\left(W;t^2\right)}{Q_{\bm{\lambda}}\left(W;t^2\right)}
\end{equation}
is the convergent $P_{\bm{\lambda}}\left(W',t\right) /Q_{\bm{\lambda}}\left(W',t\right) $ of the window $W'$ with size $\rho(W')=2\cdot \rho(W)+1$ and   top--left corner located at
\begin{align*}
\left(j(W')+1, k(W')+1\right)\;&:=\; \left(2\cdot q_{\bm{\lambda}}(W)+1, \; 2\cdot m(W)+2\cdot q_{\bm{\lambda}}(W)+1\right)  \\
& \underset{ \eqref{defqlamnw}}{=} \;  2\cdot \left(j(W), \; k(W)\right)+ (1,1),
\end{align*}
and then to show that  $W'$ is induced by $W$.\\

\noindent Note that Proposition~\ref{propconvwind} guarantees that $t$ can divide neither the polynomial $Q_{\bm{\lambda}}\left(W;t\right)$ nor the polynomial $P_{\bm{\lambda}}\left(W;t\right)+\lambda_{\left(m(W)-1\right)}\cdot Q_{\bm{\lambda}}\left(W;t\right)$. It can therefore divide neither $Q_{\bm{\lambda}}\left(W;t^2\right)$ nor
\begin{align*}
P_{\bm{\lambda}}\left(W;t^2\right)&+(-1)^k\cdot t\cdot \widehat{Q}_{\bm{\lambda}}\left(W;t^2\right)+\lambda_{\left(2\cdot m(W)-1\right)}\cdot Q_{\bm{\lambda}}\left(W;t^2\right)\\
&\underset{\eqref{2mfctmpapfld}}{=}\; \left(P_{\bm{\lambda}}\left(W;t^2\right)+\lambda_{\left(m(W)-1\right)}\cdot Q_{\bm{\lambda}}\left(W;t^2\right)\right)+ (-1)^k\cdot t\cdot \widehat{Q}_{\bm{\lambda}}\left(W;t^2\right).
\end{align*}
It remains to show that the fraction in~\eqref{fraccase1irr} is irreducible. To see this, note first that the polynomials $P_{\bm{\lambda}}\left(W;t^2\right)$ and $Q_{\bm{\lambda}}\left(W;t^2\right)$ are coprime (as can be seen by specialising  at the value $t^2$ a Bézout identity satisfied by the coprime polynomials $P_{\bm{\lambda}}\left(W;t\right)$ and $Q_{\bm{\lambda}}\left(W;t\right)$). Since $Q_{\bm{\lambda}}\left(W;t\right)$ is not di\-vi\-si\-ble by $t$, the greatest common divisor of the numerator and the denominator of the fraction  is either 1 or the polynomial $t^2+1$ (which is irreducible over $\F_\ell$ when $\ell\equiv 3\pmod{4}$). Assume the latter for a contradiction. Then, $$P_{\bm{\lambda}}\left(W;t^2\right)+t\cdot \widehat{Q}_{\bm{\lambda}}\left(W;t^2\right)\;\equiv\; P_{\bm{\lambda}}\left(W,-1\right)+t\cdot \widehat{Q}_{\bm{\lambda}}\left(W,-1\right)\;\equiv\; 0 \pmod{t^2+1}.$$ This implies that the coefficients of the above degree one polynomial vanish over the field $\F_\ell\left[t\right]/\left(\left(t^2+1\right)\cdot\F_\ell\left[t\right]\right)$, and thus in particular that the constant term $P_{\bm{\lambda}}\left(W,-1\right)$ vanishes over $\F_\ell$. This means that  $t+1$ divides  $P_{\bm{\lambda}}\left(W;t\right)$, contradicting its coprimality with $Q_{\bm{\lambda}}\left(W;t\right) $.\\

\noindent To conclude the proof of Part (a) in Proposition~\ref{genwindow}, it remains to check that the window $W$ induces $W'$. This is clear when $W$ is an empty window~:  equations~\eqref{frameinduc} are indeed evidently verified in this case. Assume therefore that $\rho(W)\ge 1$ and fix a pair of coordinates $(j,k)=\iota(m,n)$ inside
the window $W$ as determined by inequalities~\eqref{coord1} and~\eqref{coord2}~: in particular, $H_{\bm{\lambda}}\left(m, n\right)= 0$.
It then suffices to show that the pair $\left(m, n\right)=\iota^{-1}(j,k)$ induces the pair $\left(m', n'\right)=(2m, 2n)$ in the sense that the coordinates
\begin{align}
(j',k')\; &:=\; \iota\left(m', n'\right) \; =\; \iota\left(2m, 2n\right)
\underset{\eqref{indexbij}}{=}\; 2\cdot\iota(m,n) \; =\;  \left(2j, \; 2k\right)\label{correspo2m2nnm}
\end{align}
are inside the window $W'$ (in particular,  one then has that $H_{\bm{\lambda}}\left(m', n'\right)=0$).
This indeed holds as relations~\eqref{coord1} and~\eqref{coord2} yield that  $$j(W')+1\;=\; 2\cdot j(W)+1\;\le\; 2j\;\le\; 2\cdot j(W)+2\cdot \rho(W)+1\;=\; j(W')+\rho(W')$$ and $$k(W')+1\;=\;2\cdot k(W)+1\; \le\; 2k\;\le\; 2\cdot k(W)+2\cdot \rho(W)+1\;=\; k(W')+\rho(W')+1.$$ The proof of Part (a) is thus complete.

\paragraph{$\bullet$\textbf{ Case 2.}}  Assume now that $t+1$ does \emph{not} divide $Q_{\bm{\lambda}}\left(W;t\right) $. Multiply equation~\eqref{approxcovpolywinbister} by $t^2+1$ to obtain
\begin{align*}
&\left| \left(t^2+1\right)\cdot Q_{\bm{\lambda}}\left(W;t^2\right) \cdot  \Lambda_{2\cdot m(W)}(t) - \left(\left(t^2+1\right) \cdot P_{\bm{\lambda}}\left(W;t^2\right)+(-1)^k\cdot t\cdot Q_{\bm{\lambda}}\left(W;t^2\right)\right)\right|\\
&\qquad  =\; 2^{-2\cdot q_{\bm{\lambda}}(W)-2\cdot \rho(W)}\\
&\qquad  =\; 2^{-\left(2\cdot q_{\bm{\lambda}}(W)+2\right)-(2\cdot \rho(W)-3)-1}.
\end{align*}
This relation gives a convergent only when $\rho(W)>1$. Since Corollary~\ref{odwind} imposes that all nonempty windows have uneven sizes, this condition actually translates into $\rho(W)\ge 3$. Assume this inequality holds (as in the hypotheses of Part (b))
. Reproducing (a simplified version of) the arguments developed in  Case 1, one deduces that the fraction
\begin{equation}\label{fracoirrcaseb}
\frac{\left(t^2+1\right) \cdot P_{\bm{\lambda}}\left(W;t^2\right)+(-1)^k\cdot t\cdot Q_{\bm{\lambda}}\left(W;t^2\right)}{ \left(t^2+1\right)\cdot Q_{\bm{\lambda}}\left(W;t^2\right)}
\end{equation}
is irreducible and is the convergent   $P_{\bm{\lambda}}\left(W',t\right) /Q_{\bm{\lambda}}\left(W',t\right) $  of the window $W'$ with size $\rho(W'):=2\cdot \rho(W)-3$ and top--left corner
\begin{align*}
\left(j\left(W'\right)+1, k\left(W'\right)+1\right) \;&:=\; \left(2\cdot q_{\bm{\lambda}}(W)+3, \; 2\cdot m(W)+2\cdot q_{\bm{\lambda}}(W)+3\right)\\
&\underset{\eqref{defqlamnw}}{=}\;   \left(2\cdot j(W)+3, \; 2\cdot k(W)+3\right).
\end{align*}
In order to see that $W'$ is induced by $W$, fix a pair of coordinates $(j,k)=\iota(m,n)$ lying in the ranges~\eqref{coord1} and~\eqref{coord2} (so that in particular $H_{\bm{\lambda}}\left(m, n\right)=0$).
In fact,  under the assumption that $\rho(W)\ge 3$, this pair may be chosen in such a way that all inequalities in~\eqref{coord1} and  in~\eqref{coord2} are strict. These inequalities then yield that $$j(W')+1\;=\;2\cdot j(W)+3\;\le\; 2j\;\le\; 2\cdot j(W)+2\cdot \rho(W)\;=\; j(W')+\rho(W')+1$$ and $$k(W')+1\;=\; 2\cdot k(W)+3\;\le\; 2k\;\le\; 2\cdot k(W)+2\cdot \rho(W)\;=\; k(W')+\rho(W')+1,$$ which shows that the resulting coordinates
\begin{align*}
(j',k') \; :=\; \left(2j,\, 2k\right)\;\underset{\eqref{correspo2m2nnm}}{=}\;  \iota\left(2m,\,2n\right)\; :=\; \iota\left(m',n'\right)
\end{align*}
are those of an entry in $W'$. Therefore, $H_{\bm{\lambda}}\left(m', n'\right)=0$ and
the window $W'$ is induced by $W$. This  completes the proof of Part (b) in  Proposition~\ref{genwindow}.

\paragraph{$\bullet$}  Finally, to deal with Part (c) in Proposition~\ref{genwindow},  observe first that the claim therein is trivially satisfied by empty windows $W$  given the very statement of definition~\eqref{frameinduc}. Assume therefore that $W$ is nonempty and consider a  window $W'$ with size $\rho(W')\ge 3$   induced by $W$. The goal is to show that $W$ meets the assumptions of either Part (a) or Part (b) of the statement. \\

\noindent Let $(j',k')$ be the coordinates of an entry in $W'$  
meeting this property~:  there exist integers $m\in\Z$ and $n\ge 1$ such that
\begin{equation}\label{indexcocnrgu}
(j',k')\;=\;\iota(2m,2n) \qquad\; \textrm{and}\qquad \;H_{\bm{\lambda}}\left(2m, 2n\right)\;=\;  H_{\bm{\lambda}}\left(m, n\right)\;=\; 0. 
\end{equation}
Choose them in such a way that the 
corresponding entry $NW_{\bm{\lambda}}(j',k')=H_{\bm{\lambda}}(2m,2n)$ in the Number Wall is the center of the $3\times 3$ zero square located in the top--left corner of $W'$ (this is permitted by  Corollary~\ref{odwind}). \\

\noindent  Note that one has necessarily $n\ge 2$ for the last equation in~\eqref{indexcocnrgu} to hold  (when $n=1$, a quantity of the form $H_{\bm{\lambda}}\left(m, 1\right)$ reduces to a value taken by the paperfolding sequence over the alphabet $\pm 1$, which never vanishes.
It then follows from  Corollary~\ref{odwind} that $H_{\bm{\lambda}}\left(2m, 2n-1\right)= H_{\bm{\lambda}}\left(2m, 2n+1\right)= 0$
(see Figure~\ref{fig1} to confirm that  this is compatible with the indexation in the Number Wall). From identity (Z) in Corollary~\ref{labelhenkelbis}, one thus infers that
\begin{equation}\label{hankelF}
F_{\bm{\lambda}}\left(m, n\right)\;=\; F_{\bm{\lambda}}\left(m, n-1\right)\;=\; 0.
\end{equation}
Recall now that $\left(F_{\bm{\lambda}}\left(m, n\right)\right)_{m\in\Z, n\ge 1}$ is by definition the Number Wall of the sequence $\bm{\lambda'}=\left(\lambda_k+\lambda_{k+1}\right)_{k\in\Z}$ identified with the fractional part of the Laurent series \mbox{$\left(t+1\right)\cdot \Lambda(t)$}, namely the Laurent series  $\left(t+1\right)\cdot \Lambda(t)-\lambda_1$ (see Section~\ref{prelim}). Denote by $B^*_{\bm{\lambda}}\left(W^*, t\right)/A^*_{\bm{\lambda}}\left(W^*, t\right)$ the convergent of a generic window $W^*$ of size $\rho(W^*)\ge 1$ in this Number Wall, and let $$a^*_{\bm{\lambda}}\left(W^*\right)\;=\;\deg A^*_{\bm{\lambda}}\left(W^*, t\right).$$

\noindent Specify now $W^*$ to be  the window determined by the Hankel determinants appearing in~\eqref{hankelF}. In the notation introduced in~\eqref{doublylaurent}, it is readily checked that
\begin{equation}\label{LFlam}
\left(\left(t+1\right)\cdot \Lambda(t)-\lambda_1\right)_{m(W^*)}\; =\; \left(t+1\right)\cdot \Lambda_{m(W^*)}(t)-\lambda_{m(W^*)}.
\end{equation}
Then, upon defining the coordinates $$(j,k)\;:=\; \iota(m,n)\;\underset{\eqref{correspo2m2nnm}}{=}\; \left(\frac{j'}{2}, \, \frac{k'}{2}\right),$$  inequalities~\eqref{coord1} and~\eqref{coord2}
can be restated in this case with the help of the equations in~\eqref{hankelF} as
\begin{equation}\label{nw*j}
j(W^*)+1\;\le\; n-1\;<\; n\;=\; j \;\le\; j(W^*)+\rho(W^*)
\end{equation}
and
\begin{equation}\label{nw*jbis}
k(W^*)+1 \;\le\; m+n-1 \;=\; k-1 \;<\; k \;\le\; k(W^*)+\rho(W^*).
\end{equation}

\noindent Also, from Proposition~\ref{propconvwind}, the convergent of the window $W^*$ meets the relations
\begin{align}
&\left|A^*_{\bm{\lambda}}\left(W^*, t\right)\cdot \left(\left(t+1\right)\cdot \Lambda_{m(W^*)}(t)-\lambda_{m(W^*)} \right)-B^*_{\bm{\lambda}}\left(W^*, t\right)\right|\nonumber \\
&\qquad\qquad\qquad \qquad\qquad\qquad= \left| \left(A^*_{\bm{\lambda}}\left(W^*, t\right)\cdot\left(t+1\right)\right)\cdot \Lambda_{m(W^*)}(t)-D^*_{\bm{\lambda}}\left(W^*, t\right)\right|\nonumber \\
&\qquad \qquad\qquad\qquad\qquad\qquad=\; 2^{-a^*_{\bm{\lambda}}\left(W^*\right)-\rho(W^*)-1},\label{lastlb}
\end{align}
where
\begin{equation*}\label{defrdern}
D^*_{\bm{\lambda}}\left(W^*, t\right)\;:=\;B^*_{\bm{\lambda}}\left(W^*, t\right)+ \lambda_{m(W^*)}\cdot A^*_{\bm{\lambda}}\left(W^*, t\right)
\end{equation*}
and where $t$ divides neither the polynomial $A^*_{\bm{\lambda}}\left(W^*, t\right)$ nor the polynomial $B^*_{\bm{\lambda}}\left(W^*, t\right)+\left(\lambda_{\left(m(W^*)-1\right)}+\lambda_{m(W^*)}\right)\cdot A^*_{\bm{\lambda}}\left(W^*, t\right)$. Furthermore, the congruences in~\eqref{hankelF} impose that
\begin{equation}\label{sizewin*}
\rho(W^*)\ge 2.
\end{equation}
Two possibilities then arise~:

\paragraph{\textbf{(i)}} if  $(t+1)$ does not divide $D^*_{\bm{\lambda}}\left(W^*, t\right)$, then the polynomials $A^*_{\bm{\lambda}}\left(W^*, t\right)\cdot\left(t+1\right)$ and $D^*_{\bm{\lambda}}\left(W^*, t\right)$ are coprime. Furthermore, $t$ divides neither $\left(t+1\right)\cdot A^*_{\bm{\lambda}}\left(W^*, t\right)$ nor
\begin{align*}
&D^*_{\bm{\lambda}}\left(W^*, t\right)+\lambda_{\left(m(W^*)-1\right)}\cdot A^*_{\bm{\lambda}}\left(W^*, t\right)\cdot\left(t+1\right)\\
&\qquad \quad =\; B^*_{\bm{\lambda}}\left(W^*, t\right)+\left( \lambda_{\left(m(W^*)-1\right)}+\lambda_{m(W^*)}\right)\cdot A^*_{\bm{\lambda}}\left(W^*, t\right)+t\cdot \lambda_{\left(m(W^*)-1\right)}\cdot A^*_{\bm{\lambda}}\left(W^*, t\right).
\end{align*}
Equation~\eqref{lastlb} and Corollary~\ref{CSconvwin} then yield that the irreducible fraction $$\frac{P_{\bm{\lambda}}\left(W; t\right)}{Q_{\bm{\lambda}}\left(W; t\right)}\;:=\; \frac{D^*_{\bm{\lambda}}\left(W^*, t\right)}{A^*_{\bm{\lambda}}\left(W^*, t\right)\cdot \left(t+1\right)}$$ is the convergent of the window $W$ in the Number Wall of the sequence $\bm{\lambda}$ such that
\begin{equation}\label{mWW*}
q_{\bm{\lambda}}\left(W\right)\;=\;a^*_{\bm{\lambda}}\left(W^*\right)+1, \quad m(W)=m(W^*) \quad \textrm{ and }\quad \rho(W)=\rho(W^*)-1\underset{\eqref{sizewin*}}{\ge} 1.
\end{equation}
Consequently,
\begin{align*}
j(W)+1\;&\underset{\eqref{defqlamnw}}{=}\;q_{\bm{\lambda}}\left(W\right)+1\;=\;a^*_{\bm{\lambda}}\left(W^*\right)+2\;\underset{\eqref{defqlamnw}}{=}\;j(W^*)+2\\
&\underset{\eqref{nw*j}}{\le}\; j \;\;\underset{\eqref{nw*j}}{\le}\;\;  j(W^*)+\rho(W^*)\;\underset{\eqref{defqlamnw}\&\eqref{mWW*}}{=}\;  j(W)+\rho(W)
\end{align*}
and
\begin{align*}
k(W)+1\;&=\; m(W)+n(W)+1\;\underset{\eqref{defqlamnw}\&\eqref{mWW*}}{=}\; m(W^*)+n(W^*)+2\;\underset{\eqref{nw*jbis}}{=}\; k(W^*)+2\\
&\underset{\eqref{nw*jbis}}{\le}\;k\;\underset{\eqref{nw*jbis}}{\le}\; k(W^*)+\rho(W^*)\;\underset{\eqref{defqlamnw}\&\eqref{mWW*}}{=}\; k(W)+\rho(W).
\end{align*}
The integers $j$ and $k$ thus satisfy  inequalities~\eqref{coord1} and~\eqref{coord2}; as a consequence,
the point with coordinates $(j,k)=\iota(m,n)$ belongs to the window $W$. This shows, on the one hand that the considered window $W$ is the one inducing $W'$ and, on the other, that this case falls within the assumption of Part (a) of the proposition.

\paragraph{\textbf{(ii)}} if $(t+1)$ divides $D^*_{\bm{\lambda}}\left(W^*, t\right)$, then, upon setting $D^*_{\bm{\lambda}}\left(W^*, t\right)=\left(t+1\right)\cdot \widetilde{D}^*_{\bm{\lambda}}\left(W^*, t\right)$, equation~\eqref{lastlb} reads $$\left|A^*_{\bm{\lambda}}\left(W^*, t\right)\cdot \Lambda_{m(W^*)}(t)-\widetilde{D}^*_{\bm{\lambda}}\left(W^*, t\right)\right|\;=\; 2^{-a^*_{\bm{\lambda}}\left(W^*\right)-(\rho(W^*)+1)-1}.$$ Reproducing (a simpler version of) the above arguments, Corollary~\ref{CSconvwin} implies that the irreducible fraction  $$\frac{P_{\bm{\lambda}}\left(W; t\right)}{Q_{\bm{\lambda}}\left(W; t\right)}\;:=\; \frac{ \widetilde{D}^*_{\bm{\lambda}}\left(W^*, t\right)}{A^*_{\bm{\lambda}}\left(W^*, t\right)}$$ is the convergent of the window $W$ in the Number Wall of the sequence $\bm{\lambda}$ such that
\begin{equation*}\label{mWW*bis}
q_{\bm{\lambda}}\left(W\right)\;=\;a^*_{\bm{\lambda}}\left(W^*\right), \quad m(W)=m(W^*) \quad \textrm{ and }\quad \rho(W)\;=\;\rho(W^*)+1.
\end{equation*}
In fact, the known lower bound~\eqref{sizewin*} on the quantity $\rho(W^*)$ shows that, in this case, $\rho(W)\ge 3$ (this can also be recovered from Corollary~\ref{odwind} stating that the windows have odd sizes).\\

\noindent Straightforward adaptations of the above calculations yield in the same way that $j(W)+1\le j\le j(W)+\rho(W)$ and that $k(W)+1\le k\le k(W)+\rho(W)$, implying that the point with coordinates $(j,k)=\iota(m,n)$ belongs here also to the window $W$. As above, one infers on the one hand that the window $W$ induces $W'$ and, on the other, that this case falls within the assumptions of Part (b) of the statement.

\paragraph{$\bullet$} It follows from this distinction of cases that, if a window $W'$ with size $\rho(W')\ge 3$ is induced by another nonempty one, say $W$, then $W$ meets the assumptions of either Part~(a) or else of Part (b) of Proposition~\ref{genwindow}. To completes the proof of Part (c), it remains to notice that if a window $W'$ is induced by a nonempty one, then it has necessarily  size $\rho(W')\ge 3$.   This is immediate~: by the definition of the inducement process, $W'$ contains an entry with coordinates $\iota(2m,2n)$ which, from Corollary~\ref{odwind}, is then the center of a $3\times 3$ zero square.
This is enough to conclude the proof of Part~(c) and thus of Proposition~\ref{genwindow}.
\end{proof}

\subsection{An Arithmetic Condition for the Paperfolding Laurent Series  to Fail $P(t)$--LC}

\noindent Proposition~\ref{genwindow} admits the following corollary. It gives  a sufficient condition, stated in terms of  divisibility properties of the denominators of convergents,  for the paperfol\-ding Laurent series   to fail $P(t)$--LC over $\F_\ell$.

\begin{coro}[Arithmetic Condition for the Failure of $t$--LC]\label{arithCNS}
Recall that $\ell$ is a fixed prime congruent to 3 modulo 4, and  that, given $m\in\Z$, $\Lambda_m(t)$ stands for  the Laurent series derived  as in~\eqref{doublylaurent} from the paperfolding Laurent series $\Lambda(t)$ (defined in~\eqref{laurlambda}).
The coefficients of $\Lambda(t)$ are seen as elements of the finite field $\F_\ell$, and so are the entries of its number wall $NW_{\bm{\lambda}}$. 
Then, the following two statements are equivalent~:
\begin{itemize}
\item[\textbf{(i)}] the Number Wall $NW_{\bm{\lambda}}$ admits no  window of size strictly bigger than 3;
\item[\textbf{(ii)}] For all integers $h\ge 1$ and $m\in \Z$, the denominator $Q_{\bm{\lambda}}(m,h;t)$ of the convergent associated to the Laurent series $\Lambda_m(t)$ is not divisible by $t^2-1$.
\end{itemize}
If any of these conditions is met, then
\begin{equation}\label{qtvefail}
\inf_{Q(t)\in\F_\ell[T]\backslash\{\bm{0}\}}\; \left| Q(t)\right| \cdot\left|Q(t)\right|_{t} \cdot \left| \left\langle Q(t)\cdot\Lambda(t)\right\rangle\right|\;\ge \; 2^{-4}.
\end{equation}
In particular,  $\Lambda(t)$ then fails $t$--LC over $\F_\ell$.
\end{coro}


\begin{proof}
The last claim follows from the characterisation, stated in Theorem~\ref{red2},  of the failure of $t$--LC in terms of Hankel determinants~: inequality~\eqref{qtvefail} is equivalent to the non-existence of windows of size strictly bigger than 3 in the number wall $NW_{\bm{\lambda}}$ restricted to the entries with indices $m\ge 1$. From the definition of the  indexation map $\iota$ in~\eqref{indexbij}, this is the region with entries with  second coordinates $k\ge 2$. Clearly, this is implied by (i) which places no restriction on the location of the windows.

\paragraph{$\bullet$} As pointed out in Remark~\ref{rem2}, given any integers $h\ge 1$ and $m\in \Z$, there exist a (possibly empty) window $W$ and an exponent $s\ge 0$ such that   $Q_{\bm{\lambda}}(m,h;t)=t^s\cdot Q_{\bm{\lambda}}(W;t)$. Since the factor $t^s$ does not affect the property of the polynomial to be divisible by $t^2-1$, it is enough to establish the equivalence between statements (i) and (ii) in the case where one considers denominators of convergents of windows (i.e.~polynomials of the form $Q_{\bm{\lambda}}(W;t)$ for some window $W$).

\paragraph{$\bullet$} To show that (i) implies (ii), argue by contraposition assuming the existence of a window $W$   such that  $t^2-1$, and therefore $t+1$  divides $Q_{\bm{\lambda}}(W;t)$. 
It then follows from Part (a) of Proposition~\ref{convframe} that $W$ generates a window of size $2\cdot \rho(W)+1$ whose denominator (namely, $Q_{\bm{\lambda}}(W;t^2)$) is also divisible by $t^2-1$, and therefore by $t+1$. By iterating the application of Part (a) of  Proposition~\ref{convframe}, one obtains an infinite chain of induced windows $$
W_0 \;=\; W\; \longrightarrow\; W_1\; \longrightarrow\; W_2 \;\longrightarrow\; W_3 \;\longrightarrow\; \cdots
$$
with respective sizes $\left(\rho(W_i)\right)_{i\ge 0}$ such that $\rho(W_i)=2\cdot \rho(W_{i-1})+1$ for all $i\ge 1$. In particular, $\rho(W_3)=8\cdot\rho(W_0)+7\ge 7$, which shows that (i) does not hold and thus completes the proof of this implication.

\paragraph{$\bullet$} In order to show that (ii) implies (i), work again by contraposition assuming that there exists a window $W$ with size $\rho(W)\ge 4$. \\

\noindent From Remark~\ref{rem1}, every window of size at least 2 is induced by another nonempty one which, from Proposition~\ref{genwindow}, is of smaller or equal size.  Since the first row of the number wall $NW_{\bm{\lambda}}$ does not intersect any window (recall that it contains the terms of the sequence $\bm{\lambda}$ valued in the set $\{-1,1\}$), one easily  deduces  
that any  maximal chain of nonempty induced windows  starts with a window of size 1. Also, it is plain that any nonempty window, in particular the above--considered window $W$ of size $\rho(W)\ge 4$, belongs to such a chain. \\




\noindent Let then $W_0$ be a window of size $\rho_0:=\rho(W_0)=1$. By  Proposition~\ref{genwindow}, it induces a window $W_1$ if and only if $t+1$ divides $Q_{\bm{\lambda}}(W_0; t)$, in which case $W_1$ has size $\rho_1:=\rho(W_1)= 3$ and denominator
\begin{equation}\label{basecasechainwin}
Q_{\bm{\lambda}}(W_1; t)\;=\;Q_{\bm{\lambda}}\left(W_0; t^2\right).
\end{equation}
In turn, $W_1$ induces a window $W_2$ of size $\rho_2:=\rho(W_2)=7$ if $t+1$ divides  $Q_{\bm{\lambda}}(W_1; t)$ and of size $\rho_2=3$ otherwise. \\

\noindent Consider the maximal chain of nonempty, induced windows $(W_i)_{i\ge 0}$ containing $W$ (in particular, $\rho(W_0)=1$). By the above arguments, $t+1$ divides $ Q_{\bm{\lambda}}(W_0;t)$. Denote by $s = s(W_0)\ge 1$ the largest integer such that $\rho(W_1) = \rho(W_2) = \cdots=\rho(W_s)=3$. Since the chain contains the window $W$ which has size at least 4, such a value $s$ exists and is finite. From Proposition~\ref{genwindow}, the polynomial $t+1$ does not divide $Q_{\bm{\lambda}}(W_r;t)$ for all $1\le r<s$ but it  divides $ Q_{\bm{\lambda}}(W_s;t)$.\\


\noindent An easy induction inferred from Proposition~\ref{genwindow} with base case~\eqref{basecasechainwin} shows that
\begin{equation}\label{cfefe}
Q_{\bm{\lambda}}\left(W_r; t\right)\;=\;\left(\prod_{i=1}^{r-1}\left(t^{2^i}+1\right)\right)\cdot Q_{\bm{\lambda}}\left(W_0; t^{2^r}\right)
\end{equation}
for any integer $r\ge 0$ not greater than $s(W_0)$ (conventionally, an empty product is taken to equal 1). \\

\noindent Since for any value of $r$, the product in the first factor on the right--hand side of~\eqref{cfefe} is clearly not divisible by $t+1$ (it  does not vanish when evaluated at $t=-1$),  $t+1$ divides $Q_{\bm{\lambda}}\left(W_{s}; t\right)$ if and only if it divides the factor  $Q_{\bm{\lambda}}\left(W_0; t^{2^{s}}\right)$. In other words, $Q_{\bm{\lambda}}\left(W_0; (-1)^{2^{s}}\right)=0$ which, under the assumption that $s\ge 1$, means that $t-1$ divides $Q_{\bm{\lambda}}\left(W_0; t\right)$. Since, by construction, $t+1$ divides $Q_{\bm{\lambda}}\left(W_0; t\right)$ (recall that $s=s(W_0)>0$), this implies that $t^2-1$ divides $Q_{\bm{\lambda}}\left(W_0; t\right)$. This shows that (ii) does not hold and thus concludes the proof of this implication.

\end{proof}

\section{Completion of the Proof of the Main Theorem}\label{secverif}

\begin{proof}[Completion of the Proof of Theorem~\ref{mainthm}] Recall that from the reductions operated in Theorems~\ref{red1}  and~\ref{red2}, it suffices to prove the statement in the case of the irreducible polynomial $P(t)=t$ when working over  the finite field $\F_\ell$ for the fixed value of the prime $\ell\equiv 3\pmod{4}$.

\paragraph{$\bullet$}  Direct calculations show that the three determinants $H_{\bm{\lambda}}(2,3)$, $H_{\bm{\lambda}}(2,4)$ and $H_{\bm{\lambda}}(2,5)$
simultaneously vanish. From Theorem~\ref{red2}, one deduces that the paperfolding Laurent series $\Lambda(t)\in\F_\ell\left(\left(t^{-1}\right)\right)$ meets the inequality
\begin{equation}\label{tLCred}
\inf_{Q(t)\in\F_\ell[T]\backslash\{\bm{0}\}}\; \left| Q(t)\right| \cdot\left|Q(t)\right|_{t} \cdot \left| \left\langle Q(t)\cdot\Lambda(t)\right\rangle\right|\;\le\; 2^{-4}.
\end{equation}

\paragraph{$\bullet$} The goal is now to prove that the reverse inequality holds by establishing Pro\-per\-ty~(ii)  in Corollary~\ref{arithCNS}.
To this end, argue by contradiction assuming the existence of integers   $m\in\Z$ and $h\ge 1$ such that
\begin{equation}\label{eq1}
t^2-1\mid Q_{\bm{\lambda}}(m,h;t).
\end{equation}
For the sake of simplicity of notations, set $$n \;:=\; q_{\bm{\lambda}}(m,h)$$ for the corresponding normal order (as defined in~\eqref{defdegnorord}). Clearly, the above divisibility condition implies that $n\ge 2$. Assume furthermore that the pair $(m,n)$ is minimal with respect to the second entry $n$ in the following sense~: $n\ge 2$ is the smallest value for which there exist integers $m$ and $h$ satisfying~\eqref{eq1}. \\

\noindent For such a value of $n$, choose any $m$ for which this divisibility condition is met for some $h\ge 1$. If the entry $\iota(m,n)$ appears to be on the top edge of some nonempty window $W$, i.e. if $H_{\bm{\lambda}}(m-1,n+1)=0$, then shift $m$ to the left so that the entry $\iota(m,n)$ becomes the top-left corner of $W$. Such a choice of the pair $(m,n)$ thus ensures that $H_{\bm{\lambda}}(m-1,n+1)\neq 0$. 


\begin{lem}\label{lemlast}
A minimal pair $(m,n)$ achieving the divisibility relation~\eqref{eq1} satisfies the follo\-wing properties~:
\begin{itemize}
\item[\textrm{\textbf{(a)}}] the three determinants $H_{\bm{\lambda}}(m,n)$,  $H_{\bm{\lambda}}(m+1,n)$ and $H_{\bm{\lambda}}(m-1,n+1)$
do not vanish;
\item[\textbf{(b)}] the polynomial $Q_{\bm{\lambda}}(m,h;t)$ 
is not divisible by $t$. 
\item[\textbf{(c)}] If $Q_{\bm{\lambda}}(m,h;t) = Q_{\bm{\lambda}}(W;t)$ for some nonempty window $W$, then $W$ cannot be induced  by any other window (be it empty or nonempty).
\end{itemize}
\end{lem}

\begin{proof}[Proof of Lemma~\ref{lemlast}]
Each of the properties is established successively.

\paragraph{\textrm{\textbf{(a)}}} Since $n$ is the $(m,h)^{\textrm{th}}$ normal order, $H_{\bm{\lambda}}(m,n)\neq 0$ from Point~(3) in Lemma~\ref{closedromlem}. If $H_{\bm{\lambda}}(m+1,n)=0$,
then from Part~(U) in Proposition~\ref{convframe}, the denominator $Q_{\bm{\lambda}}(m,h;t)$ is a multiple of $t$ (namely, it is  $t$ times the denominator corresponding to the entry with coordinates $\iota(m+1,n-1)$  in the Number Wall), which contradicts the mini\-ma\-lity of the normal order $n$. Finally, $H_{\bm{\lambda}}(m-1,n+1)\neq 0$ by the construction of the entry $(m,n)$.

\paragraph{\textrm{\textbf{(b)}}} From Remark~\ref{rem2},  $Q_{\bm{\lambda}}(m,h;t)$  can be decomposed as  $Q_{\bm{\lambda}}(m,h;t)=t^s\cdot Q_{\bm{\lambda}}(W;t)$ for some integer $s\ge 0$ and for some (possibly empty) window $W$.  Proposition~\ref{propconvwind} (for the case where $W$ is nonempty) and Remark~\ref{empywinddiv} (when $W$ is empty) imply that $Q_{\bm{\lambda}}(W;t)$ cannot be divisible by $t$, and the minimality of the normal order $n$ furthermore imposes in all cases that $s=0$.

\paragraph{\textbf{(c)}} Suppose that $W$ is induced by some (possibly empty) window $W'$. Proceed then as in Corollary~\ref{arithCNS} and place $W$ inside a maximal chain of   induced windows
$$
W_0\to W_1\to\cdots \to W_i = W \to \cdots,
$$
where $W_0$ is a window of size at most one. From Proposition~\ref{genwindow}, the polynomial $t+1$ divides $ Q_{\bm{\lambda}}(W_0;t)$, and one has that for $i\ge 1$,
$$
\mbox{either}\qquad  Q_{\bm{\lambda}}(W_i;t) = Q_{\bm{\lambda}}(W_{i-1};t^2)\qquad\mbox{or else}\qquad Q_{\bm{\lambda}}(W_i;t) = (t^2+1) Q_{\bm{\lambda}}(W_{i-1};t^2),
$$
where $n\left(W_{i-1}\right)<n\left(W_{i}\right)$. By the choice of $W=W_i$,   the polynomial $t-1$ divides $ Q_{\bm{\lambda}}(W;t)$, which in either case implies that $Q_{\bm{\lambda}}(W_{i-1};1)=0$; that is, $t-1\mid Q_{\bm{\lambda}}(W_{i-1};t)$. By iteratively applying the same arguments to $W_{i-1}, W_{i-2}, \dots, W_1$, one finally derives that $t-1$ divides $Q_{\bm{\lambda}}(W_0;t)$, and therefore that $t^2-1\mid Q_{\bm{\lambda}}(W_0;t)$. Since $n(W_0) < n$, this contradicts the minimality of the pair $(m,n)$.
\end{proof}

\paragraph{$\bullet$}  The next step in the proof of the Main Theorem~\ref{mainthm} is to call on the explicit determinant expression for the polynomial $Q_{\bm{\lambda}}(m,h;t)$ provided by Point~(2) in Lemma~\ref{closedromlem}, namely
\begin{equation}\label{detcvgtmn}
Q_{\bm{\lambda}}(m,h;t)\;=\; S_{\bm{\lambda}} (m,n;t)
\end{equation}
(recall here the normalisation for the denominators of convergents stated in Remark~\ref{rem0}). The recurrence formulae (A) --- (D) in Proposition~\ref{proprec} are then shown to lead one to a contradiction when assuming the divisibility condition~\eqref{eq1}. Accordingly, the proof breaks, as in this proposition, into four   cases depending on the parity of the pair $(m,n)$. \\

\noindent The argument in each of the four cases follows the same pattern and requires the consideration of several subcases. The proof is meticulously detailed in the first Case~A; only the changes in it are stated for the others.\\

\noindent Before delving into Case~A, recall that $\left(F_{\bm{\lambda}}(m,n)\right)_{m\in\Z, n\ge 1}$ stands for the Number Wall associated to the fractional part of  the Laurent series $\left(1+t\right)\cdot \Lambda(t)$, namely  $\left(t+1\right)\cdot \Lambda(t)-\lambda_1$. The family of Laurent series parametrised by the integer $m\in\Z$ associated to it as in~\eqref{doublylaurent}  is defined as
\begin{equation}\label{parmlambdat+1}
\left(\left(t+1\right)\cdot \Lambda(t)-\lambda_1\right)_{m}\;=\;\left(t+1\right)\cdot \Lambda_{m}(t)-\lambda_{m}.\\
\end{equation}


\subsection{Case A :  $(m,n)\equiv (1, 0) \pmod{2}$. }

Assume that $(m,n)=(2m'+1, 2n')$ for some  $m'\in\Z$ and $n'\ge 1$. Then,
\begin{align*}
S_{\bm{\lambda}}(2m'+1,2n';t) \!\!\!&\;\;\;\; \equiv\; H_{\bm{\lambda}}(m'+1,n')\cdot S_{\bm{\lambda}}(m',n';1) \\
 &\qquad\qquad\qquad - 2\cdot G_{\bm{\lambda}}(m'+1,n'-1)\cdot V_{\bm{\lambda}}(m',n'-1;1) \\
&\qquad\qquad \qquad+ (-1)^{m+1}\cdot t\cdot F_{\bm{\lambda}}(m',n')\cdot R_{\bm{\lambda}}(m'+1,n'-1;1) \\
&\underset{\eqref{eq1} \& \eqref{detcvgtmn}}{\equiv}\; 0\pmod{t^2-1}.
\end{align*}
Since the polynomial on the right--hand side of the first relation is of degree at most 1, it has to vanish for it to be divisible by $t^2-1$. In particular, this implies that
\begin{equation}\label{fr=01}
 F_{\bm{\lambda}}(m',n')\cdot R_{\bm{\lambda}}(m'+1,n'-1;1)\;=\; 0.
 \end{equation}

\paragraph{$\bullet$} If $F_{\bm{\lambda}}(m',n')=0$ then by indentity~(Z) in  Corollary~\ref{labelhenkelbis}, one has $H_{\bm{\lambda}}(2m',2n'+1) = H_{\bm{\lambda}}(m-1,n+1)=0$ which contradicts Point~(a) in Lemma~\ref{lemlast}. Therefore,  the divisibility condition
\begin{equation}\label{divaprlem}
t-1\mid R_{\bm{\lambda}}(m'+1,n'-1;t)
\end{equation}
is met. This leads one to a   distinction of two cases~:

\paragraph{$\bullet$}\textbf{Case A.1}~:  The polynomial $ R_{\bm{\lambda}}(m'+1,n'-1;t)$ vanishes identically. It  then follows from   Point~(4) in Lemma~\ref{closedromlem} that $F_{\bm{\lambda}}(m',n') = 0.$ However, this condition has already been ruled out.

\paragraph{$\bullet$}\textbf{Case A.2}~: The polynomial $ R_{\bm{\lambda}}(m'+1,n'-1;t)$  does not  vanish and is divisible by $t-1$. Let then $$r_{\bm{\lambda}}\left(m'+1, n'-1\right)\;:=\; \deg R_{\bm{\lambda}}(m'+1,n'-1;t)\;\le\; n'-1$$ and $$R_{\bm{\lambda}}(m'+1,n'-1;t)\;=\; \left(t-1\right)\cdot \widetilde{R}_{\bm{\lambda}}(m'+1,n'-1;t),$$ where $\widetilde{R}_{\bm{\lambda}}(m'+1,n'-1;t)$ is a polynomial of degree $r_{\bm{\lambda}}\left(m'+1, n'-1\right)-1$.
It follows from Point~(2) in Lemma~\ref{closedromlem} that for some polynomials $T_{\bm{\lambda}}(m'+1, n'-1, t) $ coprime with $R_{\bm{\lambda}}(m'+1,n'-1;t)$,
\begin{align}
&\left|  \left(t-1\right)\cdot \widetilde{R}_{\bm{\lambda}}(m'+1,n'-1;t)\cdot \left(\left(t+1\right)\cdot \Lambda(t)-\lambda_1\right)_{\left(m'+1\right)}-T_{\bm{\lambda}}(m'+1, n'-1, t) \right|\nonumber\\
&\qquad\qquad\underset{\eqref{parmlambdat+1}}{=}\;  \left|  \left(t^2-1\right)\cdot \widetilde{R}_{\bm{\lambda}}(m'+1,n'-1;t)\cdot   \Lambda_{\left(m'+1\right)}(t)-T'_{\bm{\lambda}}(m'+1, n'-1, t) \right|\nonumber\\
&\qquad\qquad \le\; 2^{-n'},\label{x}
\end{align}
where
\begin{align}\label{deft'lmnt}
T'_{\bm{\lambda}}(m'+1, n'-1, t)\;=\; T_{\bm{\lambda}}(m'+1, n'-1, t) +\lambda_{\left(m'+1\right)}\cdot \left(t-1\right)\cdot \widetilde{R}_{\bm{\lambda}}(m'+1,n'-1;t).
\end{align}

\noindent Apply the transformation $t\mapsto t^2$ to relation~\eqref{x}. From the functional equation~\eqref{funceqpair} satisfied by the paperfolding Laurent series $\Lambda(t)$, it yields the existence of  polynomials \mbox{$T''_{\bm{\lambda}}(m'+1, n'-1, t)$} and $T'''_{\bm{\lambda}}(m'+1, n'-1, t)$ such that
\begin{align}
& \left|  \left(t^2-1\right)\cdot   \left(t^2+1\right)\cdot \widetilde{R}_{\bm{\lambda}}\left(m'+1,n'-1;t^2\right)\cdot   \Lambda_{\left(2m'+2\right)}(t)-T''_{\bm{\lambda}}(m'+1, n'-1, t) \right| \underset{\eqref{shiftgammam}}{=} \nonumber \\
& \left|  t\cdot   \left(t^2+1\right)\cdot \widetilde{R}_{\bm{\lambda}}\left(m'+1,n'-1;t^2\right)\cdot   \left(\left(t^2-1\right)\cdot \Lambda_{m}(t)\right)-T'''_{\bm{\lambda}}(m'+1, n'-1, t) \right|\nonumber \\
&\le\; 2^{-n}\label{ineqf1}
\end{align}
In this relation, by assumption,
\begin{equation}\label{eq2}
\deg\left( t\cdot   \left(t^2+1\right)\cdot \widetilde{R}_{\bm{\lambda}}\left(m'+1,n'-1;t^2\right)\right)\; =\; 2\cdot r_{\bm{\lambda}}\left(m'+1, n'-1\right)+1\;\le\; n-1.
\end{equation}

\noindent At the same time, the divisibility assumption~\eqref{eq1} yields a decomposition of the form $$ Q_{\bm{\lambda}}(m,h;t)\;=\; \left(t^2-1\right)\cdot \widetilde{Q}_{\bm{\lambda}}(m,h;t), \qquad \textrm{where}\qquad \deg\left(\widetilde{Q}_{\bm{\lambda}}(m,h;t)\right)\;=\;n-2.$$ From the equivalence~\eqref{convcharac} characterising the convergents of a Laurent series, it also holds that
\begin{equation}\label{ineqf2}
 \left|  \widetilde{Q}_{\bm{\lambda}}(m,h;t)  \cdot   \left(\left(t^2-1\right)\cdot \Lambda_{m}(t)\right)-P_{\bm{\lambda}}(m, h, t) \right|\;\le\; 2^{-n-1}\; =\; 2^{-(n-2)-3}.
\end{equation}

\noindent By the uniqueness of the convergent realising the normal order $n-2$ for the Laurent series $\left(t^2-1\right)\cdot \Lambda_{m}(t)$, one obtains from relations~\eqref{ineqf1} and~\eqref{ineqf2} that
\begin{equation}\label{eq3}
\frac{T'''_{\bm{\lambda}}(m'+1, n'-1, t)}{ t\cdot   \left(t^2+1\right)\cdot \widetilde{R}_{\bm{\lambda}}\left(m'+1,n'-1;t^2\right)}\;=\; \frac{P_{\bm{\lambda}}(m, h, t)}{ \widetilde{Q}_{\bm{\lambda}}(m,h;t) },
\end{equation}
where the fraction on the right--hand side is clearly irreducible. If $r_{\bm{\lambda}}(m'+1,n'-1)\le n'-2$, then it is  immediate  that the degree of the denominator of the left--hand side is at most $n-3$, which is strictly smaller than $\deg(\widetilde{Q}_{\bm{\lambda}}(m,h;t))$ ----  a contradiction. Therefore, $r_{\bm{\lambda}}(m'+1,n'-1)=n'-1$ and, given that $t^2+1$ is irreducible over $\F_\ell$ when $\ell\equiv 3\pmod{4}$, the comparison of the denominators in~\eqref{eq3} implies that $t$ divides the numerator $T'''_{\bm{\lambda}}(m'+1, n'-1, t)$. In other words, from the above expressions for this nume\-rator, one obtains working modulo $t$ that
\begin{align*}
0\;&\equiv\; T'''_{\bm{\lambda}}(m'+1, n'-1, t) \;\equiv\; T''_{\bm{\lambda}}(m'+1, n'-1, t)  +\lambda_m\cdot    R_{\bm{\lambda}}\left(m'+1,n'-1;t^2\right)\\
&\equiv\; T'_{\bm{\lambda}}\left(m'+1, n'-1, t^2\right)  +\lambda_m\cdot    R_{\bm{\lambda}}\left(m'+1,n'-1;t^2\right)\\
&\underset{\eqref{deft'lmnt}}{\equiv}\; T_{\bm{\lambda}}\left(m'+1, n'-1, t^2\right)  +\left(\lambda_{\left(m'+1\right)}+\lambda_m\right)\cdot    R_{\bm{\lambda}}\left(m'+1,n'-1;t^2\right)\\
&\underset{\eqref{reclambda}}{\equiv}\; T_{\bm{\lambda}}\left(m'+1, n'-1, t\right)  +\left(\lambda_{\left(m'+1\right)}+\lambda_{m'}\right)\cdot    R_{\bm{\lambda}}\left(m'+1,n'-1;t\right) \pmod{t}.
\end{align*}
Getting back to inequality~\eqref{x} and setting
\begin{align*}
 T_{\bm{\lambda}}\left(m'+1, n'-1, t\right)   &+\left(\lambda_{\left(m'+1\right)}+\lambda_{m'}\right)\cdot    R_{\bm{\lambda}}\left(m'+1,n'-1;t\right) \\
&\qquad  \qquad \qquad \qquad\qquad \qquad \qquad  =:\; t\cdot \widetilde{T}_{\bm{\lambda}}\left(m'+1, n'-1, t\right),
\end{align*}
one derives that
\begin{align*}
2^{-n'}\;&\ge\; \left|  R_{\bm{\lambda}}(m'+1,n'-1;t)\cdot \left(\left(t+1\right)\cdot \Lambda(t)-\lambda_1\right)_{\left(m'+1\right)}-T_{\bm{\lambda}}(m'+1, n'-1, t) \right|\\
&\underset{\eqref{shiftgammam}}{=}\; 2\cdot   \left|  R_{\bm{\lambda}}(m'+1,n'-1;t)\cdot \left(\left(t+1\right)\cdot \Lambda(t)-\lambda_1\right)_{m'}-\widetilde{T}_{\bm{\lambda}}(m'+1, n'-1, t) \right|.
\end{align*}
From Point~(3) in Lemma~\ref{closedromlem},  this implies that the entry with coordinates $\iota\left(m',n'-1\right)$ is followed by at least one zero entry in the Number Wall of the Laurent series \mbox{$\left(t+1\right)\cdot \Lambda(t)-\lambda_1$}; in other words, that $F_{\bm{\lambda}}\left(m', n'\right)=0$. This case has already been shown to lead one to a   contradiction in Case~A.1.

\paragraph{$\bullet$ Conclusion of Case~A~:} under the assumptions of this case, the polynomial $Q_{\bm{\lambda}}\left(m, h; t\right)$ cannot be divisible by $t^2-1$.\\

\subsection{Case B~: $(m,n)\equiv (1, 1) \pmod{2}$ }

Assume that $(m,n)=(2m'+1, 2n'-1)$ for some $m'\in\Z$ and $ n'\ge 2$. Reducing modulo $t^2-1$ the recurrence formula for $S_{\bm{\lambda}}(2m'+1,2n'-1;t) $ (see Proposition~\ref{proprec}) and equating it to $0$, one gets
\begin{equation}\label{fr=02}
F_{\bm{\lambda}}(m'+1,n'-1)\cdot R_{\bm{\lambda}}(m',n'-1;1)\;=\; 0.
\end{equation}


\paragraph{$\bullet$} If $F_{\bm{\lambda}}(m'+1,n'-1)=0$ then identity~(Z) in Corollary~\ref{labelhenkelbis} implies that $H_{\bm{\lambda}}\left(2m'+2, 2n'-1\right)=H_{\bm{\lambda}}\left(m+1, n\right)=0$. This contradicts Point~(a) in Lemma~\ref{lemlast}.
Therefore, the divisibility condition
\begin{equation}\label{divaprlemb}
t-1\mid R_{\bm{\lambda}}(m',n'-1;t).
\end{equation}
holds. This leads one to a   distinction of two cases~:

\paragraph{$\bullet$} \textbf{Case B.1}~: the polynomial $ R_{\bm{\lambda}}(m',n'-1;t)$ vanishes identically. From Point~(4) in Lemma~\ref{closedromlem}, this implies in particular that 
$F_{\bm{\lambda}}\left(m'+1,n'-1\right)=0$. However, this condition has already been ruled out.


\paragraph{$\bullet$} \textbf{Case B.2}~: the polynomial $ R_{\bm{\lambda}}(m',n'-1;t)$  does not  vanish and is divisible by $t-1$. Let then $$r_{\bm{\lambda}}\left(m', n'-1\right)\;:=\; \deg R_{\bm{\lambda}}(m',n'-1;t)\;\le\; n'-1$$ and $$R_{\bm{\lambda}}(m',n'-1;t)\;=\; \left(t-1\right)\cdot \widetilde{R}_{\bm{\lambda}}(m',n'-1;t),$$ where $\widetilde{R}_{\bm{\lambda}}(m',n'-1;t)$ is a polynomial of degree $r_{\bm{\lambda}}\left(m', n'-1\right)-1$. It follows from Point~(2) in Lemma~\ref{closedromlem} that for some polynomials $T_{\bm{\lambda}}(m', n'-1, t) $ coprime with $R_{\bm{\lambda}}(m',n'-1;t)$,
\begin{align}
&\left|  \left(t-1\right)\cdot \widetilde{R}_{\bm{\lambda}}(m',n'-1;t)\cdot \left(\left(t+1\right)\cdot \Lambda(t)-\lambda_1\right)_{m'}-T_{\bm{\lambda}}(m', n'-1, t) \right|\nonumber\\
&\qquad\qquad\qquad\underset{\eqref{parmlambdat+1}}{=}\;  \left|  \left(t^2-1\right)\cdot \widetilde{R}_{\bm{\lambda}}(m',n'-1;t)\cdot   \Lambda_{m'}(t)-T'_{\bm{\lambda}}(m', n'-1, t) \right|\nonumber\\
&\qquad\qquad\qquad \le\; 2^{-n'},\label{y}
\end{align}
where
\begin{align*}
T'_{\bm{\lambda}}( m', n'-1, t)\;=\; T_{\bm{\lambda}}(m', n'-1, t) +\lambda_{m'}\cdot \left(t-1\right)\cdot \widetilde{R}_{\bm{\lambda}}(m',n'-1;t).
\end{align*}


\noindent Apply the transformation $t\mapsto t^2$ to relation~\eqref{y}. From the functional equation~\eqref{funceqpair} satisfied by the paperfolding Laurent series $\Lambda(t)$, it yields the existence of  a polynomial $T''_{\bm{\lambda}}(m', n'-1, t)$ 
such that
\begin{align}
& \left|  \left(t^2-1\right)\cdot   \left(t^2+1\right)\cdot \widetilde{R}_{\bm{\lambda}}\left(m',n'-1;t^2\right)\cdot   \Lambda_{\left(2m'\right)}(t)-T''_{\bm{\lambda}}(m', n'-1, t) \right| \nonumber\\
&= \left|  \left(t^2+1\right)\cdot \widetilde{R}_{\bm{\lambda}}\left(m',n'-1;t^2\right)\cdot   \left(\left(t^2-1\right)\cdot \Lambda_{(m-1)}(t)\right)-T''_{\bm{\lambda}}(m', n'-1, t) \right|\nonumber \\
&\le\; 2^{-n-1}.\label{ineqf1b}
\end{align}
In this relation, by assumption, $$\deg\left(\left(t^2+1\right)\cdot \widetilde{R}_{\bm{\lambda}}\left(m'+1,n'-1;t^2\right)\right)\; =\; 2\cdot r_{\bm{\lambda}}\left(m'+1, n'-1\right)\;\le\; n-1.$$

\noindent Also, here again, the divisibility assumption~\eqref{eq1} yields a decomposition of the form $$ Q_{\bm{\lambda}}(m,h;t)\;=\; \left(t^2-1\right)\cdot \widetilde{Q}_{\bm{\lambda}}(m,h;t), \qquad \textrm{where}\qquad \deg\left(\widetilde{Q}_{\bm{\lambda}}(m,h;t)\right)\;=\;n-2.$$ From the equivalence~\eqref{convcharac} characterising the convergents of a Laurent series, it also holds that
\begin{align}\label{ineqf2b}
 2^{-n-1}\; &=\; 2^{-(n-1)-2}\nonumber \\
 &\ge \;  \left|  \widetilde{Q}_{\bm{\lambda}}(m,h;t)  \cdot   \left(\left(t^2-1\right)\cdot \Lambda_{m}(t)\right)-P_{\bm{\lambda}}(m, h, t) \right| \nonumber \\
 & \underset{\eqref{shiftgammam}}{=}  \left|  \left(t\cdot \widetilde{Q}_{\bm{\lambda}}(m,h;t)\right)  \cdot   \left(\left(t^2-1\right)\cdot \Lambda_{(m-1)}(t)\right)-\widetilde{P}_{\bm{\lambda}}(m, h, t) \right|,
\end{align}
where
\begin{align*}
\widetilde{P}_{\bm{\lambda}}(m, h, t)\;&=\; P_{\bm{\lambda}}(m, h, t)+ \lambda_{m-1}\cdot   \left(t^2-1\right)\cdot \widetilde{Q}_{\bm{\lambda}}(m,h;t)\\
&=\; P_{\bm{\lambda}}(m, h, t)+ \lambda_{m-1}\cdot   Q_{\bm{\lambda}}(m,h;t).
\end{align*}

\noindent By the uniqueness of the convergent realising the normal order $n-1$ for the Laurent series $\left(t^2-1\right)\cdot \Lambda_{(m-1)}(t)$, one obtains from relations~\eqref{ineqf1b} and~\eqref{ineqf2b} that
\begin{equation}\label{fraccasb}
\frac{T''_{\bm{\lambda}}(m', n'-1, t)}{   \left(t^2+1\right)\cdot \widetilde{R}_{\bm{\lambda}}\left(m',n'-1;t^2\right)}\;=\; \frac{\widetilde{P}_{\bm{\lambda}}(m, h, t)}{t\cdot  \widetilde{Q}_{\bm{\lambda}}(m,h;t) }\cdotp
\end{equation}
This equation induces a distinction of cases~:
\begin{itemize}
\item if  the numerator $\widetilde{P}_{\bm{\lambda}}(m, h, t)$ is a multiple of $t$, upon dividing inequality~\eqref{ineqf2b} by $t$, it follows from Point~(3) in Lemma~\ref{closedromlem} that the points with coordinates $\iota(m-1, n)$ and $\iota(m-1, n+1)$ all vanish in the Number Wall of the paperfolding sequence $\bm{\lambda}$. From the Square Window Lemma~\ref{sqwinthe}, this implies in particular that $H_{\bm{\lambda}}(m,n)=0$, contradicting Point~(a) in Lemma~\ref{lemlast};

\item if  the numerator $\widetilde{P}_{\bm{\lambda}}(m, h, t)$ is not divisible by $t$, then the fraction on the right-hand side of~\eqref{fraccasb} is irreducible and has a denominator of degree $n-1$. Since the fraction on the left-hand side has also denominator of degree at most $n-1$, one deduces that $t$, thence $t^2$, divides $\left(t^2+1\right)\cdot \widetilde{R}_{\bm{\lambda}}\left(m',n'-1;t^2\right)$. In turn, this implies that $t$ divides $ \widetilde{Q}_{\bm{\lambda}}(m,h;t)$, and therefore $Q_{\bm{\lambda}}(m,h;t)$, contradicting Point~(b) in Lemma~\ref{lemlast}.
\end{itemize}

\paragraph{$\bullet$ Conclusion of Case~B~:} under the assumptions of this case, the polynomial $Q_{\bm{\lambda}}\left(m, h; t\right)$ cannot be divisible by $t^2-1$.


\subsection{Case C~: $(m,n)\equiv(0, 0)\pmod{2}$}

Assume that $(m,n)=(2m', 2n')$ for some $m'\in\Z$ and $ n'\ge 1$. Reducing modulo $t^2-1$ the recurrence formula for  $S_{\bm{\lambda}}(2m',2n';t) $ (see Proposition~\ref{proprec}) and equating it to $0$, one gets
\begin{equation}\label{fr=03}
F_{\bm{\lambda}}(m',n')\cdot R_{\bm{\lambda}}(m',n'-1;1)\;=\; 0.
\end{equation}

\paragraph{$\bullet$} Suppose that $F_{\bm{\lambda}}(m',n')=0$. From identity~(Z) in Corollary~\ref{labelhenkelbis},  one has $H_{\bm{\lambda}}(2m',2n'+1) = H_{\bm{\lambda}}(m,n+1)=0$. From Point~(a) in Lemma~\ref{lemlast}, one thus infers the existence of a window $W$ in the Number Wall $NW_{\bm{\lambda}}$ with top-left corner located at $\iota(m,n)+(1,1)$. In other words,
\begin{equation}\label{qmhW}
Q_{\bm{\lambda}}(m,h;t) \; =\;  Q_{\bm{\lambda}}(W;t).
\end{equation}

\noindent Also, under the assumption that $F_{\bm{\lambda}}(m',n')=0$, Point~(3) in Lemma~\ref{closedromlem} guarantees the existence of a convergent $A(t)/B(t)$ of  the Laurent series $$((t+1)\cdot \Lambda(t)-\lambda_1)_{m'} \; =\;  ((t+1)\cdot \Lambda_{m'}(t) - \lambda_{m'})$$ such that 
\begin{equation}\label{inelasl}
\deg(B(t))\le n'-1\qquad \mbox{and}\qquad \left|B(t)\cdot ((t+1)\Lambda_{m'}(t)-\lambda_{m'}) - A(t)\right| \le 2^{-n'-1}.
\end{equation}
\noindent This induces a distinction of two further subcases~:
\begin{itemize}
\item if $t+1$ divides $A(t)+\lambda_{m'}\cdot B(t)$, then  divide the above inequality by $t+1$ to obtain that
    $$
    \deg(B(t))\le n'-1\qquad\mbox{and}\qquad \left|B(t)\cdot \Lambda_{m'}(t) - \widetilde{A}(t)\right|\le 2^{-n'-2},
    $$
where $\widetilde{A}(t) := (A(t) + \lambda_{m'}B(t))/(t+1)$ and where $\widetilde{A}(t)/B(t)$ is thus a convergent of the Laurent series $\Lambda_{m'}(t)$. This implies that $H_{\bm{\lambda}}(m',n') = H_{\bm{\lambda}}(m',n'+1) = 0$. In other words, the entry with coordinates $\iota(m',n')$ belongs to a window $W'$ of $NW_{\bm{\lambda}}$ of size at least 2, and therefore at least 3 from Corollary~\ref{odwind}. From Proposition~\ref{genwindow}, $W'$ then induces another window of size at least $3$ which is easily seen to contain the entry $2\cdot \iota(m', n')+(1,1) = \iota(m,n)+(1,1)$. This is saying that $W'$ induces the window $W$, which contradicts Point~(c) in Lemma~\ref{lemlast}.

\item if $t+1$ does not divide $A'(t):=A(t)+\lambda_{m'}\cdot B(t)$, then the fraction \mbox{$A'(t)/((t+1)\cdot B(t))$} is a convergent of $\Lambda_{m'}(t)$ such that $\deg((t+1)\cdot B(t))\le n'$. \\

\noindent If, furthermore, $\deg((t+1)\cdot B(t))=n',$ it  follows from Point~(3) in Lemma~\ref{closedromlem} and  from the second inequality in~\eqref{inelasl} that the entry with coordinates $\iota(m',n')$ lies in the top or left frame of a (possibly empty) window $W'$. From Proposition~\ref{convframe}, there exists an integer $s\ge 0$ such that $t^s\cdot Q_{\bm{\lambda}}(W';t)= (t+1)\cdot B(t)$. This shows in particular that  $t+1$ divides the denominator $Q_{\bm{\lambda}}(W';t)$, in which case Part~(a) of Proposition~\ref{genwindow} is applicable. It implies that $W'$ induces another window containing the entry with coordinates $2\cdot\iota(m',n')=\iota(m,n)$ in its top or left frame. From relation~\eqref{qmhW}, one deduces  that the induced window is $W$, which contradicts Point~(c) of Lemma~\ref{lemlast}. \\ 

\noindent If, however, $\deg((t+1)\cdot B(t))<n',$ then Point~(3) in Lemma~\ref{closedromlem} and  the second inequality in~\eqref{inelasl} guarantee that the above--defined window $W'$ is nonempty and contains the entry $\iota(m',n')$.  Case~(a) in Proposition~\ref{genwindow} is here also applicable, and it implies that $W'$ induces another window containing the entry with coordinates  $\iota(2m',2n')=\iota(m,n)$. This means that $H_{\bm{\lambda}}(m,n)=0$,  contradicting Point~(a) in Lemma~\ref{lemlast}.
\end{itemize}

%
%

\paragraph{$\bullet$} The  above distinction of cases shows  that
\begin{equation}\label{defm'n'pasdef}
F_{\bm{\lambda}}(m',n')\neq 0.
 \end{equation}
 Taking into account identity~\eqref{fr=03}, one then obtains the divisibility condition
\begin{equation}\label{divaprlemc}
t-1\mid R_{\bm{\lambda}}(m',n'-1;t),
\end{equation}
which leads one to a distinction of two cases~:

\paragraph{$\bullet$} \textbf{Case C.1}~:  the polynomial $ R_{\bm{\lambda}}(m',n'-1;t)$ vanishes identically. Point~(4) in Lemma~\ref{closedromlem} then  negates the just established relation~\eqref{defm'n'pasdef}. This case can thus not occur.


\paragraph{$\bullet$} \textbf{Case C.2}~: the polynomial $ R_{\bm{\lambda}}(m',n'-1;t)$  does not  vanish. 
Let then $$r_{\bm{\lambda}}\left(m', n'-1\right)\;:=\; \deg R_{\bm{\lambda}}(m',n'-1;t)\;\le\; n'-1$$ and $$R_{\bm{\lambda}}(m',n'-1;t)\;=\; \left(t-1\right)\cdot \widetilde{R}_{\bm{\lambda}}(m',n'-1;t),$$ where $\widetilde{R}_{\bm{\lambda}}(m',n'-1;t)$ is a polynomial of degree $r_{\bm{\lambda}}\left(m', n'-1\right)-1$. It follows from Point~(2) in Lemma~\ref{closedromlem} that for some polynomials $T_{\bm{\lambda}}(m', n'-1, t) $ coprime with $R_{\bm{\lambda}}(m',n'-1;t)$,
\begin{align}
&\left|  \left(t-1\right)\cdot \widetilde{R}_{\bm{\lambda}}(m',n'-1;t)\cdot \left(\left(t+1\right)\cdot \Lambda(t)-\lambda_1\right)_{m'}-T_{\bm{\lambda}}(m', n'-1, t) \right|\nonumber\\
&\qquad\qquad\qquad\underset{\eqref{parmlambdat+1}}{=}\;  \left|  \left(t^2-1\right)\cdot \widetilde{R}_{\bm{\lambda}}(m',n'-1;t)\cdot   \Lambda_{m'}(t)-T'_{\bm{\lambda}}(m', n'-1, t) \right|\nonumber\\
&\qquad\qquad\qquad \le\; 2^{-n'},\label{z}
\end{align}
where
\begin{equation}\label{deft'lmntc}
T'_{\bm{\lambda}}( m', n'-1, t)\;=\; T_{\bm{\lambda}}(m', n'-1, t) +\lambda_{m'}\cdot \left(t-1\right)\cdot \widetilde{R}_{\bm{\lambda}}(m',n'-1;t).
\end{equation}


\noindent Apply the transformation $t\mapsto t^2$ to relation~\eqref{z}. From the functional equation~\eqref{funceqpair} satisfied by the paperfolding Laurent series $\Lambda(t)$, it yields the existence of  a polynomial $T''_{\bm{\lambda}}(m', n'-1, t)$ 
such that
\begin{align}
& \left|  \left(t^2-1\right)\cdot   \left(t^2+1\right)\cdot \widetilde{R}_{\bm{\lambda}}\left(m',n'-1;t^2\right)\cdot   \Lambda_{\left(2m'\right)}(t)-T''_{\bm{\lambda}}(m', n'-1, t) \right| \nonumber\\
&= \left|  \left(t^2+1\right)\cdot \widetilde{R}_{\bm{\lambda}}\left(m',n'-1;t^2\right)\cdot   \left(\left(t^2-1\right)\cdot \Lambda_{m}(t)\right)-T''_{\bm{\lambda}}(m', n'-1, t) \right|\nonumber \\
&\le\; 2^{-n}.\label{ineqf1c}
\end{align}
In this relation, by assumption,
\begin{equation}\label{dertilC}
\deg\left(\left(t^2+1\right)\cdot \widetilde{R}_{\bm{\lambda}}\left(m'+1,n'-1;t^2\right)\right)\; =\; 2\cdot r_{\bm{\lambda}}\left(m'+1, n'-1\right)\;\le\; n-2.
\end{equation}
Furthermore,
\begin{align}\label{T'C}
T''_{\bm{\lambda}}(m', n'-1, t) =  T'_{\bm{\lambda}}\left(m', n'-1, t^2\right) +\left(-1\right)^{m'}\cdot t\cdot \left(t^2-1\right)\cdot \widetilde{R}_{\bm{\lambda}}\left(m',n'-1;t^2\right),
\end{align}
where $T'_{\bm{\lambda}}(m', n'-1, t)$ is defined in~\eqref{deft'lmntc}.  Also, the divisibility assumption~\eqref{eq1} yields here again a decomposition of the form $$ Q_{\bm{\lambda}}(m,h;t)\;=\; \left(t^2-1\right)\cdot \widetilde{Q}_{\bm{\lambda}}(m,h;t), \qquad \textrm{where}\qquad \deg\left(\widetilde{Q}_{\bm{\lambda}}(m,h;t)\right)\;=\;n-2.$$ From the equivalence~\eqref{convcharac} characterising the convergents of a Laurent series, it also holds that
\begin{align}\label{ineqf2c}
  \left|  \widetilde{Q}_{\bm{\lambda}}(m,h;t)  \cdot   \left(\left(t^2-1\right)\cdot \Lambda_{m}(t)\right)-P_{\bm{\lambda}}(m, h, t) \right|\;\le\;  2^{-n-1}\; &=\; 2^{-(n-2)-3}. 
\end{align}

\noindent By the uniqueness of the convergent realising the normal order $n-2$ for the Laurent series $\left(t^2-1\right)\cdot \Lambda_{m}(t)$, one obtains from relations~\eqref{ineqf1c} and~\eqref{ineqf2c} that
\begin{equation*}
\frac{T''_{\bm{\lambda}}(m', n'-1, t)}{\left(t^2+1\right)\cdot \widetilde{R}_{\bm{\lambda}}\left(m',n'-1;t^2\right)}\;=\; \frac{P_{\bm{\lambda}}(m, h, t)}{\widetilde{Q}_{\bm{\lambda}}(m,h;t) },
\end{equation*}
where the fraction of the right-hand side is clearly irreducible. Since from inequa\-lity~\eqref{dertilC}, the degree of the denominator on the left--hand side cannot be smaller than the one on the right--hand side, numerators and denominators can be identified. Thus, one obtains that $$P_{\bm{\lambda}}(m, h, t)\;=\; T''_{\bm{\lambda}}(m', n'-1, t)\quad \textrm{and}\quad Q_{\bm{\lambda}}(m,h;t)\;=\;  \left(t^2+1\right)\cdot R_{\bm{\lambda}}\left(m',n'-1;t^2\right), $$ where the second equation is obtained after multiplying out the denominators by \mbox{$t^2-1$}. Taking into account the right--hand side of inequality~\eqref{ineqf2c}, one can then rewrite inequality~\eqref{ineqf1c} as
\begin{align*}
2^{-n-1}\;&\ge\; \left|  \left(t^2+1\right)\cdot R_{\bm{\lambda}}\left(m',n'-1;t^2\right)\cdot  \Lambda_{m}(t)-T''_{\bm{\lambda}}(m', n'-1, t) \right| \\
&\underset{\eqref{funceqpair} \& \eqref{T'C}}{=}  \left|  \left(t^2+1\right)\cdot R_{\bm{\lambda}}\left(m',n'-1;t^2\right)\cdot  \Lambda_{m'}\left(t^2\right)-T'_{\bm{\lambda}}(m', n'-1, t^2) \right|.
\end{align*}
Since this last quantity is the norm of an even series and that the integer $n$ is assumed to be even, the above inequality actually implies the following stronger one~:
\begin{align*}
 \left|  \left(t^2+1\right)\cdot R_{\bm{\lambda}}\left(m',n'-1;t^2\right)\cdot  \Lambda_{m'}\left(t^2\right)-T'_{\bm{\lambda}}(m', n'-1, t^2) \right| \;\le\; 2^{-n-2}.
\end{align*}
In turn, this implies that
\begin{align}\label{inecaccon}
\left|  \left(t+1\right)\cdot R_{\bm{\lambda}}\left(m',n'-1;t\right)\cdot  \Lambda_{m'}\left(t\right)-T'_{\bm{\lambda}}(m', n'-1, t) \right| \;\le\; 2^{-n'-1},
\end{align}
which induces a distinction of cases~:
\begin{itemize}
\item if $t+1$ does not divide $T'_{\bm{\lambda}}(m', n'-1, t)$, it follows from the equivalence stated in~\eqref{convcharac} that the polynomial  $\left(t+1\right)\cdot R_{\bm{\lambda}}\left(m',n'-1;t\right)$ of degree at most $n'$ is the deno\-mi\-nator of a convergent to the Laurent series $ \Lambda_{m'}\left(t\right)$. Let $W'$ be the window the point with coordinates $\iota(m',n'+1)$ belongs to.
    From Part~(a) of Proposition~\ref{genwindow}, $W'$ induces a window, say $W$, which is easily seen to contain the point with coordinates $\iota(2m',2n')+(1,1)$. Since from Point~(b) in Lemma~\ref{lemlast},   $t$ does not divide $Q_{\bm{\lambda}}(m,h;t)$, it follows from Proposition~\ref{convframe} that $Q_{\bm{\lambda}}(m,h;t) = Q_{\bm{\lambda}}(W;t)$. This contradicts Point~(c) in Lemma~\ref{lemlast}.


\item if $T'_{\bm{\lambda}}(m', n'-1, t)$ is a multiple of  $t+1$, divide inequality~\eqref{inecaccon} by $t+1$.  From Point~(3) in  Lemma~\ref{closedromlem}, this  guarantees that the point with coordinates $\iota(m', n')$ belongs to a window in the Number Wall of the paper\-fol\-ding sequence $\bm{\lambda}$, say $W'$, which has size  $\rho(W')\ge 2$. Since all windows  have uneven size from Corollary~\ref{odwind}, this actually means that  $\rho(W')\ge 3$. From Proposition~\ref{genwindow}, $W'$ then induces a window, say $W$, with size $\rho\left(W\right)\ge 3$ which clearly contains the point $\iota(2m,',2n') + (1,1)$. Here again, this implies that $Q_{\bm{\lambda}}(m,h;t) = Q_{\bm{\lambda}}(W;t)$, which contradicts Point~(c) in Lemma~\ref{lemlast}.

\end{itemize}

\paragraph{$\bullet$ Conclusion of Case~C~:} under the assumptions of this case, the polynomial $Q_{\bm{\lambda}}\left(m, h; t\right)$ cannot be divisible by $t^2-1$.\\


\subsection{Case D~: $(m,n)\equiv(0, 1)\pmod{2}$}

Assume that $(m,n)=(2m', 2n'-1)$ for some $m'\in\Z$ and $ n'\ge 2$. Reducing modulo $t^2-1$  the recurrence formula for $S_{\bm{\lambda}}(2m',2n'-1;t) $ (see Proposition~\ref{proprec}) and equating it to $0$, one gets
\begin{equation}\label{fr=04}
F_{\bm{\lambda}}(m',n'-1)\cdot R_{\bm{\lambda}}(m',n'-1;1)\;=\; 0.
\end{equation}



\paragraph{$\bullet$} Suppose that $F_{\bm{\lambda}}(m',n'-1)=0$. Identity~(Z) in Corollary~\ref{labelhenkelbis} then implies that $H_{\bm{\lambda}}\left(2m', 2n'-1\right)=H_{\bm{\lambda}}\left(m, n\right)=0$. This contradicts Point~(a) in Lemma~\ref{lemlast}.
Therefore,
\begin{equation}\label{neqcondiF}
F_{\bm{\lambda}}(m',n'-1)\;\neq\; 0
 \end{equation}
and the divisibility condition
\begin{equation}\label{divaprlemd}
t-1\mid R_{\bm{\lambda}}(m',n'-1;t)
\end{equation}
is met. This leads to a  distinction of two cases~:

\paragraph{$\bullet$} \textbf{Case D.1}~: the polynomial $ R_{\bm{\lambda}}(m',n'-1;t)$ vanishes identically. From Point~(1) in Lemma~\ref{closedromlem}, this implies in particular that $F_{\bm{\lambda}}\left(m',n'-1\right)=0$, which is  ruled out by~\eqref{neqcondiF}.

\paragraph{$\bullet$} \textbf{Case D.2}~: the polynomial $ R_{\bm{\lambda}}(m',n'-1;t)$  does not  vanish and is divisible by $t-1$. Let in this final case also $$r_{\bm{\lambda}}\left(m', n'-1\right)\;:=\; \deg R_{\bm{\lambda}}(m',n'-1;t)\;\le\; n'-1$$ and $$R_{\bm{\lambda}}(m',n'-1;t)\;=\; \left(t-1\right)\cdot \widetilde{R}_{\bm{\lambda}}(m',n'-1;t),$$ where $\widetilde{R}_{\bm{\lambda}}(m',n'-1;t)$ is a polynomial of degree $r_{\bm{\lambda}}\left(m', n'-1\right)-1$. It follows from Point~(2) in Lemma~\ref{closedromlem} that for some polynomials $T_{\bm{\lambda}}(m', n'-1, t) $ coprime with $R_{\bm{\lambda}}(m',n'-1;t)$,
\begin{align}
&\left|  \left(t-1\right)\cdot \widetilde{R}_{\bm{\lambda}}(m',n'-1;t)\cdot \left(\left(t+1\right)\cdot \Lambda(t)-\lambda_1\right)_{m'}-T_{\bm{\lambda}}(m', n'-1, t) \right|\nonumber\\
&\qquad\qquad\qquad\underset{\eqref{parmlambdat+1}}{=}\;  \left|  \left(t^2-1\right)\cdot \widetilde{R}_{\bm{\lambda}}(m',n'-1;t)\cdot   \Lambda_{m'}(t)-T'_{\bm{\lambda}}(m', n'-1, t) \right|\nonumber\\
&\qquad\qquad\qquad \le\; 2^{-n'},\label{t}
\end{align}
where
\begin{align}\label{deft'lmntd}
T'_{\bm{\lambda}}( m', n'-1, t)\;=\; T_{\bm{\lambda}}(m', n'-1, t) +\lambda_{m'}\cdot \left(t-1\right)\cdot \widetilde{R}_{\bm{\lambda}}(m',n'-1;t).
\end{align}

%
\noindent Apply the transformation $t\mapsto t^2$ to relation~\eqref{t}. From the functional equation~\eqref{funceqpair} satisfied by the paperfolding Laurent series $\Lambda(t)$, it yields the existence of  a polynomial $T''_{\bm{\lambda}}(m', n'-1, t)$ 
such that
\begin{align}
& \left|  \left(t^2-1\right)\cdot   \left(t^2+1\right)\cdot \widetilde{R}_{\bm{\lambda}}\left(m',n'-1;t^2\right)\cdot   \Lambda_{\left(2m'\right)}(t)-T''_{\bm{\lambda}}(m', n'-1, t) \right| \nonumber\\
&= \left|  \left(t^2+1\right)\cdot \widetilde{R}_{\bm{\lambda}}\left(m',n'-1;t^2\right)\cdot   \left(\left(t^2-1\right)\cdot \Lambda_{m}(t)\right)-T''_{\bm{\lambda}}(m', n'-1, t) \right|\nonumber \\
&\le\; 2^{-n-1}.\label{ineqf1d}
\end{align}
Here, by assumption,
\begin{equation}\label{degdenomD}
\deg\left(\left(t^2+1\right)\cdot \widetilde{R}_{\bm{\lambda}}\left(m'+1,n'-1;t^2\right)\right)\; =\; 2\cdot r_{\bm{\lambda}}\left(m'+1, n'-1\right)\;\le\; n-1.
\end{equation}
Furthermore,
\begin{align}\label{T''last}
T''_{\bm{\lambda}}(m', n'-1, t) =  T'_{\bm{\lambda}}\left(m', n'-1, t^2\right) +\left(-1\right)^{m'}\cdot t\cdot \left(t^2-1\right)\cdot \widetilde{R}_{\bm{\lambda}}\left(m',n'-1;t^2\right),
\end{align}
where $T'_{\bm{\lambda}}(m', n'-1, t)$ is defined in~\eqref{deft'lmntd}. Also, recall one last time that the divisibility assumption~\eqref{eq1} yields a decomposition of the form $$ Q_{\bm{\lambda}}(m,h;t)\;=\; \left(t^2-1\right)\cdot \widetilde{Q}_{\bm{\lambda}}(m,h;t), \qquad \textrm{where}\qquad \deg\left(\widetilde{Q}_{\bm{\lambda}}(m,h;t)\right)\;=\;n-2.$$ From the equivalence~\eqref{convcharac} characterising the convergents of a Laurent series, it also holds that
\begin{align}\label{ineqf2d}
 2^{-n-1}\; &=\; 2^{-(n-2)-3}\nonumber \\
 &\ge \;  \left|  \widetilde{Q}_{\bm{\lambda}}(m,h;t)  \cdot   \left(\left(t^2-1\right)\cdot \Lambda_{m}(t)\right)-P_{\bm{\lambda}}(m, h, t) \right|.
\end{align}

\noindent By the uniqueness of the convergent realising the normal order $n-2$ for the Laurent series $\left(t^2-1\right)\cdot \Lambda_{m}(t)$, one obtains from relations~\eqref{ineqf1d} and~\eqref{ineqf2d} that
\begin{equation}\label{fraccasd}
\frac{T''_{\bm{\lambda}}(m', n'-1, t)}{   \left(t^2+1\right)\cdot \widetilde{R}_{\bm{\lambda}}\left(m',n'-1;t^2\right)}\;=\; \frac{P_{\bm{\lambda}}(m, h, t)}{\widetilde{Q}_{\bm{\lambda}}(m,h;t) },
\end{equation}
where the right--hand side is an irreducible fraction with denominator of degree \mbox{$n-2\equiv 1\pmod{2}$}. As for the left--hand side, the explicit expression for the nu\-me\-rator stated in~\eqref{deft'lmntd} and in~\eqref{T''last} makes it clear that the greatest common divisor of the numerator and the denominator is either $1$ or $1+t^2$ (which is an irreducible polynomial over $\F_\ell$ when $\ell \equiv 3\pmod{4}$). Then, the denominator of the irreducible form of the left--hand side has an even degree and can thus not coincide with the polynomial $\widetilde{Q}_{\bm{\lambda}}(m,h;t)$. This leads one to a contradiction for one last time.

\paragraph{$\bullet$ Conclusion of Case~D~:} under the assumptions of this case, the polynomial $Q_{\bm{\lambda}}\left(m, h; t\right)$ cannot be divisible by $t^2-1$.\\

\paragraph{$\bullet$ General Conclusion~:} under the assumptions of  each of the cases (A), (B), (C) and (D),  the polynomial $Q_{\bm{\lambda}}\left(m, h; t\right)$ cannot be divisible by $t^2-1$. Since the cases (A)--(D) cover all possible configurations of parity for the pair $(m,n)$, this shows that a polynomial $Q_{\bm{\lambda}}\left(m, h; t\right)$ meeting the divisibility assumption~\eqref{eq1} cannot exist. This establishes the Main Theorem~\ref{mainthm}.\\
\end{proof}

\paragraph{Declaration.}  The work of the first--named author was supported by the French ANR (project "GoDiBAp", grant number ANR-25-CE40-1961-01). The authors have no competing interests to declare that are relevant to the content of this article.

\bibliographystyle{unsrt}


\newpage

\section*{\centerline{Appendix~: Some Representations of Number Walls}}
\addcontentsline{toc}{section}{Appendix~: Some Representations of Number Walls}

\vspace{10mm}

\noindent A portion of the Number Wall of the doubly infinite paperfolding sequence over the alphabet $(-1, 1)$ (see~\eqref{reclambda} for the definition) is represented below over the first four finite fields with order congruent to 3 modulo 4, viz.~$\F_3$, $\F_7$, $\F_{11}$ and $\F_{19}$.\\

\noindent In the notations of Section~\ref{sec2}, the portion corresponds in each case to the range of indices $-100\le k\le 100$ (for the columns, from left to right) and $1\le j\le 100$ (for the rows, from top to bottom).  All pictures use the same colour scheme~: zero entries are in  red and   the nonzero values are coloured in order darkest to lightest.\\

\noindent The pictures confirm the results proved in Section~\ref{sec4}, namely that zero entries are either isolated or that they form a square of size 3. The process of generation of these squares is elucidated in Proposition~\ref{genwindow}.  \\

\vspace{10mm}

\begin{figure}[H]
         \centering
         \includegraphics[width=0.975\textwidth]{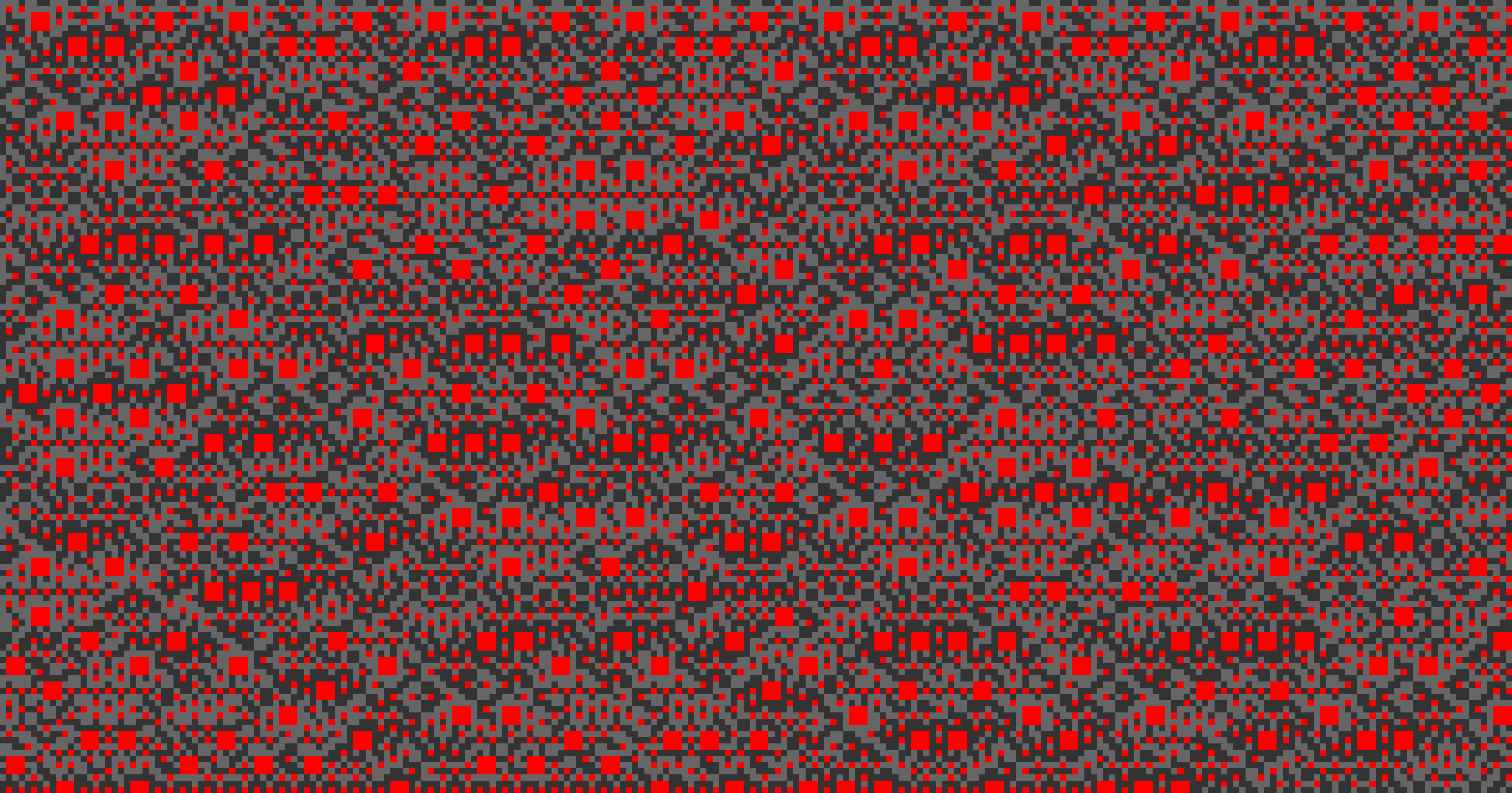}
         \caption{A portion of the Number Wall over $\F_3$}
\end{figure}

\begin{figure}[H]
         \centering
         \includegraphics[width=\textwidth]{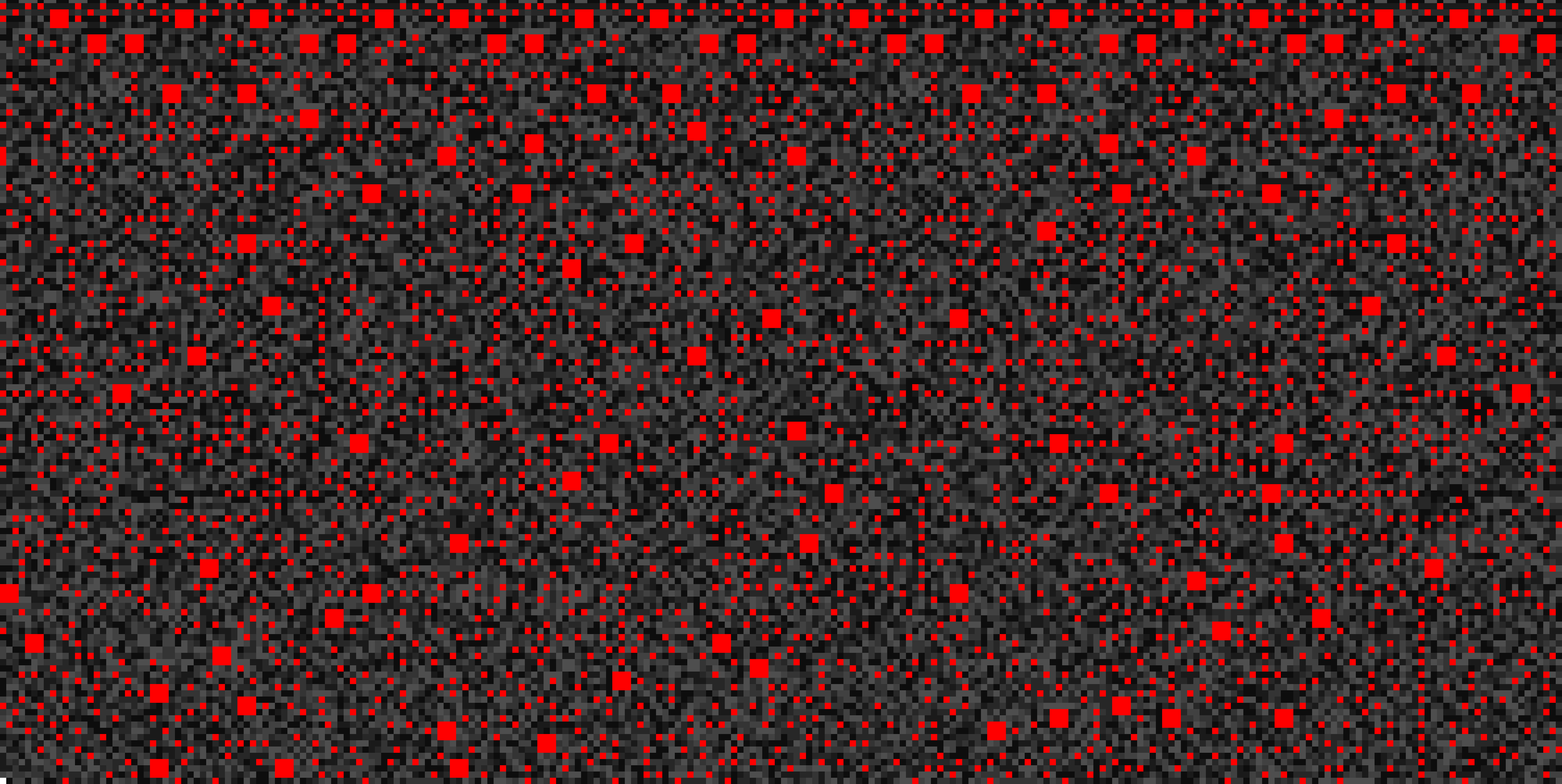}
         \caption{A portion of the Number Wall over $\F_7$}
\end{figure}

\begin{figure}[H]
         \centering
         \includegraphics[width=\textwidth]{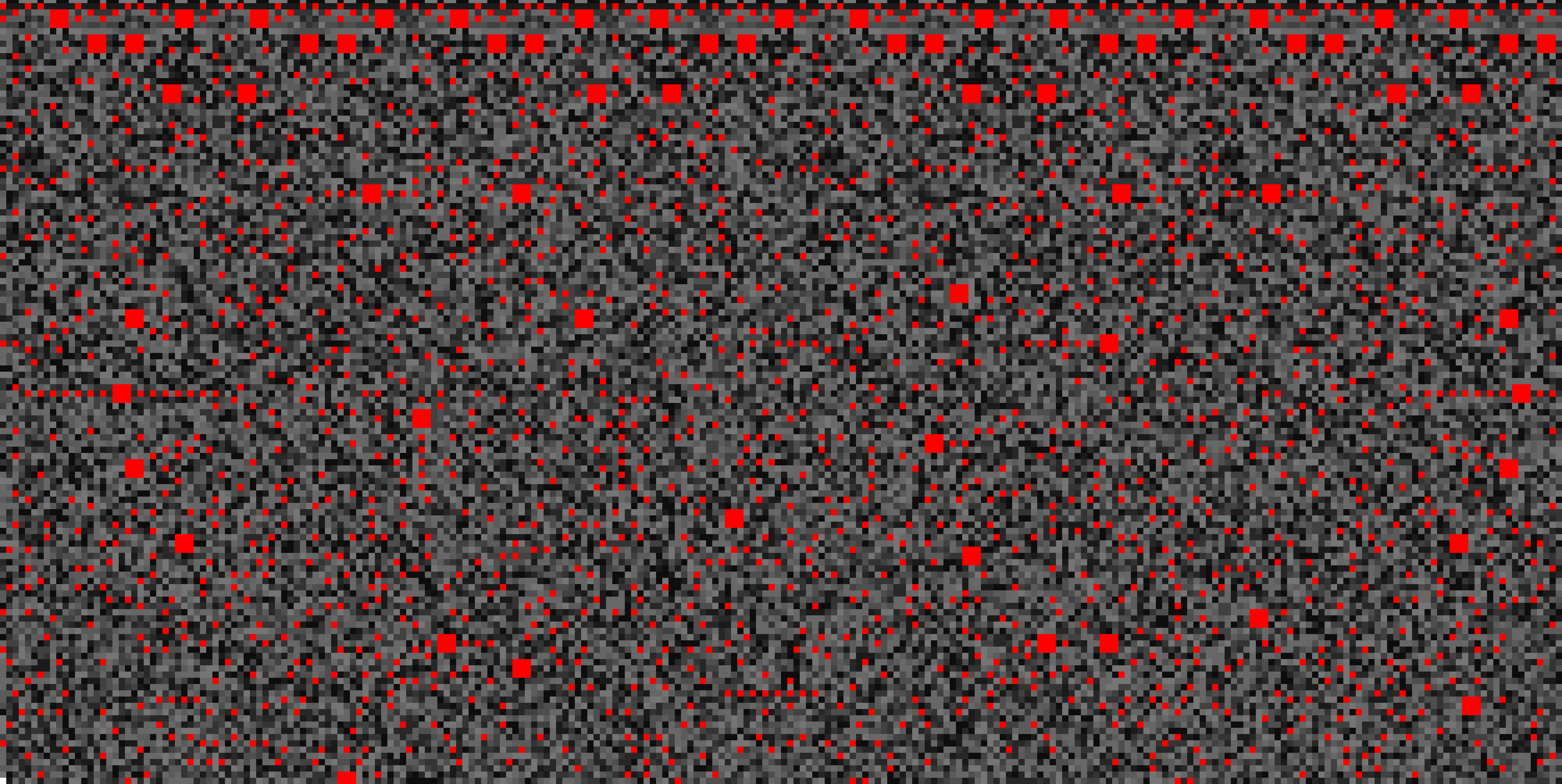}
         \caption{A portion of the Number Wall over $\F_{11}$}
\end{figure}

\begin{figure}[H]
         \centering
         \includegraphics[width=\textwidth]{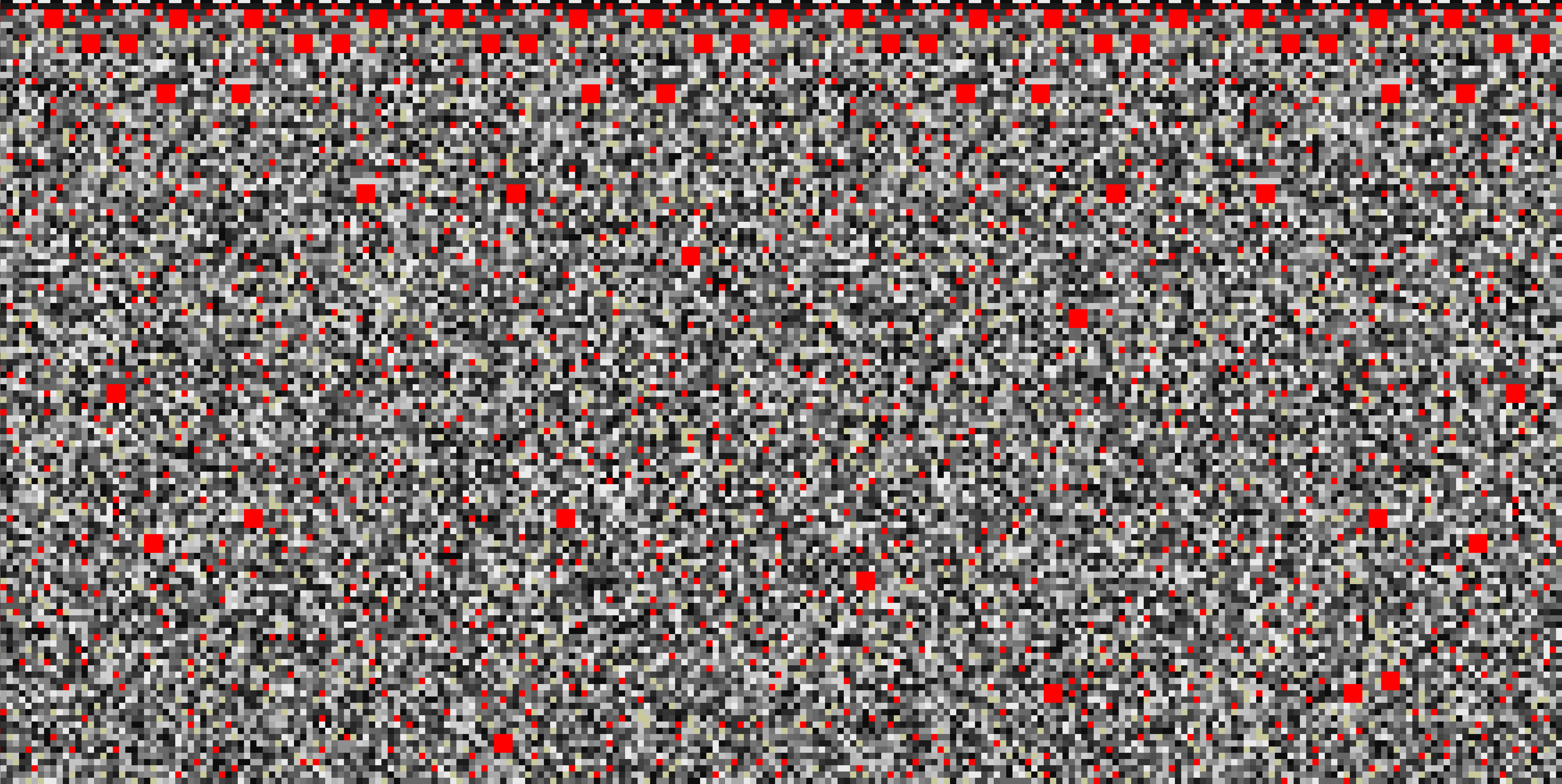}
         \caption{A portion of the Number Wall over $\F_{19}$}
\end{figure}

\paragraph{}\noindent All  pictures displayed in this appendix were generated by Steven Robertson.

\end{document}